%% file: KRZ23SW.tex
\documentclass[11pt,oneside,reqno]{amsart}

\usepackage[a4paper, total={450pt,675pt}]{geometry}
\usepackage[OT2,T1]{fontenc}
\usepackage[utf8]{inputenc}
\usepackage{amssymb,bm}

\usepackage[dvipsnames]{xcolor}
\usepackage[
    colorlinks=true,
    linkcolor=Maroon,
    citecolor=JungleGreen,
    urlcolor=NavyBlue]{hyperref}

\usepackage{xcolor,colortbl}
\usepackage[capitalize,nameinlink]{cleveref}

\usepackage{enumitem}
\usepackage{dsfont}
\usepackage{tikz}
\usetikzlibrary{patterns}
\usepackage{subcaption}
\usetikzlibrary{arrows}
\usepackage{multirow}

\usepackage{subcaption}
\newtheorem{theorem}{Theorem}[section]

\newtheorem{lemma}[theorem]{Lemma}
\newtheorem{proposition}[theorem]{Proposition}
\newtheorem{corollary}[theorem]{Corollary}

\newtheorem{remark}[theorem]{Remark}

\crefname{section}{Sect.}{section}
\numberwithin{equation}{section}

\DeclareMathOperator{\re}{Re}

\DeclareMathOperator{\SL}{SL}

\DeclareMathOperator{\supp}{supp}

\DeclareMathOperator{\const}{const.}

\newcommand*\diff{\mathop{}\!\mathrm{d}}

\makeatletter
\def\@tocline#1#2#3#4#5#6#7{\relax
  \ifnum #1>\c@tocdepth 
  \else
    \par \addpenalty\@secpenalty\addvspace{#2}%
    \begingroup \hyphenpenalty\@M
    \@ifempty{#4}{%
      \@tempdima\csname r@tocindent\number#1\endcsname\relax
    }{%
      \@tempdima#4\relax
    }%
    \parindent\z@ \leftskip#3\relax
    \advance\leftskip\@tempdima\relax
    \rightskip\@pnumwidth plus4em \parfillskip-\@pnumwidth
    #5\leavevmode\hskip-\@tempdima
      \ifcase #1
       \or\or \hskip 2em \or \hskip 2em \else \hskip 3em \fi%
      #6\nobreak\relax
    \dotfill\hbox to\@pnumwidth{\@tocpagenum{#7}}\par
    \nobreak
    \endgroup
  \fi}
\makeatother

\usepackage[
    backend=biber,
    style=alphabetic,
    giveninits=true,
    maxbibnames=10,
    sorting=nyt]{biblatex}

\DeclareFieldFormat[article,incollection]{title}{#1}
\DeclareFieldFormat{date}{#1}
\DeclareFieldFormat{pages}{{#1}}
\DeclareFieldFormat{doi}{
\href{https://www.doi.org/#1}{Doi}}
\DeclareFieldFormat{url}{
\href{#1}{Doi}}

\renewbibmacro{in:}{%
  \ifentrytype{article}{}{\printtext{\bibstring{in}\intitlepunct}}}

\begin{document}

\title[Stein-Weiss inequality on non-compact symmetric spaces]{
Stein-Weiss inequality on non-compact symmetric spaces}

\author{Vishvesh Kumar, Michael Ruzhansky, Hong-Wei Zhang}

\begin{abstract}
Let $\Delta$ be the Laplace-Beltrami operator on a non-compact symmetric space of any rank, and denote the bottom of its $L^2$-spectrum as $-|\rho|^{2}$. In this paper, we provide a comprehensive characterization of both the sufficient and necessary conditions ensuring the validity of the Stein-Weiss inequality for the entire family of operators $\lbrace{(-\Delta+b)^{-\frac{\sigma}{2}}}\rbrace_{\sigma\ge0,\,b\ge-|\rho|^{2}}$. As an application, some weighted functional inequalities, such as Heisenberg's uncertainty principle, Gagliardo-Nirenberg's interpolation inequality, Pitt's inequality, etc., become available in this context. In particular, their sets of admissible indices are larger than those in the Euclidean setting.
\end{abstract}


\makeatletter
\@namedef{subjclassname@2020}{\textnormal{2020}
    \it{Mathematics Subject Classification}}
\makeatother

\maketitle
\setcounter{tocdepth}{2}
\tableofcontents

\section{Introduction}

\subsection{Background and motivation}
In the 1950s, E. M. Stein and G. Weiss proved that for any $N\ge1$, $0<\sigma<N$, $1<p\le{q}<+\infty$, $\alpha<\frac{N}{p'}$, $\beta<\frac{N}{q}$, $\alpha+\beta\ge0$, and $\frac{1}{p}-\frac{1}{q}=\frac{N-\sigma-\alpha-\beta}{N}$, there exists a constant $C(p,q,\alpha,\beta)>0$ such that
\begin{align}\label{SW on RN}
    \int_{\mathbb{R}^{N}}\diff{x}\,
    \int_{\mathbb{R}^{N}}\diff{y}\,
    \frac{f(y)\,g(x)}{|x|^{\beta}\,|x-y|^{\sigma}\,|y|^{\alpha}}
    \le\,C(p,q,\alpha,\beta)\,
    \|f\|_{L^{p}(\mathbb{R}^{N})}\,
    \|g\|_{L^{q'}(\mathbb{R}^{N})},
\end{align}
for all $f\in{L^{p}(\mathbb{R}^{N})}$ and $g\in{L^{q'}(\mathbb{R}^{N})}$. Here, $p'$ and $q'$ are conjugate exponents of $p$ and $q$ in the sense of Hölder. This result extends the one obtained by Hardy and Littlewood in dimension $N=1$ and generalizes the Hardy-Littlewood-Sobolev inequality by incorporating two radial weights. See \cite{HL28,Sob38,SW58}. Over the last 70 years, the Stein-Weiss inequality has been proven to be a powerful tool in connecting different functional inequalities and solving nonlinear partial differential equations, which are both essential topics in the study of analysis on manifolds. This paper contributes to the literature by providing a comprehensive description of the Stein-Weiss inequality in general \textit{non-compact Riemannian symmetric spaces} (which we will refer to as \textit{symmetric spaces} in the following).  It establishes both the \textit{sufficient} and \textit{necessary} conditions for the validity of the Stein-Weiss inequality, including the two endpoints $p=1$ and $q=+\infty$. In addition, we derive other weighted functional inequalities with larger sets of admissible indices and study how the particular geometric structure of the symmetric space affects the behavior of these inequalities.

Researchers have explored different variations of the Hardy-Littlewood-Sobolev inequality and the Stein-Weiss inequality in the Euclidean setting to solve a wide range of nonlinear partial differential equations problems. These variants may involve different function spaces, weights, or specific assumptions on the underlying PDE model. We refer to \cite{Str69,Wal71,CF74,MW74,Per90,TVZ07,Bec08b,,Wan15,CLLT18,CK19,CLLT20,CC22} and the references therein. See also \cite{DZ15,Dou16}, \cite{SW21,Wan21}, and \cite{GIT21} for the Stein-Weiss inequality in the upper plane, the product spaces, and the Dunkl setting. In the non-Euclidean counterpart, the Stein-Weiss inequality has been well-established in certain Lie groups that possess the dilation property. See \cite{FS74,SW92,Bec97,GMS10,FL12,HLZ12,Bec15,KRS19,CLT19,RY19,Bec21} on the Heisenberg group, the Carnot groups, and the homogeneous groups.

Let $\mathbb{X}$ be a symmetric space of non-compact type, which is a Riemannian manifold with non-positive sectional curvature. Keep in mind that since $\mathbb{X}$ grows exponentially fast at infinity, the rescaling methods usually used to handle the homogeneous weights fail in this setting. We denote by $-\Delta$ the positive Laplace-Beltrami operator on $\mathbb{X}$ and $|\rho|^2$ the bottom of its $L^2$-spectrum. In the present paper, we focus on a large family of operators 
\begin{align}
    R(\zeta,\sigma)\,
    =\,
    (-\Delta-|\rho|^2+\zeta^2)^{-\frac{\sigma}{2}}
    \qquad\forall\,\zeta\ge0,\,\,\,\forall\,\sigma\ge0,
\end{align}
which includes several prominent examples such as the resolvent $R(\zeta,2)$, the Riesz potential $R(|\rho|,\sigma)$, the Green operator $R(|\rho|,2)$, etc. Let $k_{\zeta,\sigma}$ be its corresponding convolution kernel. The goal of this paper is to establish the following Stein-Weiss inequality
\begin{align}\label{SW on X}
    \int_{\mathbb{X}}\diff{x}\,
    \int_{\mathbb{X}}\diff{y}\,
    |x|^{-\beta}\,g(x)\,k_{\zeta,\sigma}(y^{-1}x)\,f(y)\,|y|^{-\alpha}
    \le\,C\,
    \|f\|_{L^{p}(\mathbb{X})}\,
    \|g\|_{L^{q'}(\mathbb{X})},
\end{align}
by determining all possible values for the components $\alpha$, $\beta\in\mathbb{R}$ and indices $1\le{p,\,q}\le+\infty$. The constant $C$ depends on $p,\,q,\,\alpha$, and $\beta$, but not on $f\in{L^{p}(\mathbb{X})}$ or $g\in{L^{q'}(\mathbb{X})}$.

\subsection{Statement of main results}
We omit the trivial case where $\sigma=0$. For all $\zeta\ge0$, $\sigma>0$, $\alpha\in\mathbb{R}$, and $\beta\in\mathbb{R}$, we define the operator
\begin{align*}
    T(\zeta,\sigma,\alpha,\beta)f(x)\,=\,
    \int_{\mathbb{X}}\diff{y}\,
    |x|^{-\beta}\,k_{\zeta,\sigma}(y^{-1}x)\,|y|^{-\alpha}\,f(y),
\end{align*}
where $f$ is a suitable function on $\mathbb{X}$. Note that the left side of \eqref{SW on X} is $\langle{T(\zeta,\sigma,\alpha,\beta)f,g}\rangle$ and is equal to $\langle{f,T(\zeta,\sigma,\beta,\alpha)g}\rangle$. Then, by duality of Lebesgue spaces, the Stein-Weiss inequality \eqref{SW on X} holds if and only if the operator $T(\zeta,\sigma,\alpha,\beta)$ is bounded from $L^{p}(\mathbb{X})$ to $L^{q}(\mathbb{X})$. Before presenting our main results, it is important to highlight that the behavior of the convolution kernel $k_{\zeta,\sigma}$ differs depending on its distance from the origin, as well as the values of $\sigma$ and $\zeta$. Let $\chi_{0}\in\mathcal{C}_{c}^{\infty}(\mathbb{R}_{+})$ be a cut-off function such that its support is contained in the interval $[0,1]$ and $\chi_{0}=1$ on the interval $[0,\frac12]$. We define $\chi_{\infty}=1-\chi_{0}$. For any $x\in\mathbb{X}$, we denote $\psi_{0}(x)$ as $\chi_{0}(|x|)$ (resp. $\psi_{\infty}(x)$ as $\chi_{\infty}(|x|)$), representing the bi-$K$-invariant cut-off functions defined on $\mathbb{X}$. According to \cite[Theorem 4.2.2]{AJ99}, we know that \footnote{Throughout this paper, the notation $A\lesssim{B}$ between two positive expressions means that $A\le{CB}$ for some constants $C>0$, and $A\asymp{B}$ means $A\lesssim{B}\lesssim{A}$.}
\begin{align}
    (\psi_{0}k_{\zeta,\sigma})(x)\,
    \asymp\,
    \begin{cases}
        |x|^{\sigma-n}
        &\qquad\textnormal{if}\,\,\,0<\sigma<n,\\[5pt]
        \log (\frac{1}{|x|})
        &\qquad\textnormal{if}\,\,\,\sigma=n,\\[5pt]
        1
        &\qquad\textnormal{if}\,\,\,\sigma>n,
    \end{cases}
    \label{estim k0}
\end{align}
for all $\zeta\ge0$, and 
\begin{align}
    (\psi_{\infty}k_{\zeta,\sigma})(x)\,
    \asymp\,
    \varphi_{0}(x)\,e^{-\zeta|x|}\,
    \begin{cases}
        |x|^{\frac{\sigma-\nu-1}{2}}
        &\qquad\textnormal{if}\,\,\,\zeta>0
        \,\,\,\textnormal{and}\,\,\,\sigma>0,\\[5pt]
        |x|^{\sigma-\nu}
        &\qquad\textnormal{if}\,\,\,\zeta=0
        \,\,\,\textnormal{and}\,\,\,0<\sigma<\nu,
    \end{cases}
    \label{estim kinf}
\end{align}
where $n\ge1$ and $\nu\ge3$ denote the topological dimension and pseudo-dimension of $\mathbb{X}$, and $\varphi_{0}$ is the ground spherical function which satisfies \eqref{phi0}. See the next section for more details on these notations. Notice that the jump in the polynomial factor in \eqref{estim kinf} between the cases $\zeta=0$ and $\zeta>0$, can also be observed in our Stein-Weiss inequality involving some critical indices. 

Due to the different behaviors \eqref{estim k0} and \eqref{estim kinf} of the convolution kernel, we divide the operator $T(\zeta,\sigma,\alpha,\beta)$ into two parts:
\begin{align*}
    &T(\zeta,\sigma,\alpha,\beta)f(x)\\[5pt]
    &=\,
    \underbrace{\int_{\mathbb{X}}\diff{y}\,
    |x|^{-\beta}\,(\psi_{0}k_{\zeta,\sigma})(y^{-1}x)\,|y|^{-\alpha}\,f(y)
    }_{=\,T^{0}(\zeta,\sigma,\alpha,\beta)f(x)}\,
    +\,
    \underbrace{
    \int_{\mathbb{X}}\diff{y}\,
    |x|^{-\beta}\,(\psi_{\infty}k_{\zeta,\sigma})(y^{-1}x)\,|y|^{-\alpha}\,f(y)
    }_{=\,T^{\infty}(\zeta,\sigma,\alpha,\beta)f(x)}.
\end{align*}
The following theorem clarifies both the sufficient and necessary conditions for which the operators $T^{0}$ and $T^{\infty}$ are $L^p$-$L^q$-bounded. \footnote{For the sake of simplicity, we use $T^{0}$ and $T^{\infty}$ when the values of $\zeta$, $\sigma$, $\alpha$, and $\beta$ do not need to be specified in the context. An operator is called $L^p$-$L^q$-bounded if it is defined as a bounded operator from $L^p$ to $L^q$. Furthermore, an operator is called  $L^p$-bounded if it is $L^p$-$L^p$-bounded.}

\begin{theorem}\label{mainthm1}
    Let $\mathbb{X}$ be a non-compact symmetric space of rank $\ell\ge1$, dimension $n\ge2$, and pseudo-dimension $\nu\ge3$. Suppose that $\zeta\ge0$, $\sigma>0$ ($0<\sigma<\nu$ if $\zeta=0$), $\alpha\in\mathbb{R}$, and $\beta\in\mathbb{R}$. Then we have the following.
    \begin{enumerate}[leftmargin=*,parsep=5pt]
        \item The operator $T^{0}(\zeta,\sigma,\alpha,\beta)$ is bounded from $L^{p}(\mathbb{X})$ to $L^{q}(\mathbb{X})$ if and only if $1\le{p}\le{q}\le+\infty$, $\frac{1}{p'}>\frac{\alpha}{n}$ when $\alpha>0$, $\frac{1}{q}>\frac{\beta}{n}$ when $\beta>0$, $0\le\alpha+\beta\le\sigma$, and $\frac{1}{p}-\frac{1}{q}\le\frac{\sigma-\alpha-\beta}{n}$, excluding the points $(\frac{1}{p},\frac{1}{q})=(\frac{\sigma-\alpha-\beta}{n},0)$ and $(1,1-\frac{\sigma-\alpha-\beta}{n})$ when $\sigma-\alpha-\beta\le{n}$.
    
        \item The operator $T^{\infty}(\zeta,\sigma,\alpha,\beta)$ is bounded from $L^{p}(\mathbb{X})$ to $L^{q}(\mathbb{X})$ if and only if $1\le{p}\le{q}\le+\infty$, $\frac{1}{p'}>\frac{\alpha}{n}$ when $\alpha>0$, $\frac{1}{q}>\frac{\beta}{n}$ when $\beta>0$, and one of the following conditions holds:
        \begin{enumerate}[leftmargin=*,parsep=5pt,label=(\roman*)]
            \item if $\zeta>|\rho|$;
    
            \item if $0\le\zeta\le|\rho|$, then $\frac{1}{p}>\frac{1}{2}-\frac{\zeta}{2|\rho|}$ and $\frac{1}{q}<\frac{1}{2}+\frac{\zeta}{2|\rho|}$;
    
            \item if $0<\zeta\le|\rho|$, $\frac{1}{p}=\frac{1}{2}-\frac{\zeta}{2|\rho|}>\frac{1}{q}$ or $\frac{1}{q}=\frac{1}{2}+\frac{\zeta}{2|\rho|}<\frac{1}{p}$, then 
            \begin{align}\label{mainthm1 res1}
                \min\lbrace{\alpha,\beta}\rbrace\,
                >\,
                \frac{\sigma}{2}\,-\,
                \frac{\ell+1}{2}\big(\frac{1}{2}-\frac{\zeta}{2|\rho|}\big);
            \end{align}

            \item if $0<\zeta\le|\rho|$ and $\frac{1}{p}=\frac{1}{q}=\frac{1}{2}\pm\frac{\zeta}{2|\rho|}$, then \eqref{mainthm1 res1} holds and $\alpha+\beta\ge\frac{\sigma}{2}$;
    
            \item if $\zeta=0$, $p=2<q$ or $q=2>p$, then 
            \begin{align}\label{mainthm1 res2}
                \min\lbrace{\alpha,\beta}\rbrace\,
                >\,
                \sigma-\frac{\nu}{2};
            \end{align}
    
            \item if $\zeta=0$ and $p=q=2$, then \eqref{mainthm1 res2} holds and $\alpha+\beta\ge\sigma$.
        \end{enumerate}
    \end{enumerate}
    Then, the operator $T(\zeta,\sigma,\alpha,\beta)$ is bounded from $L^{p}(\mathbb{X})$ to $L^{q}(\mathbb{X})$ if and only if $T^{0}(\zeta,\sigma,\alpha,\beta)$ and $T^{\infty}(\zeta,\sigma,\alpha,\beta)$ are both $L^p$-$L^q$-bounded.
    \end{theorem}

Theorem \ref{mainthm1} covers previously known results regarding the operator $R(\zeta,\sigma)=T(\zeta,\sigma,0,0)$, namely the Hardy-Littlewood-Sobolev inequality on $\mathbb{X}$. For the convenience of the readers, we provide a brief overview of these results. The study of $L^p$-$L^q$-boundedness of $R(\zeta,\sigma)$ emerged progressively during the 1970s. The initial investigation by Stanton and Tomas \cite{ST78} focused on the expansion of spherical functions and examined the operator $R(0,\sigma)$ ($\sigma\in\mathbb{C}$ with $\re\sigma\ge0$) in symmetric spaces of rank one. In the same context, Lohoué and Rychener investigated the resolvent $R(\zeta,2)$ in \cite{LR82} (also explored in \cite{DST88}). Anker and Lohoué extended these rank one results to symmetric spaces with normal real form  \cite{AL86}. Lohoué further generalized Varopoulos's findings from the operator $R(0,1)$ to the entire family of operators $R(0,\sigma)$ for all $\sigma\in\mathbb{C}$ with $\re{\sigma}\ge0$ on general symmetric spaces of non-compact type \cite{Var88,Loh89}. Subsequently, Taylor \cite{Tay89} and Anker \cite{Ank90,Ank92} provided a nearly complete description of the $L^p$-$L^q$-boundedness for the entire family of operators $R(\zeta,\sigma)$, excluding certain critical indices that were finally addressed in \cite{CGM93}. Some of the aforementioned results are also available on more general manifolds beyond symmetric spaces, as seen in \cite{Ste70,Loh80,Loh85,Cow83,Str83}.

Let us illustrate the properties in Theorem \ref{mainthm1} with figures. In the following figures, the blue zone represents the range of admissible values of $(\frac{1}{p},\frac{1}{q})$ (the dashed lines and the open points are not included), within which the operator is bounded from $L^{p}(\mathbb{X})$ to $L^{q}(\mathbb{X})$. The red lines or points mean that the operator could be $L^p$-$L^q$-bounded, but with additional assumptions. Outside these admissible regions, the corresponding operator cannot be $L^p$-$L^q$-bounded. 

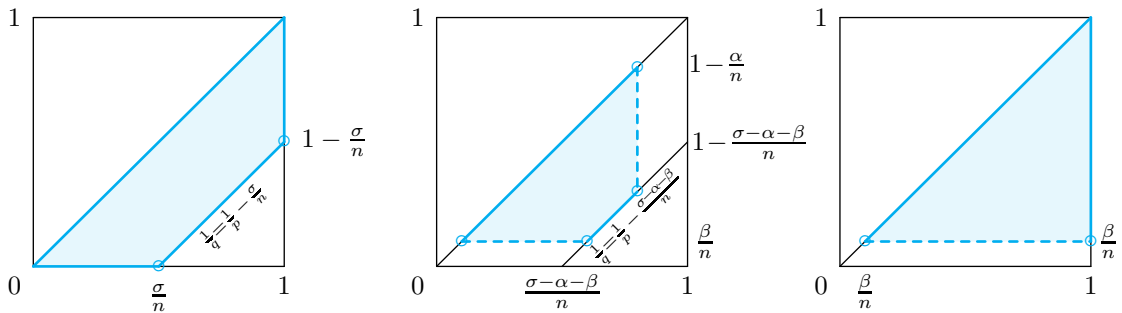
\begin{figure}[b]
    \centering
        \begin{subfigure}{0.33\textwidth}
            \input{HLS02.tex}
            \noindent\subcaption{$\sigma\!<\!n$ and $\alpha\!=\!\beta\!=\!0$}
        \end{subfigure}
        ~
        \begin{subfigure}{0.33\textwidth}
            \input{SW01.tex}
            \subcaption{$\sigma\!-\!\alpha\!-\!\beta\!<\!n$ and $\alpha\!>\!\beta>0$}
        \end{subfigure}
        ~
        \begin{subfigure}{0.33\textwidth}
            \input{SW03.tex}
            \subcaption{$\sigma\!-\!\alpha\!-\!\beta\!\ge\!{n}$ and $\alpha\!\le\!0\!<\!\beta$}
        \end{subfigure}
        \caption{Examples of the $L^{p}$- $L^{q}$-boundedness of the operator $T^{0}(\zeta,\sigma,\alpha,\beta)$.}
        \label{SW0Fig}
    \end{figure}

At first, Figure \ref{SW0Fig} shows several examples of local properties in Theorem \ref{mainthm1}.\textit{(1)}. The segment $\frac{1}{q}=\frac{1}{p}-\frac{\sigma-\alpha-\beta}{n}$ contracts to the point $(\frac{1}{p},\frac{1}{q})=(1,0)$ when $\sigma-\alpha-\beta$ tends to $n$. This means particularly that the admissible region in Figure \ref{SW0Fig}.(i) can be widened to the entire triangle when $\sigma>{n}$. It is noteworthy that the constraints $\frac{1}{p'}>\frac{\alpha}{n}$ and $\frac{1}{q}>\frac{\beta}{n}$ when $\alpha,\,\beta>0$ in the Stein-Weiss inequality arise naturally from the local behavior of the homogeneous weights. By combining them with the condition $\alpha+\beta\ge0$, one notices that the point $(\frac{1}{p},\frac{1}{q})=(1,0)$ is always excluded if $(\alpha,\beta)\neq(0,0)$. In other words, the operator $T^{0}(\zeta,\sigma,\alpha,\beta)$ cannot be bounded from $L^{1}(\mathbb{X})$ to $L^{\infty}(\mathbb{X})$ if a weight exists, even for $\sigma$ large. However, the operator $T^{0}(\zeta,\sigma,0,0)$ is known to be $L^{1}$-$L^{\infty}$-bounded for all $\zeta\ge0$ when $\sigma>n$, according to Theorem \ref{mainthm1}.\textit{(1)} with $\alpha=\beta=0$.

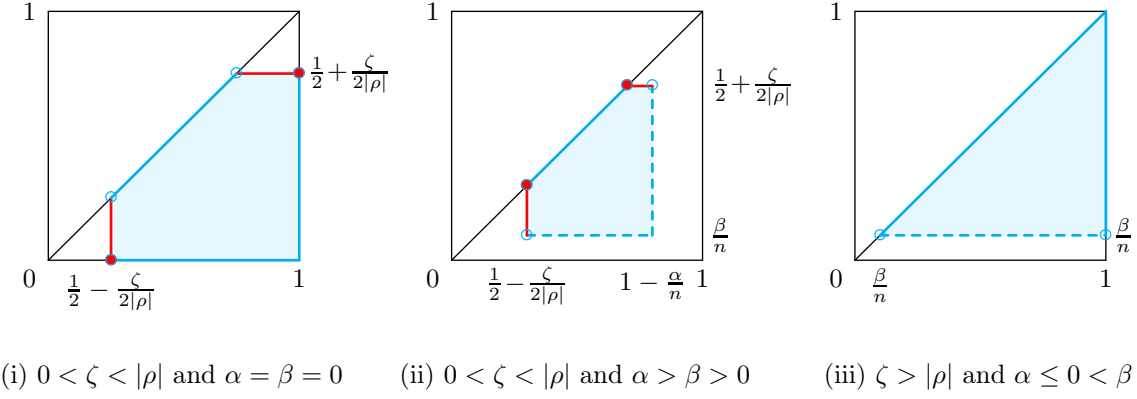
\begin{figure}
    \centering
    \begin{subfigure}{0.33\textwidth}
        \input{HLSinf2.tex}
        \subcaption{$0<\zeta<|\rho|$ and $\alpha=\beta=0$}
    \end{subfigure}
    ~
    \begin{subfigure}{0.33\textwidth}
        \input{SWinf2.tex}
        \subcaption{$0<\zeta<|\rho|$ and $\alpha>\beta>0$}
    \end{subfigure}
    ~
    \begin{subfigure}{0.33\textwidth}
        \input{SWinf1.tex}
        \subcaption{$\zeta>|\rho|$ and $\alpha\le0<\beta$}
    \end{subfigure}
       \caption{Examples of the $L^p$-$L^q$-boundedness of $T^{\infty}(\zeta,\sigma,\alpha,\beta)$ for $\zeta>0$.}
    \label{SWinfFig}
   \end{figure}

A peculiarity of symmetric spaces is their exponential growth at infinity, which has significant implications for the behavior of the operator $T^{\infty}(\zeta,\sigma,\alpha,\beta)$ and exhibits phenomena very different from those in the Euclidean setting. In simple terms, in certain situations, when we have a sufficiently strong exponential decay from the kernel estimate \eqref{estim kinf}, the role played by polynomial factors becomes less essential. This relaxes the balance condition $\frac{1}{p}-\frac{1}{q}=\frac{\sigma-\alpha-\beta}{n}$, which is typically required in Euclidean analysis, or more general, on manifolds with the dilation property, and provides a wider range of  admissible pairs $(\frac{1}{p},\frac{1}{q})$, see Figure \ref{SWinfFig}. On the critical lines (represented as red segments in Figure \ref{SWinfFig}) where the decay obtained from the kernel precisely offsets the exponential growth of the volume, the operator $T^{\infty}(\zeta,\sigma,\alpha,\beta)$ can still be $L^p$-$L^q$-bounded by carefully determining the weights.

Once we fix two weights with $\alpha>0$ and $\beta>0$, and take into account their local behaviors, it becomes evident that the region of the admissible pairs $(\frac{1}{p},\frac{1}{q})$ ensuring the $L^p$-$L^q$-boundedness of the operator $T^{\infty}(\zeta,\sigma,\alpha,\beta)$ is more constrained compared to that of $T^{\infty}(\zeta,\sigma,0,0)$. However, we observe that $T^{\infty}(\zeta,\sigma,\alpha,\beta)$ can be bounded in the diagonal critical case where $p=q$, in contrast to the operator $T^{\infty}(\zeta,\sigma,0,0)$, as depicted in subfigures \textit{(i)} and \textit{(ii)} of Figure \ref{SWinfFig}. In fact, if $\alpha=\beta=0$ are taken in Theorem \ref{mainthm1}.\textit{(2).(iv)}, there will be no valid $\sigma$. That is, the operator $T^{\infty}(\zeta,\sigma,0,0)$ cannot be $L^p$-bounded when $\frac{1}{p}$ locates on the red critical lines in Figure \ref{SWinfFig}.(i), see also \cite[Proposition 4.3]{CGM93}. In this sense, the Stein-Weiss inequality can be seen as a diagonal improvement of the Hardy-Littlewood-Sobolev inequality.

Observe that the two red segments in Subfigures \textit{(i)} and \textit{(ii)} of Figure \ref{SWinfFig} contract to the points $(0,0)$ and $(1,1)$ respectively as $\zeta$ approaches $|\rho|$. However, it is impossible for the operator $T^{\infty}(|\rho|,\sigma,\alpha,\beta)$ to be $L^{1}$-bounded or $L^{\infty}$-bounded in both the weighted and unweighted cases. On the one hand, the operator $T^{\infty}(|\rho|,\sigma,0,0)$ is not $L^p$-bounded for $p=1$ or $p=+\infty$ as we just mentioned. On the other hand, the condition \eqref{mainthm1 res1} tells us that $\alpha$ and $\beta$ should both be positive when $\zeta=|\rho|$. Therefore, the operator $T^{\infty}(|\rho|,\sigma,\alpha,\beta)$ cannot be $L^1$-bounded or $L^{\infty}$-bounded due to the conditions $\frac{1}{p'}\ge\frac{\alpha}{n}$ and $\frac{1}{p}\ge\frac{\beta}{n}$ that are required for the local parts. Note also that Subfigure \ref{SW0Fig}.\textit{(iii)} and Subfigure \ref{SWinfFig}.\textit{(iii)} exhibit identical shapes. However, the boundedness of $T^{0}(\zeta,\sigma,\alpha,\beta)$ is affected by the parameter $\sigma$, whereas the boundedness of $T^{\infty}(\zeta,\sigma,\alpha,\beta)$ depends on the parameter $\zeta$. 

\begin{figure}[b]
    \centering
        \begin{subfigure}{0.5\textwidth}
            \hspace{50pt}
            \input{HLSinf3.tex}
            \noindent\subcaption{$\alpha=\beta=0$}
        \end{subfigure}
        ~
        \begin{subfigure}{0.5\textwidth}
            \hspace{35pt}
            \input{SWinf3.tex}
            \subcaption{$\alpha>0$ and $\beta>0$}
        \end{subfigure}
        \caption{Examples of the $L^{p}$- $L^{q}$-boundedness of $T^{\infty}(0,\sigma,\alpha,\beta)$}
        \label{SWlim}
\end{figure}
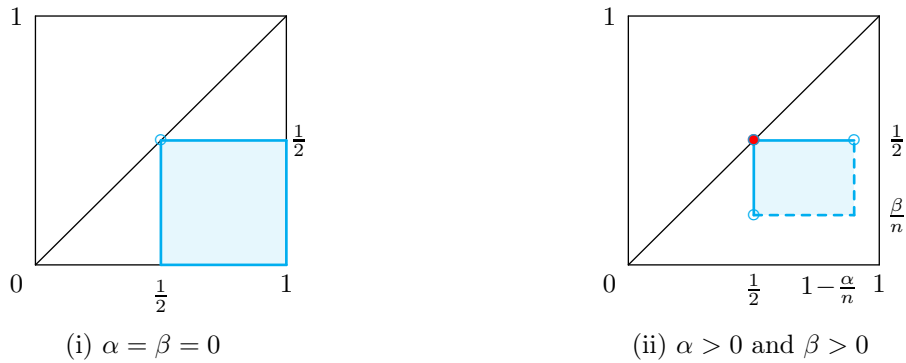

The same phenomena occur in the limiting case where $\zeta=0$. Theorem \ref{mainthm1}.\textit{(2).(vi)} demonstrates that the operator $T^{\infty}(0,\sigma,\alpha,\beta)$ can be bounded from $L^2{(\mathbb{X})}$ to $L^2{(\mathbb{X})}$, provided that $\alpha+\beta\ge\sigma$. This condition rules out the possibility of the operator $R^{\infty}(0,\sigma)=T^{\infty}(0,\sigma,0,0)$ being $L^2$-bounded, see Figure \ref{SWlim}. Recall that, in this limiting case, the intersection of the domain of $R^{\infty}(0,\sigma)$ on $L^2(\mathbb{X})$ with $L^p(\mathbb{X})$ is not dense in $L^p(\mathbb{X})$ when $p<2$ and $\sigma\ge\frac{\nu}{2}$. In the study of the multiplier theory, one defines $R^{\infty}(0,\sigma)$ on $L^p(\mathbb{X})$ by means of analytic continuation by considering complex exponent $\sigma$. In this case, if $\re{\sigma}$ is non-negative and $\re{\sigma}-\nu$ is not a non-negative even integer, the operator $R^{\infty}(0,\sigma)$ is $L^p$-$L^q$-bounded for all $1\le{p}<2<q\le+\infty$, see \cite[Theorem 6.1]{CGM93}. Note also that $R^{\infty}(0,\sigma)$ is not $L^p$-$L^q$-bounded for any $\sigma\ge\frac{\nu}{2}$ if $p=2$ or $q=2$, and not $L^2$-bounded for any $\sigma>0$, see \cite{Loh89}. This paper focuses on the Stein-Weiss inequality, we will consider only real exponents. In this case, Theorem \ref{mainthm1} with $\alpha=\beta=0$ is compatible with these earlier results.

As a consequence of Theorem \ref{mainthm1}, we present the following Stein-Weiss inequality.
\begin{theorem}[Stein-Weiss inequality on non-compact symmetric spaces]\label{mainthm2}
    Let $\mathbb{X}$ be a non-compact symmetric space of rank $\ell\ge1$, dimension $n\ge2$, and pseudo-dimension $\nu\ge3$. Suppose that $\zeta\ge0$, $\sigma>0$ ($0<\sigma<\nu$ if $\zeta=0$), $\alpha\in\mathbb{R}$, and $\beta\in\mathbb{R}$. Then, the inequality
        \begin{align}\label{SW}
        \||\cdot|^{-\beta}\,(-\Delta-|\rho|^{2}+\zeta^{2})^{-\frac{\sigma}{2}}(|\cdot|^{-\alpha}f)\|_{L^{q}(\mathbb{X})}\,
        \lesssim\,
        \|f\|_{L^{p}(\mathbb{X})}
        \qquad\forall\,f\in{L^{p}(\mathbb{X})},
        \end{align}
    holds if and only if $1\le{p}\le{q}\le+\infty$, $\frac{1}{p'}>\frac{\alpha}{n}$ when $\alpha>0$, $\frac{1}{q}>\frac{\beta}{n}$ when $\beta>0$, $0\le\alpha+\beta\le\sigma$, and $\frac{1}{p}-\frac{1}{q}\le\frac{\sigma-\alpha-\beta}{n}$ (excluding the points $(\frac{1}{p},\frac{1}{q})=(\frac{\sigma-\alpha-\beta}{n},0)$ and $(1,1-\frac{\sigma-\alpha-\beta}{n})$ when $\sigma-\alpha-\beta\le{n}$), and one of the following conditions is met:
        \begin{enumerate}[leftmargin=*,parsep=5pt,label=(\roman*)]
            \item if $\zeta>|\rho|$;

            \item if $0\le\zeta\le|\rho|$, then $\frac{1}{p}>\frac{1}{2}-\frac{\zeta}{2|\rho|}$ and $\frac{1}{q}<\frac{1}{2}+\frac{\zeta}{2|\rho|}$;
    
            \item if $0<\zeta\le|\rho|$, $\frac{1}{q}<\frac{1}{2}-\frac{\zeta}{2|\rho|}=\frac{1}{p}$ or $\frac{1}{q}=\frac{1}{2}+\frac{\zeta}{2|\rho|}<\frac{1}{p}$, then $\min\lbrace{\alpha,\beta}\rbrace>\frac{\sigma}{2}-\frac{\ell+1}{2}(\frac{1}{2}-\frac{\zeta}{2|\rho|})$;
            
            \item if $0<\zeta\le|\rho|$ and $\frac{1}{p}=\frac{1}{q}=\frac{1}{2}\pm\frac{\zeta}{2|\rho|}$, then $\min\lbrace{\alpha,\beta}\rbrace>\frac{\sigma}{2}-\frac{\ell+1}{2}(\frac{1}{2}-\frac{\zeta}{2|\rho|})$ and $\alpha+\beta\ge\frac{\sigma}{2}$;
            
            \item if $\zeta=0$, $p=2<q$ or $p<2=q$, then $\min\lbrace{\alpha,\beta}\rbrace>\sigma-\frac{\nu}{2}$;
            
            \item if $\zeta=0$ and $p=q=2$, then $\min\lbrace{\alpha,\beta}\rbrace>\sigma-\frac{\nu}{2}$ and $\alpha+\beta\ge\sigma$.
        \end{enumerate}
\end{theorem}

\begin{remark}
    The real hyperbolic space $\mathbb{H}^{n}$ is the most common example of negatively curved manifolds. In particular, it is a symmetric space of rank $\ell=1$, which always has pseudo-dimension $\nu=3$. In this case, $\rho=\frac{n-1}{2}$ is a positive real number. Substituting these values into Theorem \ref{mainthm2} yields the Stein-Weiss inequality in real hyperbolic spaces.
\end{remark}

\begin{remark}
    Together with the assumptions of the local part, the condition $\textit{(vi)}$ in Theorem \ref{mainthm2} can be written as: if $\zeta=0$ and $p=q=2$, then $\alpha+\beta=\sigma$ and $\max\lbrace{\alpha,\beta}\rbrace<\frac{\nu}{2}$. This $L^2$-version of the Stein-Weiss inequality was previously established in \cite{KRZ23}. 
\end{remark}    

\begin{figure}
    \centering
       \input{MainTHM.tex}
       \caption{Case $0<\zeta<|\rho|$, $0<\sigma<n$, $\alpha>0$, and $\beta>0$.}
    \label{SWFig}
\end{figure}
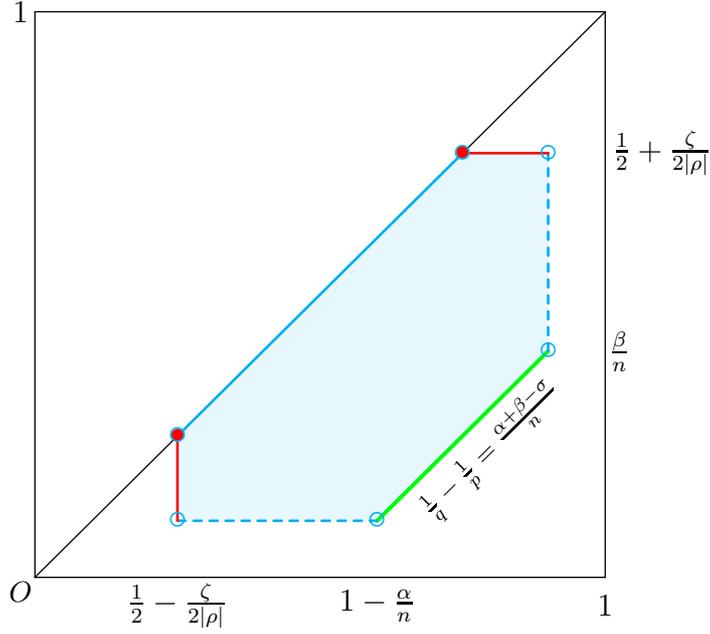

\begin{remark}
    In the recent paper \cite{KKR22}, the authors studied the operator $T(\zeta,\sigma,\alpha,\beta)$ for $n\ge3$,  $0<\sigma<n$ and $\zeta>0$ sufficiently large, which corresponds to a subcase of Theorem \ref{mainthm2}.\textit{(i)}. In this particular scenario, the associated convolution kernel enjoys a very large exponential decay, see estimate \eqref{estim kinf}. This allows one to avoid the main difficulty arising from the geometry at infinity. This also explains why dealing with the large shift case \textit{(i)} and the non-critical case \textit{(ii)} in Theorem \ref{mainthm2} is relatively easier, see Section \ref{Section Proof}. Moreover, the balance condition $\frac{1}{p}-\frac{1}{q}=\frac{\sigma-\alpha-\beta}{n}$ is required in \cite[Theorem 1.2]{KKR22}, which means that their set of admissible pairs $(\frac{1}{p},\frac{1}{q})$ is restricted to a segment as in the Euclidean setting, see Figure \ref{SWFig}. Furthermore, to the best of our knowledge, none of the aforementioned works has investigated the endpoints $p=1$ and $q=+\infty$ and thoroughly considered the necessity aspect.
\end{remark}

There is another example to consider. Suppose that $0<\zeta<|\rho|$, $0<\sigma<n$, $\alpha>0$, and $\beta>0$. In this case, the Stein-Weiss inequality \eqref{SW} holds if and only if the pair $(\frac{1}{p},\frac{1}{q})$ belongs to the crystal shape region of Figure \ref{SWFig} (the dashed lines and the open points are not included), with more constraints on the two red critical lines. We highlight that this admissible region is notably larger compared to the one in the Euclidean setting or on manifolds with dilation property. In those contexts, the admissible set is restricted to the green line due to the aforementioned balance condition, which is no longer applicable in our setting because of the special geometry at infinity of symmetric spaces.

\subsection{Strategy of proof}\label{Sect Strategy}
Contrasting with previous works on the Hardy-Littlewood-Sobolev inequality, our approach does not require relating the operator to the heat kernel through the use of Mellin-type transforms or employing Fourier analysis techniques. Instead, we establish Theorem \ref{mainthm1} by skillfully combining pointwise kernel estimates with various formulations of the Kunze-Stein phenomena. 

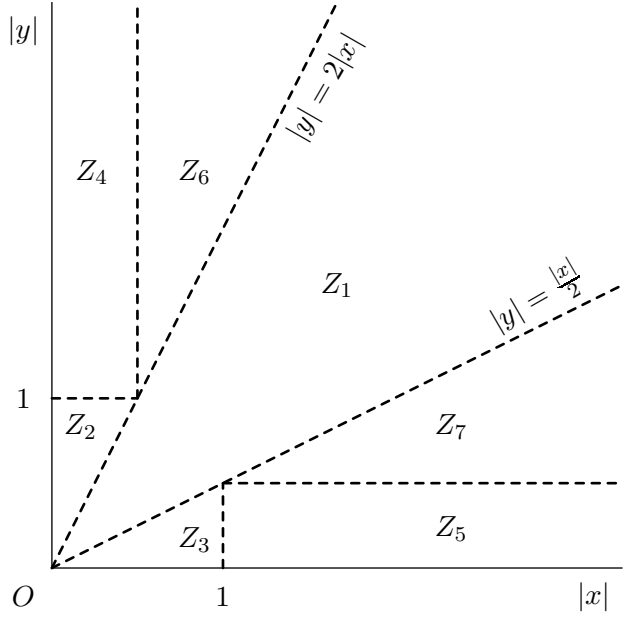
\begin{figure}
    \centering
       \input{Regions.tex}
    \caption{Different regions.}
    \label{Regions}
\end{figure}

To achieve this, we divide the product space $\mathbb{X}\times\mathbb{X}$ into the following disjoint regions (refer to Figure \ref{Regions}):
\begin{itemize}[parsep=5pt]
    \item 
    $Z_{1}=\lbrace{(x,y)\in\mathbb{X}\times\mathbb{X}
    \,|\,\frac{|x|}{2}<|y|<2|x|}\rbrace$;

    \item 
    $Z_{2}=\lbrace{(x,y)\in\mathbb{X}\times\mathbb{X}
    \,|\,2|x|\le|y|\le1}\rbrace$ and 
    $Z_{3}=\lbrace{(x,y)\in\mathbb{X}\times\mathbb{X}
    \,|\,2|y|\le|x|\le1}\rbrace$;

    \item 
    $Z_{4}=\lbrace{(x,y)\in\mathbb{X}\times\mathbb{X}
    \,|\,2|x|\le1\le|y|}\rbrace$ and
    $Z_{5}=\lbrace{(x,y)\in\mathbb{X}\times\mathbb{X}
    \,|\,2|y|\le1\le|x|}\rbrace$;

    \item 
    $Z_{6}=\lbrace{(x,y)\in\mathbb{X}\times\mathbb{X}
    \,|\,1\le2|x|\le|y|}\rbrace$ and
    $Z_{7}=\lbrace{(x,y)\in\mathbb{X}\times\mathbb{X}
    \,|\,1\le2|y|\le|x|}\rbrace$.
\end{itemize}

Due to the different behaviors of the convolution kernel, this partition is more refined compared to the one considered in \cite{SW58} for the Euclidean setting. This paper aims to determine all possible values of $\zeta\ge0$, $\sigma>0$, $\alpha\in\mathbb{R}$, $\beta\in\mathbb{R}$, and $1\le{p},\,{q}\le+\infty$, for which the operators 
\begin{align*}
    T_{j}^{0}(\zeta,\sigma,\alpha,\beta)f(x)\,=\,
    \int_{\mathbb{X}}\diff{y}\,
    \mathds{1}_{Z_{j}}(x,y)\,
    \psi_{0}(y^{-1}x)\,
    |x|^{-\beta}\,k_{\zeta,\sigma}(y^{-1}x)\,
    |y|^{-\alpha}\,f(y),
\end{align*}
and
\begin{align*}
    T_{j}^{\infty}(\zeta,\sigma,\alpha,\beta)f(x)\,=\,
    \int_{\mathbb{X}}\diff{y}\,
    \mathds{1}_{Z_{j}}(x,y)\,\psi_{\infty}(y^{-1}x)\,
    |x|^{-\beta}\,k_{\zeta,\sigma}(y^{-1}x)\,
    |y|^{-\alpha}\,f(y),
\end{align*}
are bounded from $L^p(\mathbb{X})$ into $L^q(\mathbb{X})$. Then, the Stein-Weiss inequality follows as a consequence. Recall that $\psi_{0}$ is a bi-$K$-invariant cut-off function on $\mathbb{X}$ such that its support is contained in the subset $\lbrace{x\in\mathbb{X}\,|\,|x|\le1}\rbrace$, and $\psi_{\infty}=1-\psi_{0}$. Therefore, for all $f\in{L^{p}(\mathbb{X})}$, we have $T_{j}^{0}(\zeta,\sigma,\alpha,\beta)f=0$ when $4\le{j}\le7$ and $T_{j}^{\infty}(\zeta,\sigma,\alpha,\beta)f=0$ when $j=2,3$. In the local part, the analysis shows a Euclidean-like nature, and we extend the argument presented in \cite{SW58} to symmetric spaces. The main ingredient is a generalization of the Hardy-Littlewood-P\'olya type estimate. While the study of $T_{j}^{\infty}$ builds upon the Kunze-Stein phenomena, which we will review in the next section. The proof of Theorem \ref{mainthm1} is given in Section \ref{Section Proof}. 

Because of the implicit constants in both the kernel estimates and the Kunze-Stein phenomena, we will not search for the optimal constant in \eqref{SW}.

\subsection{Applications}
Let us conclude this introductory section with some applications. As mentioned above, the Stein-Weiss inequality serves as a bridge to many other weighted functional inequalities. In this section we will highlight some prominent examples. For $\sigma>0$ and $1<p<+\infty$, we define the Sobolev space $H^{\sigma,p}(\mathbb{X})$ as the image of $L^{p}(\mathbb{X})$ under the operator $(-\Delta)^{\frac{\sigma}{2}}$, equipped with the norm $\|f\|_{H^{\sigma,p}(\mathbb{X})}=\|(-\Delta)^{\frac{\sigma}{2}}f\|_{L^{p}(\mathbb{X})}$. Our first example is the following.
\begin{corollary}[Hardy-Sobolev inequality]
    Let $\mathbb{X}$ be a symmetric space of dimension $n\ge2$.
    Suppose that $\sigma>0$, $\beta\in\mathbb{R}$, and $1<p\le{q}<+\infty$. Then, the inequality
    \begin{align}
        \||\cdot|^{-\beta}f\|_{L^{q}(\mathbb{X})}\,
        \lesssim\,\|f\|_{H^{\sigma,p}(\mathbb{X})}
        \qquad\forall\,f\in{H^{\sigma,p}(\mathbb{X})}
        \label{HardySobolev}
    \end{align}
    holds if and only if $0\le\beta<\frac{n}{q}$ and $\frac{1}{p}-\frac{1}{q}\le\frac{\sigma-\beta}{n}$.
\end{corollary}
This corollary is a direct consequence of Theorem \ref{mainthm2} if one takes $\alpha=0$ and $\zeta=|\rho|$. The inequality \eqref{HardySobolev} is known as the \textit{Hardy inequality} if $p=q$ and as the \textit{Sobolev inequality} if $\beta=0$. As explained in Figure \ref{SWFig}, the balance condition $\frac{1}{p}-\frac{1}{q}=\frac{\sigma-\beta}{n}$ has been relaxed because of the special geometry at infinity of the symmetric space. Consequently, the Hardy inequality 
\begin{align}
    \||\cdot|^{-\beta}f\|_{L^{p}(\mathbb{X})}\,
    \lesssim\,\|f\|_{H^{\sigma,p}(\mathbb{X})}
    \label{Hardy}
\end{align}
holds for a wider range $0\le\beta\le\sigma$ (and $\beta<\frac{n}{q}$). Recall that an estimate similar to \eqref{Hardy} is valid in the Euclidean setting only if $\beta=\sigma<\frac{n}{q}$. This inequality can be seen as a mathematical formulation of the uncertainty principle. It states that if a function $f$ is localized around the origin in $\mathbb{R}^{N}$, then its momentum is necessarily large. However, in our geometry, the stronger singularity $|x|^{-\sigma}$ in $\mathbb{R}^{N}$ can be improved as $|x|^{-\beta}$ with $0\le\beta\le\sigma$, which means that the momentum $\|f\|_{H^{\sigma,p}(\mathbb{X})}$ could be not that large.

From the Hardy inequality \eqref{Hardy} and the Hölder inequality we can also derive the following \textit{Heisenberg's uncertainty principle}.
\begin{corollary}[Heisenberg's uncertainty principle]
    Let $\mathbb{X}$ be a symmetric space of dimension $n\ge2$.
    Suppose that $\sigma>0$, $1<p<+\infty$, $0\le\beta\le\sigma$, and $\beta<\frac{n}{p}$. Then, the inequality
    \begin{align}
        \|f\|_{H^{\sigma,p}(\mathbb{X})}\,
        \||\cdot|^{\beta}f\|_{L^{p'}(\mathbb{X})}\,
        \gtrsim\,
        \|f\|_{L^{2}(\mathbb{X})}^{2}
        \label{Heisenberg}
    \end{align}
    holds for any $f\in{H^{\sigma,p}(\mathbb{X})}$ such that $(|\cdot|^{\beta}f)\in{L^{p'}(\mathbb{X})}$.
\end{corollary}
Inequality \eqref{Heisenberg} shows that if the momentum $\|f\|_{H^{\sigma,p}(\mathbb{X})}$ is not too large, then the followed term is not too small. However, here we have more flexible options for the term $\||\,\cdot\,|^{\beta}f\|_{L^{p'}(\mathbb{X})}$ compared to the Euclidean setting where a similar estimate holds only with $\beta=\sigma$.

Weighted inequalities for the Fourier transform can also serve as a natural approach to understanding uncertainty. We denote by $\mathcal{F}:\mathbb{X}\rightarrow\mathfrak{a}\times{B}$ the Helgason-Fourier transform that maps a Schwartz function from the symmetric space $\mathbb{X}$ to $\mathfrak{a}\times{B}$, where $B=K/C_{K}(\exp\mathfrak{a})$ and $C_{K}(\exp\mathfrak{a})$ is the centralizer of $\exp\mathfrak{a}$ in $K$. According to the Plancherel theorem (see, for instance, \cite[Theorem 1.5, p.227]{Hel94}), this transform $\mathcal{F}$ extends to an isometry of $L^{2}(\mathbb{X})$ to $L^2(\mathfrak{a}\times{B},\,|W|^{-1}|\mathbf{c}(\lambda)|^{-2}\diff{\lambda})$, where $W$ is the Weyl group and $\mathbf{c}(\lambda)$ denotes the Harish-Chandra $\mathbf{c}$-function. Using Theorem \ref{mainthm2} with $\zeta=0$, $\beta=0$, and $q=2$, and the Plancherel theorem, we obtain the following corollary.
\begin{corollary}[Pitt's inequality]
    Let $\mathbb{X}$ be a symmetric space of dimension $n\ge2$ and pseudo-dimension $\nu\ge3$.
    Suppose that $0<\sigma<\nu$ and $f\in\mathcal{S}(\mathbb{X})$. Then, the inequality 
    \begin{align}
        \||\cdot|^{-\sigma}\mathcal{F}f\|_{L^{2}(\mathfrak{a}\times{B},\,|W|^{-1}|\mathbf{c}(\lambda)|^{-2}\diff{\lambda})}\,
        \lesssim\,
        \||\cdot|^{\alpha}f\|_{L^{p}(\mathbb{X})}
        \label{Pitt}
    \end{align}
    holds for $1\le{p}<2$ if and only if $0\le\alpha\le\sigma<\frac{\nu}{2}$, $\alpha<\frac{n}{p'}$, and $\frac{1}{p}-\frac{1}{2}\le\frac{\sigma-\alpha}{n}$. Moreover, the inequality \eqref{Pitt} holds for $p=2$ if and only if $\alpha=\sigma<\min\lbrace{\frac{n}{2},\,\frac{\nu}{2}}\rbrace$.
\end{corollary}
Although we need the extra condition $\sigma<\frac{\nu}{2}$ in the $L^2$-case, since the topological dimension and the pseudo-dimension both play their roles, note that inequality \eqref{Pitt} holds for a wider range of $p$ when $p<2$. In the Euclidean case, the balance condition $\frac{1}{p}-\frac{1}{2}=\frac{\sigma-\alpha}{n}$ is always required, see, for instance, \cite{Bec95,Bec08a,DCGT17,DNTW23}.

On the other hand, by taking $\beta=0$ in \eqref{HardySobolev}, we know that the Sobolev inequality 
\begin{align}
    \|f\|_{L^{q}(\mathbb{X})}\,
    \lesssim\,\|f\|_{H^{\sigma,p}(\mathbb{X})}
    \label{Sobolev}
\end{align}
holds for $\sigma>0$ and $1<p\le{q}<+\infty$ if and only if $\frac{1}{p}-\frac{1}{q}\le\frac{\sigma}{n}$. Note in particular that we do not need to assume $\sigma<n$, and that the condition $\frac{1}{p}-\frac{1}{q}\le\frac{\sigma}{n}$ is always fulfilled if $\sigma\ge{n}$. Therefore, the set of admissible pairs $(\frac{1}{p},\frac{1}{q})$ that satisfy the Sobolev inequality \eqref{Sobolev} is also much larger in comparison to the Euclidean space $\mathbb{R}^{N}$, in which the balance condition $\frac{1}{p}-\frac{1}{q}=\frac{\sigma}{N}$ with $\sigma<{N}$ is necessary. Furthermore, the inequality \eqref{Sobolev} remains valid for $1<p<+\infty$ and $\sigma\ge\frac{n}{p}$ in the endpoint case where $q=+\infty$. Such a Sobolev embedding has already been established in \cite[Lemma 2.2]{Ank92} for $p=2$. We also draw the reader's attention to its bi-$K$-invariant improvement that was proven in \cite[Lemma 2.3]{Ank92}.

By reapplying both the Sobolev inequality \eqref{Sobolev} and the Hölder inequality, we arrive at a more general form of the inequality, known as the \textit{Gagliardo-Nirenberg interpolation inequality}.
\begin{corollary}[Gagliardo-Nirenberg interpolation inequality]
    Let $\mathbb{X}$ be a symmetric space of dimension $n\ge2$.
    Suppose that $\sigma>0$, $0\le\theta\le1$, $1<s<+\infty$, and $1\le{p,\,q}\le+\infty$ satisfy the condition
    \begin{align*}
        \frac{1}{q}\,
        \ge\,
        \begin{cases}
            \theta(\frac{1}{s}-\frac{\sigma}{n})\,
            +\,\frac{1-\theta}{p}
            &\qquad\textnormal{if}\,\,\,\sigma<n,\\[5pt]
            \frac{1-\theta}{p}
            &\qquad\textnormal{if}\,\,\,\sigma\ge{n}.
        \end{cases}
    \end{align*}
    Then, the inequality
    \begin{align}
        \|f\|_{L^{q}(\mathbb{X})}\,
        \lesssim\,
        \|f\|_{H^{\sigma,s}(\mathbb{X})}^{\theta}\,
        \|f\|_{L^{p}(\mathbb{X})}^{1-\theta}
    \end{align}
    holds for any $f\in{L^{p}(\mathbb{X})}$ such that $f\in{H^{\sigma,s}(\mathbb{X})}$.
\end{corollary}
Similarly, for the validity of the Gagliardo-Nirenberg inequality in $\mathbb{R}^N$, it is essential that the balance condition $\frac{1}{q}=\theta(\frac{1}{s}-\frac{\sigma}{N})+\frac{1-\theta}{p}$ should be satisfied. This constraint is relaxed due to the special geometry at infinity of the symmetric spaces, as we have mentioned several times above.

Another example is the following Poincaré-Sobolev inequality for the fractional Laplace-Beltrami operator, which has been of recent interest. The following corollary is obtained from Theorem \ref{mainthm2} by taking $\zeta=0$, $\sigma=\frac{\sigma'}{2}$, $\alpha=\beta=0$, and $p=2$.
\begin{corollary}[Poincaré-Sobolev inequality]
    Let $\mathbb{X}$ be a symmetric space of dimension $n\ge2$ and pseudo-dimension $\nu\ge3$. 
    Then, the inequality
    \begin{align}
        \|(-\Delta-|\rho|^{2})^{\frac{\sigma'}{4}}f\|_{L^{2}(\mathbb{X})}\,
        \gtrsim\,
        \|f\|_{L^{q}(\mathbb{X})}
        \qquad\forall\,f\in{H^{\frac{\sigma'}{2},2}(\mathbb{X})}
        \label{Poincare}
    \end{align}
    holds if and only if $0<\sigma'<\nu$ and $max\lbrace{0,\,\frac{1}{2}-\frac{\sigma'}{2n}}\rbrace\le\frac{1}{q}<\frac{1}{2}$.
\end{corollary}
It has recently been proved in \cite[Theorem 6.2]{LLY20} and \cite[Theorem 1.12]{BP22} that inequality \eqref{Poincare} holds when $n\ge3$, $0<\sigma'<\min\lbrace{n,\nu}\rbrace$, and $\frac{1}{2}-\frac{\sigma'}{2n}\le\frac{1}{q}<\frac{1}{2}$. Since the conditions in Theorem \ref{mainthm2} are both sufficient and necessary, this corollary states that the inequality \eqref{Poincare} remains valid on the hyperbolic plane where $n=2$, and the condition $0<\sigma'<\nu$ is necessary, but $\sigma'<n$ is not. When $\sigma'\ge{n}$, the inequality \eqref{Poincare} holds for all $2<q\le+\infty$.

The last example is a direct application in the study of partial differential equations. It is well known that the Stein-Weiss inequality is a key tool for establishing the smoothing properties of dispersive equations, see for example \cite{Sug03,RS06,Chi08,RS12,Kai14,RS16}. Consider the Schrödinger equation
\begin{align}\label{Schrodinger} 
    (i\partial_{t}+\Delta)\,u(t,x)\,=\,0,
    \qquad\,u(0,x)\,=u_{0}(x),
    \qquad\forall\,t\in\mathbb{R}^{*},\,\,\,
    \forall\,x\in\mathbb{X}.
\end{align}
We know that the solution operator $e^{it\Delta}$ preserves the $L^2$-norm for any fixed time. The smoothing property is used to obtain additional regularity by integrating the solution to \eqref{Schrodinger} in time. The following corollary, based on the $L^2$ Stein-Weiss inequality, was given in \cite{KRZ23}. In particular, it implies that, on a hyperbolic plane where $n=2$, the smoothing property holds for all $-\frac{1}{2}<\gamma<\frac{1}{2}$. Recall that, on a Euclidean plane, a similar property is possible if and only if $0<\gamma<\frac{1}{2}$.
\begin{corollary}[Kato-type smoothing property]
    Let $\mathbb{X}$ be a symmetric space of dimension $n\ge2$ and pseudo-dimension $\nu\ge3$. Suppose that either $1-\min\lbrace{\frac{n}{2},\,\frac{\nu}{2}}\rbrace<\gamma<\frac{1}{2}$ if $n\ge3$ or $-\frac{1}{2}<\gamma<\frac{1}{2}$ if $n=2$. Then, the solution to the Schrödinger equation \eqref{Schrodinger} satisfies the smoothing property
    \begin{align*}
        \||\cdot|^{\gamma-1}\,
        (-\Delta-|\rho|^{2})^{\frac{\gamma}{2}}\,u\|_{L^{2}(\mathbb{R}_{t}\times\mathbb{X})}\,
        \lesssim\,
        \|u_{0}\|_{L^{2}(\mathbb{X})}
        \qquad\forall\,u_{0}\in{L^2(\mathbb{X})}.
    \end{align*}
\end{corollary}

These examples of inequalities above are important tools for studying nonlinear partial differential equations. Since their sets of admissible indices are larger than those in the Euclidean context, we are able to obtain some stronger results. We refer to \cite{Ban07,AP09,IS09,AP14,AZ20,Zha21,AMPVZ23} where the better dispersive properties and Strichartz inequalities have been established in symmetric spaces, and the authors have used them to obtain well-posedness for a more general family of nonlinearities and some scattering results that are known to fail in the Euclidean setting.

\section{Preliminaries}\label{Section.2 Prelim}
We denote by $\mathbb{X}=G/K$ a Riemannian symmetric space of non-compact type, where $G$ is a connected, non-compact, semisimple Lie group with finite center and $K$ is a maximal compact subgroup of $G$. Let $\mathfrak{g}=\mathfrak{k}\oplus\mathfrak{p}$ be the Cartan decomposition of the Lie algebra $\mathfrak{g}$ of $G$. The Killing form of $\mathfrak{g}$ induces a $K$-invariant inner product $\langle\cdot\,,\,\cdot\rangle$ on $\mathfrak{p}$ and therefore a $G$-invariant Riemannian metric on $\mathbb{X}$. Denote by $\mathfrak{a}\subset\mathfrak{p}$ the Cartan subspace and $\mathfrak{a}^{+}\subset\mathfrak{a}$ the positive Weyl chamber. We say that $\mathbb{X}$ is of rank $\ell$ if $\mathfrak{a}$ has dimension $\ell$. We identify $\mathfrak{a}$ with its dual $\mathfrak{a}^{*}$ by means of the inner product inherited from $\mathfrak{p}$. Let $\Sigma\subset\mathfrak{a}$ be the root system of $( \mathfrak{g},\mathfrak{a})$ and $W$ be the Weyl group associated with $\Sigma$. We denote by $\Sigma^{+}$, $\Sigma_{r}^{+}$, and $\Sigma_{s}^{+}$ the set of positive, positive reduced, and simple roots. For each $\alpha\in\Sigma^{+}$, let $m_{\alpha}$ be the dimension of the positive root subspace $\mathfrak{g}_{\alpha}$. We denote by $\rho=\frac{1}{2}\sum_{\alpha\in\Sigma^{+}}m_{\alpha}\alpha$ the half sum of all positive roots counted with their multiplicities. Recall that the relation between the dimension and the pseudo-dimension (or dimension at infinity) 
    \begin{align}\label{Dimensions}
        n\,=\,
        \ell\,+\,\sum_{\alpha \in \Sigma^{+}}\,m_{\alpha}
        \quad\textnormal{and}\quad
        \nu\,=\,\ell\,+\,2|\Sigma_{r}^{+}|
    \end{align}
depends on the geometric structure of $\mathbb{X}$. For example, $\nu=3$ while $n\ge2$ is arbitrary in rank one and $\nu=n$ if $G$ is complex. In general, we know that $\ell\ge1$, $n\ge2$ and $\nu\ge3$.

For every $x\in{G}$, we denote by $A(x)$ its unique $\mathfrak{a}^{+}$-component in the Iwasawa decomposition $G=K(\exp\mathfrak{a}^{+})K$, and by $x^{+}$ its $\overline{\mathfrak{a}^{+}}$-component in the Cartan decomposition $G=K(\exp\overline{\mathfrak{a}^{+}})K$. Then we can write the Haar measure $\diff{x}$ on $G$ as
    \begin{align*}
        \int_{G}\diff{x}\,f(x)\,
        =\,\const\,\int_{K}\,\diff{k}_{1}\,
        \int_{\mathfrak{a}^{+}}\,\diff{x}^{+}\,\delta(x^{+})\,
        \int_{K}\,\diff{k}_{2}\,f(k_{1}(\exp x^{+})k_{2}),
    \end{align*}
where the density $\delta(x^{+})=\prod_{\alpha\in\Sigma^{+}}(\sinh\langle{\alpha,x^{+}}\rangle)^{m_{\alpha}}$ satisfies
     \begin{align}\label{density}
        \delta(x^{+})\,
        \asymp\,
            \prod_{\alpha\in\Sigma^{+}}
            \Big\lbrace 
            \frac{\langle\alpha,x^{+}\rangle}
            {1+\langle\alpha,x^{+}\rangle}
            \Big\rbrace^{m_{\alpha}}\,
            e^{\langle2\rho,x^{+}\rangle}\,
        \lesssim\,
            \begin{cases}
                |x^{+}|^{n-\ell}\,
                &\quad\textnormal{if}\;\;|x^{+}|\le\,1,\\[5pt]
                e^{\langle2\rho,x^{+}\rangle}\,
                &\quad\textnormal{for all}\;\;
                x^{+}\in\overline{\mathfrak{a}^{+}}.
            \end{cases}
    \end{align}

A function $f$ is called bi-$K$-invariant if $f(x)=f(\exp{x^{+}})$ for all $x\in{G}$. We denote by $\varphi_{\lambda}$ the \textit{elementary spherical function}, which is a bi-$K$-invariant eigenfunction for each $\lambda\in\mathfrak{a}$ and all $G$-invariant differential operators on $\mathbb{X}$. In the non-compact setting, we have the integral formula 
    \begin{align}
        \varphi_{\lambda}(x) 
        = \int_{K}\diff{k}\,
            e^{\langle{i\lambda+\rho,\,A(kx)}\rangle},
        \label{sphericalfct}
    \end{align}
see \cite[Theorem 4.3, p.418]{Hel00}. Recall that $\varphi_{\lambda}$ is symmetric in the sense that $\varphi_{\lambda}(x^{-1})=\varphi_{-\lambda}(x)$ for all $\lambda\in\mathfrak{a}$ and $x\in{G}$. For any bi-$K$-invariant function $f$ on $G$, one can write
\begin{align}
    \int_{G}\diff{y}\,
    f(y)\,\varphi_{\lambda}(y^{-1}x)\,
    &=\,
    \int_{G}\diff{y}\,
    f(y)\,\int_{K}\diff{k}\,
    \varphi_{\lambda}(y^{-1}kx)\notag\\[5pt]
    &=\,
    \varphi_{\lambda}(x)\,
    \int_{G}\diff{y}\,
    f(y)\,\varphi_{-\lambda}(y),
    \label{trick}
\end{align}
since $\varphi_{\lambda}$ is spherical and $\diff{k}$ is a normalized measure on $K$.

In our proofs, the Kunze-Stein phenomenon is a key tool. Let $f$ and $g$ be two integrable functions on $G$, the convolution product of $f$ and $g$ is defined by 
\begin{align*}
    (f*g)(x)\,
    =\,
    \int_{G}\diff{y}\,f(y)\,g(y^{-1}x).
\end{align*}
Recall first that the Young convolution inequality holds on the unimodular group $G$: for all $1\le{s,\,t_{1},\,t_{2}}\le+\infty$ such that $\frac{1}{t_{2}}+1=\frac{1}{s}+\frac{1}{t_{1}}$, we have 
\begin{align}
    L^{s}(G)*L^{t_{1}}(G)\,
    \subset\,L^{t_{2}}(G).
    \label{Young}    
\end{align}
Note that \eqref{Young} holds with $t_{1}=t_{2}=2$ if and only if $s=1$. In \cite{KS60}, Kunze and Stein discovered that, in the particular case where $G=\SL(2,\mathbb{R})$, the following holds
\begin{align}
    L^{s}(G)*L^{2}(G)\,
    \subset\,L^{2}(G)
    \qquad\forall\,1\le{s}<2,
    \label{KSL2}    
\end{align}
for all $1\le{s}<2$. This remarkable property has been extended by Cowling \cite{Cow78} to connected semisimple Lie groups with finite center. By interpolation between \eqref{KSL2} and $L^{\infty}(G)*L^{1}(G)\subset{L^{\infty}(G)}$, we deduce that $L^{t}(G)*L^{s}(G)\subset{L^{t}(G)}$, for $2\le{t}<+\infty$ and $1\le{s}<t'$. By duality, we obtain
    \begin{align}
        L^{t}(G)*L^{t'}(G)\subset{L^{s'}(G)},
        \qquad\forall\,1<s<t\le{2}.
        \label{KS1}
    \end{align}
Notice that one can permute $L^{t}(G)$ and $L^{t'}(G)$ in \eqref{KS1} since $G$ is unimodular. 

In the case where $g\in{L_{\textrm{loc}}^{1}(G)}$ is additionally bi-$K$ invariant, we know more precisely that
        \begin{align}
            \|f*|g|\|_{L^{s}(G)}\,
            =\,
            \|f\|_{L^{s}(G)}\,
            \int_{G}\diff{x}\,
            \varphi_{i(\frac{2}{s}-1)\rho}(x)\,|g(x)|,
            \label{KSbK1}
        \end{align}
for all $1\le{s}\le+\infty$, see \cite{Her70,Cow97}. Since the vector $-(\frac{2}{s}-1)\rho$ belongs to the Weyl chamber $\mathfrak{a}^{+}$ for $s>2$ and the set $\lbrace{\alpha\in\Sigma_{r}^{+}\,|\,\langle{\alpha,-(\frac{2}{s}-1)\rho}\rangle=0}\rbrace$ is empty, it is known that 
\begin{align}
    \varphi_{i(\frac{2}{s}-1)\rho}(\exp{x^{+}})\,
    \asymp\,
    e^{-\frac{2}{s}\langle{\rho,x^{+}}\rangle}\,
    \qquad\forall\,x^{+}\in\mathfrak{a}^{+},\,\,\,
    \forall\,s>2,
    \label{phirho}
\end{align}
according to \cite[Remark 3.1]{Sch08} and \cite[Lemma 3.4]{NPP14} (note the different notations of their spherical functions and \eqref{sphericalfct}). For $s=2$, an additional polynomial term will appear in the estimate of the ground spherical function:
\begin{align}\label{phi0}
    \varphi_{0}(\exp{x^{+}})\,
    \asymp\,\Big\lbrace{
    \prod_{\alpha\in\Sigma_{r}^{+}}(1+\langle{\alpha,x^{+}\rangle})\,
    }\Big\rbrace\,e^{-\langle{\rho,x^{+}}\rangle}
    \qquad\forall\,x^{+}\in\mathfrak{a}^{+},
\end{align}
see, for instance, \cite[Proposition 2.2.12]{AJ99}. By interpolation between \eqref{KSbK1} with $s=2$ and $L^{1}(G)*L^{\infty}(K\backslash{G/K})\subset{L^{\infty}(G)}$, we obtain the following estimate:
    \begin{align}
        \|f*|g|\|_{L^{s}(G)}\,
        \le\,
        \|f\|_{L^{s'}(G)}\,
        \Big\lbrace{
        \int_{G}\diff{x}\,\varphi_{0}(x)\,
        |g(x)|^{\frac{s}{2}}
        }\Big\rbrace^{\frac{2}{s}}
        \qquad\forall\,2\le{s}<+\infty.
        \label{KSbK2}
    \end{align}

\section{Proof of the main theorem}\label{Section Proof}
We give the proof of Theorem \ref{mainthm1} in this section. We start with the operator $T^{0}$. A crucial component of our analysis relies on Lemma \ref{SWLemma}, a Hardy-Littlewood-P\'olya type estimate. Moving on to the operator $T^{\infty}$, we need a delicate argument which combines the kernel estimate and the Kunze-Stein Phenomenon. The core of that part is a careful analysis for the critical indices. All necessity considerations are explored in the last subsection.

\subsection{$L^{p}$-$L^q$-boundedness of $T^{0}$}
Recall that the study of $T^{0}$ only involves operators $T_{1}^{0}$, $T_{2}^{0}$ and $T_{3}^{0}$, where
\begin{align*}
    T_{j}^{0}(\zeta,\sigma,\alpha,\beta)f(x)\,=\,
    \int_{\mathbb{X}}\diff{y}\,
    \mathds{1}_{Z_{j}}(x,y)\,
    \psi_{0}(y^{-1}x)\,
    |x|^{-\beta}\,k_{\zeta,\sigma}(y^{-1}x)\,
    |y|^{-\alpha}\,f(y)
    \qquad\forall\,j\,=\,1,2,3.
\end{align*}
Here, the zones $Z_{j}$ are defined in Sect.\ref{Sect Strategy} and can be visualized in Figure \ref{Regions}.  We initiate our discussion with the operator $T_{1}^{0}$. To establish its $L^{p}$-$L^q$-boundedness, we apply a dyadic decomposition technique. Recall that $\psi_{0}\in\mathcal{C}_{c}^{\infty}(K\backslash{G/K})$ is a bi-$K$-invariant function on $\mathbb{X}$ such that its support is contained in the set $\lbrace{x\in\mathbb{X}\,|\,|x|\le1}\rbrace$ and $\psi_{0}(x)=1$ if $|x|\le\frac12$. We define 
\begin{align*}
        \Psi_{k}(x)\,
        =\,
        \psi_{0}(2^{-k}x)-\psi_{0}(2^{-k+1}x)
        \qquad\forall\,k\in\mathbb{Z},\,\,\,
        \forall\,x\in\mathbb{X},
\end{align*}
which are all bi-$K$-invariant cut-off functions on $\mathbb{X}$. For each integer $k$, the function $\Psi_{k}$ is compactly supported in $\lbrace{x\in\mathbb{X}\,|\,2^{k-2}\le|x|\le2^{k}}\rbrace$. In particular, we have 
\begin{align}\label{dyadicZ}
    \sum_{k\in\mathbb{Z}}\Psi_{k}(x)\,=\,1
    \qquad\forall\,x\in\mathbb{X},
\end{align}
and
\begin{align*}
    \|f\|_{L^{p}(\mathbb{X})}^{p}\,
    \asymp\,
    \sum_{k\in\mathbb{Z}}
    \|\Psi_{k}f\|_{L^{p}(\mathbb{X})}^{p}\,
    \qquad\forall\,f\in{L^{p}(\mathbb{X})},
\end{align*}
for all $1\le{p}<+\infty$. Now we prove the following result concerning the $L^p$-$L^q$-boundedness of $T_{1}^{0}$.

\begin{proposition}\label{PropT01}
Suppose that $\zeta\ge0$, $\sigma>0$ ($0<\sigma<\nu$ if $\zeta=0$), $\alpha\in\mathbb{R}$, $\beta\in\mathbb{R}$, and $0\le\alpha+\beta\le\sigma$. Then the operator $T_{1}^{0}(\zeta,\sigma,\alpha,\beta)$ is bounded from $L^{p}(\mathbb{X})$ to $L^{q}(\mathbb{X})$ for all $1\le{p}\le{q}\le{\infty}$ satisfying $\frac{1}{p'}>\frac{\alpha}{n}$ when $\alpha>0$, $\frac{1}{q}>\frac{\beta}{n}$ when $\beta>0$, $\frac{1}{p}-\frac{1}{q}\le\frac{\sigma-\alpha-\beta}{n}$ when $\sigma-\alpha-\beta<n$, and $(\frac{1}{p},\frac{1}{q})\neq(\frac{\sigma-\alpha-\beta}{n},0)$ or $(1,1-\frac{\sigma-\alpha-\beta}{n})$.
\end{proposition}

\begin{proof}
First, note that this proposition is true for $\alpha=\beta=0$, according to \cite[Corollary 4.2.(i)]{Ank92}. We start with the case where $0\le\alpha+\beta<\sigma$. On the one hand, if $\sigma<n$, we have
    \begin{align*}
        T_{1}^{0}(\zeta,\sigma,\alpha,\beta)f(x)\,
        \lesssim\,
        \int_{\mathbb{X}}\diff{y}\,
        \mathds{1}_{Z_{1}}(x,y)\,
        \psi_{0}(y^{-1}x)\,
        |y^{-1}x|^{\sigma-\alpha-\beta-n}\,f(y),
    \end{align*}
according to the kernel estimate \eqref{estim k0} and the fact that $|x|^{-\beta}|y|^{-\alpha}\lesssim|y^{-1}x|^{-\alpha-\beta}$, since $|y|\asymp|x|$ for $(x,y)\in{Z_{1}}$ and $\alpha+\beta\ge0$. Therefore, we have $|T_{1}^{0}(\zeta,\sigma,\alpha,\beta)f|\lesssim{|T_{1}^{0}(\zeta,\sigma-\alpha-\beta,0,0)f|}$ where $0<\sigma-\alpha-\beta\le\sigma<n$ as assumed. It is already known that the operator $T_{1}^{0}(\zeta,\sigma-\alpha-\beta,0,0)$ is bounded from $L^{p}(\mathbb{X})$ to $L^{q}(\mathbb{X})$ for all $1\le{p}\le{q}\le{\infty}$ satisfying $\frac{1}{p}-\frac{1}{q}\le\frac{\sigma-\alpha-\beta}{n}$, except for the points $(\frac{\sigma-\alpha-\beta}{n},0)$ and $(1,1-\frac{\sigma-\alpha-\beta}{n})$.

On the other hand, note that $\frac{1}{p'}>\frac{\alpha}{n}$ and $\frac{1}{q}>\frac{\beta}{n}$ imply that there exists $\varepsilon>0$ such that $\alpha+\beta\le{n(\frac{1}{p'}+\frac{1}{q}-\varepsilon)}$. According to the kernel estimate \eqref{estim k0} with $\sigma\ge{n}$, we have, for all $(x,y)\in{Z_1}$,
\begin{align*}
    |x|^{-\beta}\,\psi_{0}(y^{-1}x)\,k_{\zeta,\sigma}(y^{-1}x)\,|y|^{-\alpha}\,
    &\lesssim\,
    |x|^{-\alpha-\beta}\,\psi_{0}(y^{-1}x)\,|y^{-1}x|^{-\frac{\varepsilon n}{2}}\\[5pt]
    &\lesssim\,
    \psi_{0}(y^{-1}x)\,|y^{-1}x|^{-n(\frac{1}{p'}+\frac{1}{q}-\frac{\varepsilon}{2})}\\[5pt]
    &=\,
    \psi_{0}(y^{-1}x)\,|y^{-1}x|^{n(\frac{1}{p}-\frac{1}{q}+\frac{\varepsilon}{2})\,-\,n},
\end{align*}
where $0<\frac{1}{p}-\frac{1}{q}+\frac{\varepsilon}{2}<1$. Hence, we obtain
\begin{align*}
    \|T_{1}^{0}(\zeta,\sigma,\alpha,\beta)f\|_{L^q(\mathbb{X})}\,
    \lesssim\,
    \Big\|T_{1}^{0}\Big(\zeta,n\Big(\frac{1}{p}-\frac{1}{q}+\frac{\varepsilon}{2}\Big),0,0\Big)f\Big\|_{L^q(\mathbb{X})}\,
    \lesssim\,
    \|f\|_{L^p(\mathbb{X})},
\end{align*}
for all $f\in{L^{p}(\mathbb{X})}$. Therefore, the operator $T_{1}^{0}(\zeta,\sigma,\alpha,\beta)$ is $L^p$-$L^q$-bounded for all $\zeta\ge0$ and $\sigma>0$ in the case where $0\le\alpha+\beta<\sigma$.

When $\alpha+\beta=\sigma$, we observe, together with conditions $1\le{p}\le{q}\le{\infty}$ and $\frac{1}{p}-\frac{1}{q}\le\frac{\sigma-\alpha-\beta}{n}$, that $1\le{p=q}\le+\infty$. Moreover, since $\frac{1}{p'}>\frac{\alpha}{n}$ and $\frac{1}{q}>\frac{\beta}{n}$, we know that $0<\sigma=\alpha+\beta<n$. Using again the kernel estimate \eqref{estim k0} and the dyadic decomposition \eqref{dyadicZ}, we have, for all $1\le{p}<+\infty$,
\begin{align*}
    \|T_{1}^{0}(\zeta,\sigma,\alpha,\beta)
    &f\|_{L^{p}(\mathbb{X})}^{p}\\[5pt]
    &\lesssim\,
    \sum_{k\in\mathbb{Z}}\,
    \int_{\mathbb{X}}\diff{x}\,
    \Big|
    \int_{\frac{|x|}{2}\le|y|\le2|x|}\diff{y}\,\Psi_{k}(x)\,
    |x|^{-\beta}\,\psi_{0}(y^{-1}x)\,
    |y^{-1}x|^{\sigma-n}\,|y|^{-\alpha}\,f(y)
    \Big|^{p}.
\end{align*}
Notice that $2^{k-2}\le|x|\le2^{k}$ and $\frac{|x|}{2}\le|y|\le2|x|$ imply $2^{k-3}\le|y|\le2^{k+2}$ and $|y^{-1}x|\le2^{k+2}$. By using the Young convolution inequality \eqref{Young}, we obtain
\begin{align*}
    \|T_{1}^{0}(\zeta,\sigma,\alpha,\beta)f\|_{L^{p}(\mathbb{X})}^{p}\,
    &\lesssim\,
    \sum_{k\in\mathbb{Z}}\,2^{-p\sigma{k}}\,
    \int_{\mathbb{X}}\diff{x}\,
    \Big|(\Psi_{k}f)\,
    *\,
    \Big(|\cdot|^{\sigma-n}\,
    \chi_{0}\Big(\frac{|\,\cdot\,|}{2^{k+3}}\Big)\,\psi_{0}\Big)(x)\Big|^{p}\\[5pt]
    &\le\,
    \sum_{k\in\mathbb{Z}}\,2^{-p\sigma{k}}\,
    \|\Psi_{k}f\|_{L^{p}(\mathbb{X})}^{^p}\,
    \Big\||\cdot|^{\sigma-n}\,
    \chi_{0}\Big(\frac{|\,\cdot\,|}{2^{k+3}}\Big)\,\psi_{0}\Big\|_{L^{1}(\mathbb{X})}^{p},
\end{align*}
where
\begin{align*}
    \Big\||\cdot|^{\sigma-n}\,
    \chi_{0}\Big(\frac{|\,\cdot\,|}{2^{k+3}}\Big)\,\psi_{0}\Big\|_{L^{1}(\mathbb{X})}\,
    &\le\,
    \int_{\mathfrak{a}^{+}}\diff{x^{+}}\,
    \chi_{0}(|x^{+}|)\,|x^{+}|^{\sigma-n}|x^{+}|^{n-\ell}\\[5pt]
    &\le\,
    \int_{0}^{1}\diff{r}\,r^{\sigma-1}\,=\,\sigma
\end{align*}
if $k>-3$, and 
    \begin{align*}
        \Big\||\cdot|^{\sigma-n}\,
        \chi_{0}\Big(\frac{|\,\cdot\,|}{2^{k+3}}\Big)\,\psi_{0}\Big\|_{L^{1}(\mathbb{X})}\,
        &\le\,
        \int_{\mathfrak{a}^{+}}\diff{x^{+}}\,
        \chi_{0}\Big(\frac{|x^{+}|}{2^{k+3}}\Big)\,
        |x^{+}|^{\sigma-n}|x^{+}|^{n-\ell}\\[5pt]
        &\le\,
        \int_{0}^{2^{k+3}}\diff{r}\,r^{\sigma-1}\,
        =\,\sigma\,2^{(k+3)\sigma}
    \end{align*}
for all $k\le-3$. Therefore, for all $1\le{p}<+\infty$,
\begin{align*}
    \|T_{1}^{0}(\zeta,\sigma,\alpha,\beta)f\|_{L^{p}(\mathbb{X})}^{p}\,
    \lesssim\,
    \sum_{k\in\mathbb{Z}}\,
    \|\Psi_{k}f\|_{L^{p}(\mathbb{X})}^{^p}\,
    \asymp\,
    \|f\|_{L^{p}(\mathbb{X})}^{^p},
\end{align*}
in the case where $\sigma=\alpha+\beta$. The sup norm case can be handled in the same way and the proof is complete.
\end{proof}

Now let us turn to operators $T_{2}^{0}$ and $T_{3}^{0}$. For this part, we will extend the argument carried out in \cite{SW58} from the Euclidean setting to non-compact symmetric spaces. The main step is the following Hardy-Littlewood-P\'olya type estimate. Compared to similar results in the Euclidean setting, e.g. in \cite{SW58,Str69,Wal71}, it is noteworthy that the condition \eqref{kappa12} in the following lemma arises naturally due to the particular geometry at infinity of symmetric spaces.

\begin{lemma}\label{SWLemma}
    Let $\mathbb{X}$ be a non-compact symmetric space of rank $\ell\ge1$ and $\mathcal{K}:\mathbb{R}_{+}\times\mathbb{R}_{+}\rightarrow\mathbb{R}_{+}$ a homogeneous function of degree $-\ell$ such that, for all $1\le{p}\le+\infty$,
        \begin{align*}
            C_{\mathcal{K}}\,
            =\,
            \int_{0}^{+\infty}\diff{s}\,s^{\frac{\ell}{p'}-1}\,\mathcal{K}(1,s)\,
            <\,+\infty.
        \end{align*}
    Let $\kappa_{1}$ and $\kappa_{2}$ be two bi-$K$-invariant functions on $G$ such that, for all $x\in{G}$,
        \begin{align}
            |\kappa_{1}(\exp{x}^{+})|\,\delta^{\frac{1}{p}}(x^{+})\,
            \le\,C_{\kappa_{1}}\,
            \qquad\textnormal{and}\qquad
            |\kappa_{2}(\exp{x}^{+})|\,\delta^{\frac{1}{p'}}(x^{+})\,
            \le\,C_{\kappa_{2}}
            \label{kappa12}
        \end{align}
    where $x^+$ is the $\overline{\mathfrak{a}^{+}}$-component in the Cartan decomposition of $x$ and the density $\delta(x^{+})$ satisfies \eqref{density}. Then, the operator $S:L^{p}(G)\rightarrow{L^{p}(G)}$ defined by
        \begin{align*}
            Sf(x)\,
            =\,\kappa_{1}(x)\,
                \int_{G}\diff{y}\,
                \mathcal{K}(|x|,|y|)\,
                \kappa_{2}(y)\,f(y)
        \end{align*}
    is bounded for all $1\le{p}\le+\infty$.
    \end{lemma}
    
\begin{proof}
Since $\kappa_{1}$ and $\kappa_{2}$ are bi-$K$-invariant functions on $G$, there exists a constant $C_0>0$ such that 
        \begin{align*}
            Sf(x)\,
            =\,C_0\,\kappa_{1}(\exp{x^{+}})\,
                \int_{K}\diff{k_{1}}\,\int_{K}\diff{k_{2}}\,
                \underbrace{
                \int_{\mathfrak{a}^{+}}\diff{y^{+}}\,\delta(y^{+})\,
                    \mathcal{K}(|x^{+}|,|y^{+}|)\,
                    \kappa_{2}(\exp{y^{+}})\,f(k_{1}(\exp{y^{+}})k_{2})
                }_{=\, I_{k_{1},k_{2}}(|x^{+}|)},
        \end{align*}
    according to the Cartan decomposition. Assume that $\ell\ge2$. Using polar coordinates $x^{+}=R\xi$ and $y^{+}=r\eta$ with $R,r>0$ and $\xi,\eta\in\mathbb{S}^{\ell-1}$, we write
        \begin{align}
            I_{k_{1},k_{2}}(R)\,
            =\,
            \int_{\mathbb{S}^{\ell-1}}\diff{\sigma_{\eta}}\,
                    \int_{0}^{+\infty}\diff{r}\,r^{\ell-1}\,
                    \delta(r\eta)\,\mathcal{K}(R,r)\,
                    \kappa_{2}(\exp(r\eta))\,f(k_{1}(\exp{(r\eta))k_{2}}).
            \label{polarI}
        \end{align}
    Next, substituting $r=sR$ and using the homogeneity of $\mathcal{K}$, we obtain, for all $k_{1},k_{2}\in{K}$,
        \begin{align}\label{SW jensen}
           I_{k_{1},k_{2}}(R)\,
            =\,\int_{\mathbb{S}^{\ell-1}}\diff{\sigma_{\eta}}
                \underbrace{
                \int_{0}^{+\infty}\diff{s}\,s^{\ell-1}\,
                    \delta(sR\eta)\,\mathcal{K}(1,s)\,
                    \kappa_{2}(\exp(sR\eta))\,f(k_{1}(\exp{(sR\eta))k_{2}})
                }_{=\,\widetilde{I}_{\eta,k_{1},k_{2}}(R)}.
        \end{align}
    
    In the case where $p=\infty$ and $p'=1$, we deduce from \eqref{SW jensen} and \eqref{kappa12} that 
    \begin{align*}
        |Sf(x)|\,
        \le\,
        C_0\,C_{\kappa_{1}}\,C_{\kappa_{2}}\,|\omega_{\ell-1}|\,
        \|f\|_{L^{\infty}(\mathbb{X})}\,
        \underbrace{\int_{0}^{+\infty}\diff{s}\,s^{\ell-1}\,\mathcal{K}(1,s)
        }_{=\, C_{\mathcal{K}}}\,
        \lesssim\,
        \|f\|_{L^{\infty}(\mathbb{X})},
    \end{align*}
    since $\diff{k_{1}}$ and $\diff{k_{2}}$ are normalized measures in $K$. When $1\le{p}<+\infty$, we know from the duality that there exists a function $h$ on $\mathbb{R}_{+}$ such that $\int_{0}^{+\infty}\diff{R}R^{\ell-1}|h(R)|^{p'}=1$ and 
        \begin{align}
            \Big\lbrace
            \int_{0}^{+\infty}\diff{R}\,R^{\ell-1}\,
            |\widetilde{I}_{\eta,k_{1},k_{2}}(R)|^{p}
            \Big\rbrace^{\frac{1}{p}}\,
            =\,
            \int_{0}^{+\infty}\diff{R}\,
            R^{\ell-1}\,h(R)\,\widetilde{I}_{\eta,k_{1},k_{2}}(R).
            \label{aux h}
        \end{align}
    Using successively \eqref{aux h}, the Fubini theorem, and the Hölder inequality, we have
        \begin{align*}
            &\Big\lbrace
            \int_{0}^{+\infty}\diff{R}\,R^{\ell-1}\,
            |\widetilde{I}_{\eta,k_{1},k_{2}}(R)|^{p}
            \Big\rbrace^{\frac{1}{p}}\\[5pt]
            &=\,
            \int_{0}^{+\infty}\diff{s}\,
                s^{\ell-1}\,\mathcal{K}(1,s)\,
                \int_{0}^{+\infty}\diff{R}\,R^{\ell-1}\,
                h(R)\,\delta(sR\eta)\,
                \kappa_{2}(\exp(sR\eta))\,
                f(k_{1}(\exp{(sR\eta))k_{2}})\\[5pt]
            &\le\,
            \int_{0}^{+\infty}\diff{s}\,
                s^{\ell-1}\,\mathcal{K}(1,s)\,
                \Big\lbrace{
                \int_{0}^{+\infty}\diff{R}\,R^{\ell-1}\,
                \delta(sR\eta)\,
                |f(k_{1}(\exp{(sR\eta))k_{2}})|^{p}
                }\Big\rbrace^{\frac{1}{p}}\\[5pt]
                &\hspace{106pt}
                \underbrace{\Big\lbrace{
                \int_{0}^{+\infty}\diff{R}\,R^{\ell-1}\,
                \delta(sR\eta)\,
                |\kappa_{2}(\exp(sR\eta))|^{p'}\,
                |h(R)|^{p'}\,
                }\Big\rbrace^{\frac{1}{p'}}}_{\le\,C_{\kappa_{2}}},
        \end{align*}
    due to \eqref{kappa12} and the definition of $h$. Substituting $R=rs^{-1}$ in the rest inner integral, we obtain 
    \begin{align}
            \Big\lbrace
            \int_{0}^{+\infty}\diff{R}\,R^{\ell-1}\,
            |\widetilde{I}_{\eta,k_{1},k_{2}}(R)|^{p}
            \Big\rbrace^{\frac{1}{p}}
            \le\,
            C_{\mathcal{K}}\,C_{\kappa_{2}}\,
            \Big\lbrace{
                \int_{0}^{+\infty}\diff{r}\,r^{\ell-1}\,
                \delta(r\eta)\,
                |f(k_{1}(\exp{(r\eta))k_{2}})|^{p}
                }\Big\rbrace^{1/p}.
            \label{aux h2}
        \end{align} 
    Applying Jensen's inequality to \eqref{SW jensen}, we deduce from \eqref{aux h2} that
        \begin{align}
        \int_{0}^{+\infty}\diff{R}\,R^{\ell-1}\,|I_{k_{1},k_{2}}(R)|^{p}\,
        &\le\,|\omega_{\ell-1}|^{p-1}\,
            \int_{0}^{+\infty}\diff{R}\,R^{\ell-1}\,
                \int_{\mathbb{S}^{\ell-1}}\,\diff{\sigma_{\eta}}\,
                |\widetilde{I}_{\eta,k_{1},k_{2}}(R)|^{p}\notag\\[5pt]
        &\le\,
                C_{\mathcal{K}}^{p}\,C_{\kappa_{2}}^{p}\,
                |\omega_{\ell-1}|^{p-1}\,
                \int_{\mathbb{S}^{\ell-1}}\,\diff{\sigma_{\eta}}\,
                \int_{0}^{+\infty}\diff{r}\,r^{\ell-1}\delta(r\eta)\,
                |f(k_{1}(\exp(r\eta))k_{2})|^{p}\notag\\[5pt]
        &=\, 
            C_{\mathcal{K}}^{p}\,C_{\kappa_{2}}^{p}\,
            |\omega_{\ell-1}|^{p-1}\,
            \int_{\mathfrak{a}^{+}}\diff{y^{+}}\,
                \delta(y^{+})\,|f(k_{1}(\exp{y^{+}})k_{2})|^{p}.
            \label{SW jensen2}
        \end{align}
    
    Using successively the Cartan decomposition, the assumption \eqref{kappa12}, the bi-$K$-invariance of $I_{k_{1},k_{2}}$, and the estimate \eqref{SW jensen2}, we finally obtain, for all $1\le{p}<+\infty$ and $\ell\ge2$,
        \begin{align*}
        \|Sf\|_{L^{p}(\mathbb{X})}^{p}\,
        &=\,C_0\,
            \int_{\mathfrak{a}^{+}}\diff{x^{+}}\,
            \delta(x^{+})\,
            |Sf(\exp{x^{+}})|^{p}\\[5pt]
        &\le\,C_0^{p+1}\,  
            \int_{\mathfrak{a}^{+}}\diff{x^{+}}\,
            \delta(x^{+})\,
            |\kappa_{1}(\exp{x^{+}})|^{p}\,
                \int_{K}\diff{k_{1}}\,\int_{K}\diff{k_{2}}\,
                |I_{k_{1},k_{2}}(|x^{+}|)|^{p}\\[5pt]
        &\le\,C_{\kappa_{1}}^{p}\,C_0^{p+1}\,|\omega_{\ell-1}|\,
            \int_{K}\diff{k_{1}}\,\int_{K}\diff{k_{2}}\,
            \int_{0}^{+\infty}\,\diff{R}\,R^{\ell-1}\,
            |I_{k_{1},k_{2}}(R)|^{p}\\[5pt]
        &\le\,\widetilde{C}^{p}\,C_0\,
            \int_{K}\diff{k_{1}}\,\int_{K}\diff{k_{2}}\,
            \int_{\mathfrak{a}^{+}}\diff{y^{+}}\,
                \delta(y^{+})\,|f(k_{1}(\exp{y^{+}})k_{2})|^{p}\\[5pt]
        &=\,\widetilde{C}^{p}\,
            \|f\|_{L^{p}(\mathbb{X})}^{p},
        \end{align*}
    where $\widetilde{C}=C_{0}C_{\mathcal{K}}C_{\kappa_{1}}C_{\kappa_{2}}|\omega_{\ell-1}|$. When $\mathbb{X}$ is of rank $\ell=1$, we skip step \eqref{polarI} and notice that $I_{k_{1},k_{2}}(R)=\widetilde{I}_{\eta,k_{1},k_{2}}(R)$ for all $k_{1},k_{2}\in{K}$ in step \eqref{SW jensen}, then conclude in the same way.
    \end{proof}

Based on this lemma, we prove the following proposition regarding the $L^p$-$L^q$-boundedness of the operators $T_{2}^{0}$ and $T_{3}^{0}$.

\begin{proposition}\label{PropT023}
    Suppose that $\zeta\ge0$, $\sigma>0$ ($0<\sigma<\nu$ if $\zeta=0$), $\alpha\in\mathbb{R}$, $\beta\in\mathbb{R}$, and $0\le\alpha+\beta\le\sigma$. Then, the operators $T_{2}^{0}(\zeta,\sigma,\alpha,\beta)$ and $T_{3}^{0}(\zeta,\sigma,\alpha,\beta)$ are bounded from $L^{p}(\mathbb{X})$ to $L^{q}(\mathbb{X})$ for all $1\le{p}\le{q}\le{\infty}$ satisfying $\frac{1}{p'}>\frac{\alpha}{n}$ when $\alpha>0$, $\frac{1}{q}>\frac{\beta}{n}$ when $\beta>0$, and $\frac{1}{p}-\frac{1}{q}\le\frac{\sigma-\alpha-\beta}{n}$ when $\sigma-\alpha-\beta<n$.
\end{proposition}

\begin{proof}
We will prove that $T_{2}^{0}$ is bounded from $L^p(\mathbb{X})$ to $L^q(\mathbb{X})$ for suitable $p$ and $q$. Then the $L^p$-$L^q$-boundedness of $T_{3}^{0}$ follows similarly by symmetry, see Figure \ref{Regions}. By duality of Lebesgue spaces, we know that the operator $T_{2}^{0}$ is bounded from $L^p(\mathbb{X})$ to $L^q(\mathbb{X})$ if and only if the double integral $\mathcal{I}_{2}^{0}$ defined by
    \begin{align*}
        \mathcal{I}_{2}^{0}\,
        =\,
        \int\int_{Z_{2}}\diff{x}\diff{y}\,
        |x|^{-\beta}\,(\psi_{0}k_{\zeta,\sigma})(y^{-1}x)\,
        |y|^{-\alpha}\,f(y)\,g(x)
    \end{align*}
satisfies $|\mathcal{I}_{2}^{0}|\lesssim\|f\|_{L^{p}(\mathbb{X})}\|g\|_{L^{q'}(\mathbb{X})}$ for all $f\in{L^p(\mathbb{X})}$ and $g\in{L^{q'}(\mathbb{X})}$. Recall that the kernel $\psi_{0}k_{\zeta,\sigma}$ behaves differently depending on $\sigma$, see \eqref{estim k0}. We will first assume that $0<\sigma<{n}$ and deal with the case $\sigma\ge{n}$ at the end of the proof. For all $(x,y)$ in ${Z_{2}}$, we have $2|x|\le|y|\le1$, $|y^{-1}x|\ge\frac{|y|}{2}$, and then $|y^{-1}x|^{\sigma-n}\lesssim|y|^{\sigma-n}$. Hence,
    \begin{align}
        |\mathcal{I}_{2}^{0}|\,
        &\lesssim\,
        \int\int_{Z_{2}}\diff{x}\diff{y}\,
        |x|^{-\beta}\,|y|^{\sigma-n-\alpha}\,
        \psi_{0}(y^{-1}x)\,f(y)\,g(x)\notag\\[5pt]
        &\lesssim\,
        \int_{|y|\le1}\diff{y}\,f(y)\,
        |y|^{\sigma-\alpha-\beta}\,
        \underbrace{
        |y|^{-n+\beta}\,\int_{|x|\le\frac{|y|}{2}}\diff{x}\,
        |x|^{-\beta}\,g(x)
        }_{:=\,(V_{\beta}g)(y)}.
        \label{V beta}
    \end{align}
We claim that, for all $1\le{q}\le+\infty$ such that $\frac{1}{q}>\frac{\beta}{n}$ when $\beta>0$, the following two estimates hold:
    \begin{align}
        |\chi_{0}(y)V_{\beta}g(y)|\,
        \lesssim\,|y|^{-\frac{n}{q'}}\,
        \|g\|_{L^{q'}(\mathbb{X})},
        \label{V1}
    \end{align}
    and
    \begin{align}
        \|\chi_{0}V_{\beta}g\|_{L^{q'}(\mathbb{X})}\,
        \lesssim\,
        \|g\|_{L^{q'}(\mathbb{X})}.
        \label{V2}
    \end{align}

We deduce from \eqref{V beta} and the Hölder inequality that
    \begin{align*}
        |\mathcal{I}_{2}^{0}|\,
        \lesssim\,
        \|f\|_{L^{p}(\mathbb{X})}\,
        \||\cdot|^{n(\frac{1}{p}-\frac{1}{q})}\,\chi_{0}V_{\beta}g
        \|_{L^{p'}(\mathbb{X})},
    \end{align*}
provided that $\frac{1}{p}-\frac{1}{q}\le\frac{\sigma-\alpha-\beta}{n}$. It is clear that, when $1\le{p=q}\le+\infty$, we have
    \begin{align*}
        |\mathcal{I}_{2}^{0}|\,
        \lesssim\,
        \|f\|_{L^{p}(\mathbb{X})}\,
        \|\chi_{0}V_{\beta}g
        \|_{L^{p'}(\mathbb{X})}\,
        \lesssim\,
         \|f\|_{L^{p}(\mathbb{X})}\,
        \|g\|_{L^{q'}(\mathbb{X})},
    \end{align*}
according to \eqref{V2}. In the case where $p\neq{q}$ (then $q'\neq\infty$) and $p>1$, we deduce from \eqref{V1} and \eqref{V2} that
    \begin{align*}
        \||\,.\,|^{n(\frac{1}{p}-\frac{1}{q})}\,\chi_{0}V_{\beta}g
        \|_{L^{p'}(\mathbb{X})}\,
        &=\,
        \Big\lbrace{
        \int_{\mathbb{X}}\diff{y}\,
        |y|^{n(\frac{1}{p}-\frac{1}{q})p'}\,
        |(\chi_{0}V_{\beta}g)(y)|^{p'-q'}\,
        |(\chi_{0}V_{\beta}g)(y)|^{q'}
        }\Big\rbrace^{\frac{1}{p'}}\\[5pt]
        &\lesssim\,
        \Big\lbrace{
        \int_{\mathbb{X}}\diff{y}\,
        |y|^{n(\frac{1}{p}-\frac{1}{q})p'}\,
        |y|^{-\frac{n(p'-q')}{q'}}\,
        |(\chi_{0}V_{\beta}g)(y)|^{q'}
        }\Big\rbrace^{\frac{1}{p'}}\,
        \|g\|_{L^{q'}(\mathbb{X})}^{\frac{p'-q'}{p'}}\\[5pt]
        &\lesssim\,
        \|g\|_{L^{q'}(\mathbb{X})}^{\frac{q'}{p'}}\,
        \|g\|_{L^{q'}(\mathbb{X})}^{\frac{p'-q'}{p'}}\,
        =\,\|g\|_{L^{q'}(\mathbb{X})},
    \end{align*}
since $n(\frac{1}{p}-\frac{1}{q})p'-\frac{n(p'-q')}{q'}=0$. If $p\neq{q}$ and $p=1$, the sup norm estimate follows directly from \eqref{V1}:
\begin{align*}
    \big||y|^{n(\frac{1}{p}-\frac{1}{q})}\,(\chi_{0}V_{\beta})g(y)\big|\,
    \lesssim\,
    |y|^{n(1-\frac{1}{q})}\,|y|^{-\frac{n}{q'}}\,
    \|g\|_{L^{q'}(\mathbb{X})}\,
    =\,\|g\|_{L^{q'}(\mathbb{X})}.
\end{align*}
Therefore, for all $1\le{p}\le{q}\le{\infty}$ such that $\frac{1}{q}>\frac{\beta}{n}$ when $\beta>0$ and $\frac{1}{p}-\frac{1}{q}\le\frac{\sigma-\alpha-\beta}{n}$, we have $|\mathcal{I}_{2}^{0}|\lesssim\|f\|_{L^{p}(\mathbb{X})}\|g\|_{L^{q'}(\mathbb{X})}$, that is, the operator $T_{2}^{0}(\zeta,\sigma,\alpha,\beta)$ is bounded from $L^{p}(\mathbb{X})$ to $L^{q}(\mathbb{X})$, when $0\le\alpha+\beta\le\sigma<n$.

Now, we prove the claims \eqref{V1} and \eqref{V2}. Using the Hölder inequality, we have straightforwardly
    \begin{align*}
        |\chi_{0}(y)\,V_{\beta}g(y)|\,
        &\lesssim\,
        |y|^{-n+\beta}\,
        \Big\lbrace{
        \int_{|x|\le\frac{|y|}{2}}\diff{x}\,
        |x|^{-\beta{q}}
        }\Big\rbrace^{\frac{1}{q}}\,
        \|g\|_{L^{q'}(\mathbb{X})}\\[5pt]
        &\lesssim\,
        |y|^{-n+\beta}\,
        \Big\lbrace{
        \int_{0}^{\frac{|y|}{2}}\diff{r}\,
        r^{-\beta{q}+n-1}
        }\Big\rbrace^{\frac{1}{q}}\,
        \|g\|_{L^{q'}(\mathbb{X})}\,
        \lesssim\,
        |y|^{-\frac{n}{q'}}\,
        \|g\|_{L^{q'}(\mathbb{X})},
    \end{align*}
for all $1<q<+\infty$, provided that $\beta<\frac{n}{q}$. If $q=\infty$, the condition $\beta<\frac{n}{q}$ implies that $\beta<0$. Therefore, the estimate \eqref{V1} follows immediately. When $q=1$, \eqref{V1} is equivalent to \eqref{V2}, which is a consequence of Lemma \ref{SWLemma}. More precisely,
    \begin{align*}
        |(\chi_{0}V_{\beta}g)(y)|\,
        &\lesssim\,
        \underbrace{\vphantom{\int}
        \chi_{0}(|y|)\,|y|^{\frac{\ell-n}{q'}}
        }_{\kappa_{1}(y)}\,
        \int_{\mathbb{X}}\diff{x}\,
        \underbrace{\vphantom{\int}
        \chi_{0}\Big(\frac{|x|}{|y|}\Big)
        |y|^{-n+\beta+\frac{n-\ell}{q'}}\,
        |x|^{-\beta+\frac{n-\ell}{q}}}_{\mathcal{K}(|y|,|x|)}
        \underbrace{\vphantom{\int}
        \chi_{0}(|x|)\,|x|^{\frac{\ell-n}{q}}
        }_{\kappa_{2}(x)}\,|g(x)|
    \end{align*}
where $(-n+\beta+\frac{n-\ell}{q'})+(-\beta+\frac{n-\ell}{q})=-\ell$,
    \begin{align*}
        \kappa_{1}(\exp{y^{+}})\delta^{\frac{1}{q'}}(y^{+})=\mathrm{O}(1)
        \qquad\textnormal{and}\qquad
        \kappa_{2}(\exp{x^{+}})\delta^{\frac{1}{q}}(x^{+})=\mathrm{O}(1).
    \end{align*}
We deduce \eqref{V2} from Lemma \ref{SWLemma} by noticing that
    \begin{align*}
        \int_{0}^{+\infty}\diff{s}\,
        s^{\frac{\ell}{q}-1}\,\mathcal{K}(1,s)\,
        =\,
        \int_{0}^{1}\diff{s}\,
        s^{\frac{\ell}{q}-1}\,s^{-\beta+\frac{n-\ell}{q}}\,
        <\,+\infty,
    \end{align*}
provided that $\beta<\frac{n}{q}$ when $\beta>0$. Then we complete the proof for the $L^{p}$-$L^{q}$-boundedness of the operator $T_{2}^{0}(\zeta,\sigma,\alpha,\beta)$ in the case where $0<\sigma<n$.

Using a standard modification based on the above arguments, we can handle the case where $\sigma\ge{n}$. Note that conditions $\frac{1}{p'}>\frac{\alpha}{n}$, $\frac{1}{q}>\frac{\beta}{n}$, and $\sigma\ge{n}$ imply that
\begin{align*}
    n\,\Big(\frac{1}{p}-\frac{1}{q}\Big)\,
    <\,n-\alpha-\beta\,\le\,\sigma-\alpha-\beta.
\end{align*}
Let $\varepsilon_1>0$ be a constant such that $\sigma-\varepsilon_1=n$ if $\sigma>n$ and $\varepsilon_2>0$ such that  $n(\frac{1}{p}-\frac{1}{q})\le{\sigma-\alpha-\beta-\varepsilon_{2}}$ if $\sigma=n$. Therefore, for $\varepsilon=\varepsilon_{1}$ or $\varepsilon_2$, we have
    \begin{align*}
        |\mathcal{I}_{2}^{0}|\,
        \lesssim\,
        \int_{|y|\le1}\diff{y}\,f(y)\,
        \underbrace{\vphantom{\int_{|x|\le\frac{|y|}{2}}}\,
        |y|^{\sigma-\alpha-\beta-\varepsilon}
        }_{\le\,|y|^{n(\frac{1}{p}-\frac{1}{q})}}\,
        \underbrace{
        |y|^{-\sigma+\beta+\varepsilon}\,\int_{|x|\le\frac{|y|}{2}}\diff{x}\,
        (\psi_{0}\,k_{\zeta,\sigma})(y^{-1}x)\,
        |x|^{-\beta}\,g(x)
        }_{:=\,(\tilde{V}_{\beta}g)(y)}
    \end{align*}
in both cases. Moreover, the term $\chi_{0}\tilde{V}_{\beta}g$ has same properties as \eqref{V1} and \eqref{V2} for all $\sigma\ge{n}$. Hence, the $L^p$-$L^q$-boundedness of $T_{2}^{0}(\zeta,\sigma,\alpha,\beta)$ follows by the similar arguments as in the case where $0<\sigma<n$. Then the proof is complete.
\end{proof}

Combining Proposition \ref{PropT01} and Proposition \ref{PropT023}, we obtain the following corollary for the $L^{p}$-$L^{q}$-boundedness of the operator $T^{0}(\zeta,\sigma,\alpha,\beta)$, which is the sufficiency part of Theorem \ref{mainthm1}.\textit{(1)}. 
\begin{corollary}
    Suppose that $\zeta\ge0$, $\sigma>0$ ($0<\sigma<\nu$ if $\zeta=0$), $\alpha\in\mathbb{R}$, $\beta\in\mathbb{R}$, and $0\le\alpha+\beta\le\sigma$. Then, the operator $T^{0}(\zeta,\sigma,\alpha,\beta)$ is bounded from $L^{p}(\mathbb{X})$ to $L^{q}(\mathbb{X})$ for all $1\le{p}\le{q}\le{\infty}$ satisfying $\frac{1}{p'}>\frac{\alpha}{n}$ when $\alpha>0$, $\frac{1}{q}>\frac{\beta}{n}$ when $\beta>0$, $\frac{1}{p}-\frac{1}{q}\le\frac{\sigma-\alpha-\beta}{n}$ when $\sigma-\alpha-\beta<n$, and $(\frac{1}{p},\frac{1}{q})\neq(\frac{\sigma-\alpha-\beta}{n},0)$ or $(1,1-\frac{\sigma-\alpha-\beta}{n})$.    
\end{corollary}

\subsection{$L^{p}$-$L^q$-boundedness of $T^{\infty}$}
The analysis changes dramatically in this case because the geometry at infinity now plays its part. We will establish the $L^p$-$L^q$-boundedness of the operators $T_{j}^{\infty}$, which are defined by
\begin{align*}
    T_{j}^{\infty}(\zeta,\sigma,\alpha,\beta)f(x)\,=\,
    \int_{\mathbb{X}}\diff{y}\,
    \mathds{1}_{Z_{j}}(x,y)\,\psi_{\infty}(y^{-1}x)\,
    |x|^{-\beta}\,k_{\zeta,\sigma}(y^{-1}x)\,
    |y|^{-\alpha}\,f(y),
\end{align*}
for $j=1,4,5,6,7$. Recall that $T_{j}^{\infty}(\zeta,\sigma,\alpha,\beta)f$ vanishes for any $f$ if $j=2$ or $3$. Since the results depend on the values of the parameters $\zeta$, $\sigma$, $p$ and $q$, we introduce the following notation to unify different subcases and simplify the statements. For all $\zeta\ge0$, $1\le{s}\le+\infty$, $\omega\in\mathbb{R}$, and $1\le{a}<{b}\le+\infty$, we denote by
\begin{align}
    \Phi(\zeta,s,\omega,a,b)\,
    =\,
    \begin{cases}
        b^{-M},\,\forall\,M>0
        &\quad\textnormal{if}\,\,\,\zeta>|\rho|\,\,\,
        \textnormal{or}\,\,\,\frac{1}{s}<\frac{1}{2}+\frac{\zeta}{2|\rho|},\\[10pt]
        \displaystyle \Big\lbrace{\int_{a}^{b}\diff{r}\,r^{s(\omega+\frac{\sigma}{2}-\frac{\ell+1}{2}\frac{1}{s'})-1}}\Big\rbrace^{1/s}
        &\quad\textnormal{if}\,\,\,\frac{1}{s}=\frac{1}{2}+\frac{\zeta}{2|\rho|}\,\,\,
        \textnormal{with}\,\,\,0<\zeta\le|\rho|,\\[10pt]
        \displaystyle \Big\lbrace{\int_{a}^{b}\diff{r}\,r^{2(\omega+\sigma-\frac{\nu}{2})-1}}\Big\rbrace^{1/2}
        &\quad\textnormal{if}\,\,\,\zeta=0\,\,\,
        \textnormal{and}\,\,\,s=2.
    \end{cases}
    \label{DefPhi}
\end{align}
For any $1\le{a}<{b}\le+\infty$, let $\chi_{a,b}\in\mathcal{C}_{c}^{\infty}(\mathbb{R}_{+}^{*})$ be a cut-off function compactly supported in the interval $[a,b]$. We denote by $\psi_{a,b}(x)=\chi_{a,b}(|x|)$ the corresponding bi-$K$-invariant cut-off function on $\mathbb{X}$. Then its support is contained in the annulus $\lbrace{x\in\mathbb{X}\,|\,a\le|x|\le{b}}\rbrace$. The crucial components of our analysis in this part are the following two lemmas.
\begin{lemma}\label{lemma omega1}
Suppose that $\zeta\ge0$, $\sigma>0$ ($0<\sigma<\nu$ if $\zeta=0$), $\omega\in\mathbb{R}$, and $1\le{a}<{b}\le+\infty$. Then, for all $1\le{s}\le+\infty$ such that $\frac{1}{s}\le\frac{1}{2}+\frac{\zeta}{2|\rho|}$, we have
    \begin{align*}
        \||\cdot|^{\omega}\,\psi_{a,b}\,k_{\zeta,\sigma}
        \|_{L^{s}(\mathbb{X})}\,
        \lesssim\,
        \Phi(\zeta,s,\omega,a,b).
    \end{align*}
\end{lemma}

\begin{proof}
For the sake of simplicity, let us denote by
\begin{align}
    m_{\zeta}=\,
    \begin{cases}
        \frac{\sigma-\nu-1}{2}
        &\qquad\textnormal{if}\,\,\,\zeta>0,\\[5pt]
        \sigma-\nu
        &\qquad\textnormal{if}\,\,\,\zeta=0
        \,\,\,\textnormal{and}\,\,\,0<\sigma<\nu.
    \end{cases}
    \label{mzeta}
\end{align}
Then, the kernel estimate \eqref{estim kinf} becomes $(\psi_{\infty}k_{\zeta,\sigma})(x)\asymp|x|^{m_{\zeta}}\varphi_{0}(x)\,e^{-\zeta|x|}$ for all $x\in\mathbb{X}$. By using successively this kernel estimate, the Cartan decomposition, the estimates \eqref{density} of the density and \eqref{phi0} of the ground spherical function, we obtain
    \begin{align}
        \||\,.\,|^{\omega}\,\psi_{a,b}k_{\zeta,\sigma}
        \|_{L^{s}(\mathbb{X})}^{s}\,
        &\asymp\,
        \int_{\mathbb{X}}\diff{x}\,|x|^{s\omega}\,\psi_{a,b}^{s}(x)\,
        |x|^{sm_{\zeta}}\,\varphi_{0}^{s}(x)\,e^{-s\zeta|x|}\notag\\[5pt]
        &=\,
        \int_{\mathfrak{a}^{+}}\diff{x^{+}}\,
        \delta(x^{+})\,\psi_{a,b}^{s}(x^{+})\,|x^{+}|^{s(\omega+m_{\zeta})}\,
        \varphi_{0}^{s}(\exp{x^{+}})\,e^{-s\zeta|x^{+}|}\notag\\[5pt]
        &\lesssim\,
         \int_{\mathfrak{a}^{+}}\diff{x^{+}}\,\psi_{a,b}(x^{+})\,
         |x^{+}|^{s(\omega+m_{\zeta}+\frac{\nu-\ell}{2})}\,
         e^{(2-s)\langle{\rho,x^{+}}\rangle}\,e^{-s\zeta|x^{+}|},
         \label{lemma1Lr}
    \end{align}
for all $1\le{s}<+\infty$. By noting that $s\zeta|x^{+}|\ge\frac{s\zeta}{|\rho|}\langle{\rho,x^{+}}\rangle$, the estimate \eqref{lemma1Lr} reads
\begin{align*}
    \||\,.\,|^{\omega}\,\psi_{a,b}k_{\zeta,\sigma}
    \|_{L^{s}(\mathbb{X})}^{s}\,
    &\lesssim\,
    \int_{\mathfrak{a}^{+}}\diff{x^{+}}\,\psi_{a,b}(x^{+})\,
    |x^{+}|^{s(\omega+m_{\zeta}+\frac{\nu-\ell}{2})}\,
    e^{-2(s(\frac{1}{2}+\frac{\zeta}{2|\rho|})-1)\langle{\rho,x^{+}}\rangle},
\end{align*}
which is $\mathrm{O}(b^{-M})$ for all $M>0$, if $s(\frac{1}{2}+\frac{\zeta}{2|\rho|})-1$ is strictly positive, in other words, when $\zeta>|\rho|$ or $\frac{1}{s}<\frac{1}{2}+\frac{\zeta}{2|\rho|}$. Similar result holds for $s=\infty$ since the kernel $\psi_{\infty}k_{\zeta,\sigma}$ includes an exponential decay. 

On the other hand, when $\zeta=0$ and $s=2$, we deduce straightforwardly from \eqref{lemma1Lr} that 
    \begin{align*}
        \||\,.\,|^{\omega}\,\psi_{a,b}k_{\zeta,\sigma}
        \|_{L^{2}(\mathbb{X})}^{2}\,
        &\lesssim\,
        \int_{\mathfrak{a}^{+}}\diff{x^{+}}\,\psi_{a,b}(x^{+})\,
        |x^{+}|^{2(\omega+\sigma-\nu+\frac{\nu-\ell}{2})}\\[5pt]
        &\le\,
        \int_{a}^{b}\diff{r}\,
        r^{2(\omega+\sigma-\frac{\nu}{2})-1}\,
        =\,
        \Phi^{2}(0,2,\omega,a,b).
    \end{align*}

We are left with the critical case where $\frac{1}{s}=\frac{1}{2}+\frac{\zeta}{2|\rho|}$ with $0<\zeta\le|\rho|$. By introducing polar coordinates in \eqref{lemma1Lr}, we write, for rank $\ell\ge2$,
\begin{align}
    \||\,.\,|^{\omega}\,\psi_{a,b}k_{\zeta,\sigma}
    \|_{L^{s}(\mathbb{X})}^{s}\,
    &\lesssim\,
    \int_{a}^{b}\diff{r}\,r^{\ell-1}
    r^{s(\omega+\frac{\sigma-\ell-1}{2})}\,e^{-s\zeta{r}}\,
    \int_{0}^{\theta_{\rho}}\diff{\theta}\,(\sin\theta)^{\ell-2}\,
    e^{(2-s)(\cos\theta)|\rho|r}\,
    \label{Lemma1theta1}
\end{align}
where $0\le\theta_{\rho}<\frac{\pi}{2}$ denotes the angle between vectors $\rho$ and $x^{+}$. Using the elementary formula $\cos\theta=1-2\sin^{2}\frac{\theta}{2}$ and the fact that $\sin\theta\asymp\theta$ around the origin, we have
\begin{align}
    \int_{0}^{\theta_{\rho}}\diff{\theta}\,(\sin\theta)^{\ell-2}\,
    e^{(2-s)(\cos\theta)|\rho|r}\,
    &\asymp\,
    e^{(2-s)|\rho|r}\,
    \int_{0}^{\theta_{\rho}}\diff{\theta}\,\theta^{\ell-2}\,
    e^{-\frac{(2-s)|\rho|r}{2}\theta^{2}}\notag\\[5pt]
    &\asymp\,
    r^{-\frac{\ell-1}{2}}\,e^{(2-s)|\rho|r}\,
    \int_{0}^{C\theta_{\rho}}\diff{t}\,t^{\ell-2}\,e^{-t^{2}},
    \label{Lemma1theta2}
\end{align}
where $C=\sqrt{\frac{(2-s)|\rho|r}{2}}$. Here, we have substituted the variable $\theta=C^{-1}t$. Note that $s<2$ in the present case. We deduce from \eqref{Lemma1theta1} and \eqref{Lemma1theta2} that
    \begin{align*}
        \||\,.\,|^{\omega}\,\psi_{a,b}k_{\zeta,\sigma}
        \|_{L^{s}(\mathbb{X})}^{s}\,
        &\lesssim\,
        \int_{a}^{b}\diff{r}\,
        r^{s(\omega+\frac{\sigma}{2}-\frac{\ell+1}{2}\frac{1}{s'})-1}\,
        e^{-2|\rho|(s(\frac{1}{2}+\frac{\zeta}{2|\rho|})-1)r},
    \end{align*}    
which also holds in rank one since $\rho$ and $x^{+}$ are two positive numbers in \eqref{lemma1Lr}. Note that the exponential factor in the last integral has no contributions in the critical case where $\frac{1}{s}=\frac{1}{2}+\frac{\zeta}{2|\rho|}$. So the proof is complete.
\end{proof}

\begin{remark}
Due to Lemma \ref{lemma omega1}, we know that the bi-$K$-invariant function $(|\cdot|^{\omega}\psi_{1,\infty}k_{\zeta,\sigma})\in{L^{s}(\mathbb{X})}$ for all $\omega\in\mathbb{R}$ if $\zeta>|\rho|$ or $\frac{1}{s}<\frac{1}{2}+\frac{\zeta}{2|\rho|}$. On the critical line $\frac{1}{s}=\frac{1}{2}+\frac{\zeta}{2|\rho|}$ with $0\le\zeta\le|\rho|$, the function $(|\cdot|^{\omega}\psi_{1,\infty}k_{\zeta,\sigma})$ is still in $L^s(\mathbb{X})$ provided that $\omega<-\frac{\sigma}{2}+\frac{\ell+1}{2}\frac{1}{s'}$ when $0<\zeta\le|\rho|$ or $\omega<-\sigma-\frac{\nu}{2}$ when $\zeta=0$. 
\label{Remark omega}
\end{remark}

\begin{lemma}\label{lemma omega2}
Suppose that $\zeta\ge0$, $\sigma>0$ ($0<\sigma<\nu$ if $\zeta=0$), $\omega\in\mathbb{R}$, and $1\le{a}<{b}\le+\infty$. Then we have the following.
    \begin{enumerate}[leftmargin=*,parsep=5pt]
        \item (Off-diagonal case) For all $1\le{p}<{q}\le+\infty$ satisfying $\frac{1}{p}\ge\frac{1}{2}-\frac{\zeta}{2|\rho|}$ and $\frac{1}{q}\le\frac{1}{2}+\frac{\zeta}{2|\rho|}$, we have
        \begin{align*}
            \|f*(|\cdot|^{\omega}\,\psi_{a,b}\,
            k_{\zeta,\sigma})\|_{L^{q}(\mathbb{X})}\,
            \lesssim\,
            \max\lbrace{
                \Phi(\zeta,p',\omega,a,b),\,\Phi(\zeta,q,\omega,a,b)\,
            }\rbrace\,
            \|f\|_{L^{p}(\mathbb{X})}.
        \end{align*}

        \item (Diagonal case) For all $1\le{p}\le+\infty$, we have
            \begin{align*}
                \|f*(|\cdot|^{\omega}\,\psi_{a,b}\,
                k_{\zeta,\sigma})\|_{L^{p}(\mathbb{X})}\,
                \lesssim\,
                \Phi^{\#}(\zeta,p,\omega,a,b)\,
                \|f\|_{L^{p}(\mathbb{X})},
            \end{align*}
        where
        \begin{align}\label{DefPhiD}
            \Phi^{\#}(\zeta,p,\omega,a,b)\,=\,
            \begin{cases}
                b^{-M},\,\forall\,M>0
                &\quad\textnormal{if}\,\,\,
                \zeta>|\rho|\,\,\,\textnormal{or}\,\,\,
                \frac{1}{2}-\frac{\zeta}{2|\rho|}<\frac{1}{p}<\frac{1}{2}+\frac{\zeta}{2|\rho|},\\[10pt]
                \displaystyle \int_{a}^{b}\diff{r}\,r^{\omega+\frac{\sigma}{2} -1}
                &\quad\textnormal{if}\,\,\,\frac{1}{p}=\frac{1}{2}\pm\frac{\zeta}{2|\rho|}\,\,\,
                \textnormal{with}\,\,\,0<\zeta\le|\rho|,\\[10pt]
                \displaystyle \int_{a}^{b}\diff{r}\,r^{\omega+\sigma-1}
                &\quad\textnormal{if}\,\,\,\zeta=0\,\,\,
                \textnormal{and}\,\,\,p=2.
            \end{cases}
        \end{align}
    \end{enumerate}
\end{lemma}

\begin{remark}
    Observe that when there is no weight, i.e., $\omega=0$, the term
    \begin{align*}
        \Phi^{\#}\Big(\zeta,\frac{2|\rho|}{|\rho|\pm\zeta},0,a,\infty\Big)\,
        =\,
        \int_{a}^{+\infty}\diff{r}\,r^{\frac{\sigma}{2}+(\frac{1}{2}-\frac{\zeta}{2|\rho|})(\nu-\ell)-1}
    \end{align*}
    is not finite for any $0<\zeta\le|\rho|$, since $\sigma>0$ and $\nu-\ell=2|\Sigma_{r}^{+}|>0$. This highlights the breakdown of the diagonal Hardy-Littlewood-Sobolev inequality at critical points when $\zeta$ is small. Similarly, the term $\Phi^{\#}(\zeta,2,0,a,\infty)$ diverges in the limiting case where $\zeta=0$, and we know that the operator $(-\Delta-|\rho|^{2})^{-\frac{\sigma}{2}}$ is not $L^2$-bounded.
\end{remark}

\begin{figure}
     \centering
        \input{Lemma1.tex}
        \caption{This figure corresponds to the case where $0<\zeta<|\rho|$. The critical lines, given by $\frac{1}{p} = \frac{1}{2} - \frac{\zeta}{2|\rho|}$ and $\frac{1}{q} = \frac{1}{2} + \frac{\zeta}{2|\rho|}$, contract to the endpoints $(0,0)$ and $(1,1)$ as $\zeta$ approaches $|\rho|$. They disappear for $\zeta>|\rho|$. For the diagonal case where $\frac{1}{p} = \frac{1}{q}$, the valid range is represented by the red segment. This segment contracts to the point $(\frac{1}{2}, \frac{1}{2})$ as $\zeta$ reaches $0$.}
        \label{lemma2figure}
\end{figure}
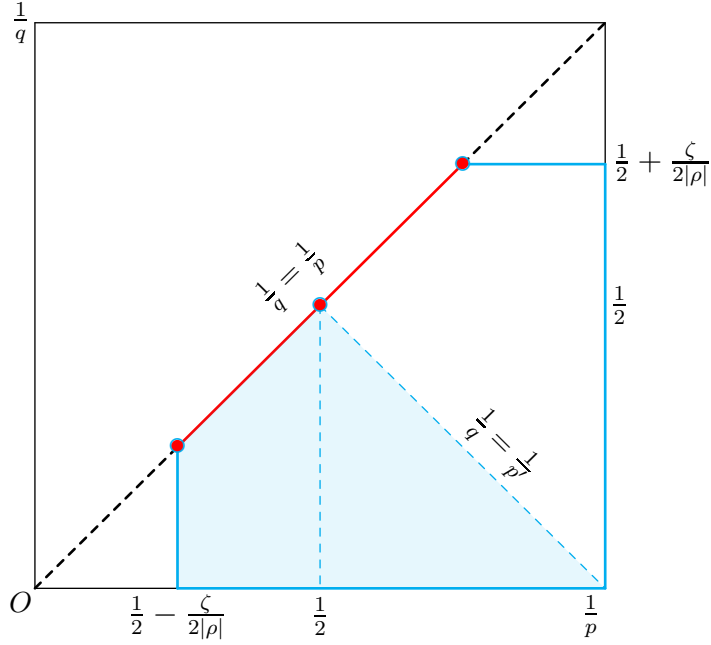

\begin{proof}
This lemma follows from Lemma \ref{lemma omega1} and the Kunze-Stein phenomenon. For the sake of simplicity, let us denote by $B(x)=|x|^{\omega}\psi_{a,b}(x)k_{\zeta,\sigma}(x)$ which is a bi-$K$-invariant function on $\mathbb{X}$. 

\textit{(1)} Note that $p\neq\infty$ and $q\neq1$ in the off-diagonal case. When $q=\infty$, we derive directly from Young's inequality that
\begin{align}
    \|f*B\|_{L^{\infty}(\mathbb{X})}\,
    \lesssim\,
    \|f\|_{L^{p}(\mathbb{X})}\,
    \|B\|_{L^{p'}(\mathbb{X})}
    \label{conv1}
\end{align}
for all $1\le{p}<+\infty$. When $1<q<+\infty$, we first consider the case where $\frac{1}{q}<\frac{1}{p'}$, see the region shaded in blue in Figure \ref{lemma2figure}. On the one hand, when $0<\frac{1}{q}<\frac{1}{p'}\le\frac{1}{2}$, that is $1<q'<p\le2$, we deduce from \eqref{KS1} with $t=p$ and $s=q'$ that
   \begin{align}
        \|f*B
        \|_{L^{q}(\mathbb{X})}\,
        \le\,
        \|f\|_{L^{p}(\mathbb{X})}\,
        \|B\|_{L^{p'}(\mathbb{X})}.
        \label{conv2}
    \end{align}
On the other hand, by taking $t=p'$ and $s=q'$ after permuting $L^{t}(\mathbb{X})$ and $L^{t'}(\mathbb{X})$ in \eqref{KS1}, one can show that the estimate \eqref{conv2} still holds when $\frac{1}{2}\le\frac{1}{p'}\le\frac{1}{2}+\frac{\zeta}{2|\rho|}$, that is $\frac{1}{2}-\frac{\zeta}{2|\rho|}\le\frac{1}{p}\le\frac{1}{2}$.

In the case where $\frac{1}{q}>\frac{1}{p'}$, we deduce from \eqref{conv2} and the duality that 
    \begin{align}
         \|f*B\|_{L^{q}(\mathbb{X})}\,
         &\le\,
         \sup_{\|h\|_{L^{q'}(\mathbb{X})}\le1}\,
         \|f\|_{L^{p}(\mathbb{X})}\,
         \|h*\check{B}\|_{L^{p'}(\mathbb{X})}\notag\\[5pt]
         &\le\,
         \sup_{\|h\|_{L^{q'}(\mathbb{X})}\le1}\,
          \|f\|_{L^{p}(\mathbb{X})}\,
          \|h\|_{L^{q'}(\mathbb{X})}\,
          \|\check{B}\|_{L^{q}(\mathbb{X})}\,
         \le\,
         \|f\|_{L^{p}(\mathbb{X})}\,
         \|B\|_{L^{q}(\mathbb{X})},\,
          \label{conv3}
    \end{align}
where $\check{B}(x)=B(x^{-1})$ satisfies $\|\check{B}\|_{L^{p}(\mathbb{X})}=\|B\|_{L^{p}(\mathbb{X})}$ for all $p\ge1$.
    
In the remaining case where $\frac{1}{q}=\frac{1}{p'}$, we note from Figure \ref{lemma2figure} that $2<q\le+\infty$. The endpoint case $q=+\infty$ is already handled as in \eqref{conv1}. For $2<q<+\infty$, we obtain, by using \eqref{KSbK2} with $s=q=p'$, that
    \begin{align*}
        \|f*B\|_{L^{q}(\mathbb{X})}\,
        \le\,
        \|f\|_{L^{p}(\mathbb{X})}\,
        \Big\lbrace{
        \int_{G}\diff{x}\,\varphi_{0}(x)\,
        |B(x)|^{\frac{q}{2}}
        }\Big\rbrace^{\frac{2}{q}},
    \end{align*}
where
    \begin{align*}
        \int_{G}\diff{x}\,\varphi_{0}(x)\,
        |B(x)|^{\frac{q}{2}}\,
        &\asymp\,
        \int_{\mathfrak{a}^{+}}\diff{x^{+}}\,\delta(x^{+})\,
        \psi_{a,b}^{\frac{q}{2}}(x^{+})\,
        |x^{+}|^{\frac{q(\omega+m_{\zeta})}{2}}\,
        \varphi_{0}^{\frac{q}{2}+1}(\exp{x^{+}})\,
        e^{-\frac{q\zeta}{2}|x^{+}|}\\[5pt]
        &\lesssim\,
        \int_{\mathfrak{a}^{+}}\diff{x^{+}}\,\psi_{a,b}(x^{+})\,
        |x^{+}|^{\frac{q(\omega+m_{\zeta})}{2}
        +(\frac{q}{2}+1)\frac{\nu-\ell}{2}}\,
        e^{(1-\frac{q}{2})\langle{\rho,x^{+}}\rangle}\,
        e^{-\frac{q\zeta}{2}|x^{+}|},
    \end{align*}
according to estimates \eqref{estim kinf}, \eqref{phi0}, and \eqref{density}. Since $q>2$, we obtain, when $\frac{1}{q}=\frac{1}{p'}$,
\begin{align}
    \|f*B\|_{L^{q}(\mathbb{X})}\,
   \lesssim\,
   b^{-M}\,\|f\|_{L^{p}(\mathbb{X})}.
   \label{conv4}
\end{align}
for any $M>0$.

By combining \eqref{conv1}, \eqref{conv2}, \eqref{conv3}, and \eqref{conv4} together, the statement \textit{(1)} follows from Lemma \ref{lemma omega1}.

\textit{(2)} Now, let us turn to the diagonal case where $1\le{p}=q\le+\infty$. By using the bi-$K$-invariant Kunze-Stein phenomenon \eqref{KSbK1}, we obtain
\begin{align}
    \|f*B\|_{L^{p}(\mathbb{X})}\,
    =\,
    \|f\|_{L^{p}(\mathbb{X})}\,
    \underbrace{\vphantom{\int}
    \int_{G}\diff{x}\,\psi_{a,b}(x)\,
    |x|^{\omega}\,k_{\zeta,\sigma}(x)\,
    \varphi_{i(\frac{2}{p}-1)\rho}(x),
    }_{:=\,I}
    \label{integralI}
\end{align}
where $\varphi_{i(\frac{2}{p}-1)\rho}(\exp{x^{+}})\asymp{e^{-\frac{2}{p}\langle{\rho,x^{+}}\rangle}}$ for $p>2$, according to \eqref{phirho}. On the one hand, for $\frac{1}{2}-\frac{\zeta}{2|\rho|}\le\frac{1}{p}<\frac{1}{2}$ with $0<\zeta\le|\rho|$, we have
\begin{align*}
    I\,
    &\asymp\,
    \int_{\mathfrak{a^{+}}}\diff{x^{+}}\,
    \delta(x^{+})\,\psi_{a,b}(x^{+})\,
    |x^{+}|^{\omega+\frac{\sigma-\nu-1}{2}}\,
    \varphi_{0}(\exp{x^{+}})\,
    e^{-\zeta|x^{+}|}\,e^{-\frac{2}{p}\langle\rho,x^{+}\rangle}\\[5pt]
    &\lesssim\,
    \int_{\mathfrak{a^{+}}}\diff{x^{+}}\,\psi_{a,b}(x^{+})\,
    |x^{+}|^{\omega+\frac{\sigma-\ell-1}{2}}\,
    e^{(1-\frac{2}{p})\langle\rho,x^{+}\rangle}\,
    e^{-\zeta|x^{+}|},
\end{align*}
which can be handled by using the same argument as in Lemma \ref{lemma omega1}. In more detail, when $\zeta>|\rho|$ or $\frac{1}{2}-\frac{\zeta}{2|\rho|}<\frac{1}{p}<\frac{1}{2}$, we know that $(1-\frac{2}{p})\langle\rho,x^{+}\rangle-\zeta|x^{+}|\le(1-\frac{2}{p}-\frac{\zeta}{|\rho|})\langle\rho,x^{+}\rangle$ is strictly negative, then $I$ is $\mathrm{O}(b^{-M})$ for all $M>0$ and $\omega\in\mathbb{R}$. In the critical case where $\frac{1}{p}=\frac{1}{2}-\frac{\zeta}{2|\rho|}$, there is no more exponential decay in the integral $I$, we conclude by using \eqref{Lemma1theta2} again. Similar result remains valid for $\frac{1}{2}<\frac{1}{p}\le\frac{1}{2}+\frac{\zeta}{2|\rho|}$ by duality.

On the other hand, in the limiting case where $\zeta=0$, by taking $p=2$ in \eqref{integralI}, we obtain similarly 
\begin{align*}
    \|f*B\|_{L^{2}(\mathbb{X})}\,
    \lesssim\,
    \|f\|_{L^{2}(\mathbb{X})}\,
    \int_{\mathfrak{a^{+}}}\diff{x^{+}}\,\psi_{a,b}(x^{+})\,
    |x^{+}|^{\omega+\sigma-\ell},
\end{align*}
according to \eqref{estim kinf} and \eqref{phi0}. Then the proof is complete.
\end{proof}

Now, let us establish the $L^{p}$-$L^{q}$-boundedness for the operators $T_{j}^{\infty}$ with $j=1,\,4,\,5,\,6,\,7$. For $j=1$, the following statement is proved by combining Lemma \ref{lemma omega2} with the dyadic decomposition techniques.
\begin{proposition}\label{T1inf prop}
Suppose that $\zeta\ge0$, $\sigma>0$ ($0<\sigma<\nu$ if $\zeta=0$), $\alpha\in\mathbb{R}$, and $\beta\in\mathbb{R}$. Then the operator 
$T_{1}^{\infty}(\zeta,\sigma,\alpha,\beta)$ is bounded from $L^p(\mathbb{X})$ to $L^q(\mathbb{X})$ for all $1\le{p}\le{q}\le+\infty$ such that $\frac{1}{p}\ge\frac{1}{2}-\frac{\zeta}{2|\rho|}$ and $\frac{1}{q}\le\frac{1}{2}+\frac{\zeta}{2|\rho|}$, with the constraints
\begin{align}\label{T1ResOff}
    \alpha+\beta\,>\,
    \begin{cases}
        \frac{\sigma}{2}-(\frac{1}{2}-\frac{\zeta}{2|\rho|})\frac{\ell+1}{2}
        &\qquad\textnormal{if}\,\,\,0<\zeta\le|\rho|,\\[5pt]
        \sigma-\frac{\nu}{2}
        &\qquad\textnormal{if}\,\,\,\zeta=0.
    \end{cases}
\end{align}
either in the critical and off-diagonal cases where $\frac{1}{p}=\frac{1}{2}-\frac{\zeta}{2|\rho|}>\frac{1}{q}$ or $\frac{1}{q}=\frac{1}{2}+\frac{\zeta}{2|\rho|}<\frac{1}{p}$, or
\begin{align}\label{T1ResOn}
    \alpha+\beta\,\ge\,
    \begin{cases}
        \frac{\sigma}{2}
        &\qquad\textnormal{if}\,\,\,0<\zeta\le|\rho|,\\[5pt]
        \sigma
        &\qquad\textnormal{if}\,\,\,\zeta=0.
    \end{cases}
\end{align}
in the diagonal case where $\frac{1}{p}=\frac{1}{q}=\frac{1}{2}\pm\frac{\zeta}{2|\rho|}$.
\end{proposition}

\begin{remark}
The condition \eqref{T1ResOn} induces \eqref{T1ResOff}, which means that the preservation of the $L^p$-boundedness comes with stricter assumptions. In particular, the condition \eqref{T1ResOn} is invalid when $\alpha=\beta=0$, which corresponds to the fact that the diagonal Hardy-Littlewood-Sobolev inequality fails in the critical cases.  
\end{remark}

\begin{proof}
As we turn to the cases that are distant from the origin, we make slight adjustments to the cut-off functions defined in \eqref{dyadicZ}. From now on, we denote by $\Psi_{0}(y)=\psi_{0}(y)$ and 
    \begin{align*}
        \Psi_{k}(y)\,
        =\,
        \chi_{0}(2^{-k}|y|)-\chi_{0}(2^{-k+1}|y|)
        \qquad\forall\,k\ge1,\,\,\,
        \forall\,y\in\mathbb{X},
    \end{align*}
which are all bi-$K$-invariant cut-off functions on $\mathbb{X}$. Note that $\Psi_{k}$ is compactly supported in $\lbrace{y\in\mathbb{X}\,|\,2^{k-2}\le|y|\le2^{k}}\rbrace$ for every $k\ge1$, and we have the partition of unity
    \begin{align}
        \sum_{k\in\mathbb{N}}\Psi_{k}(y)\,=\,1
        \qquad\forall\,y\in\mathbb{X}.
        \label{dyadicN}
    \end{align}

Recall that, for all $(x,y)\in{Z_{1}}$ and $y\in\supp{\Psi_{k}}$, we have $|x|^{-\beta}|y|^{-\alpha}\asymp|y|^{-(\alpha+\beta)}\asymp2^{-(\alpha+\beta)k}$ and $|y^{-1}x|\lesssim2^{k}$. Hence, for all $1\le{q}<+\infty$, 
\begin{align}
    \|T_{1}^{\infty}(\zeta,\sigma,\alpha,\beta)f\|_{L^{q}(\mathbb{X})}^{q}\,
    &=\,
    \int_{\mathbb{X}}\diff{x}\,
    \Big|\int_{\mathbb{X}}\diff{y}\,
    \mathds{1}_{Z_{1}}(x,y)\,|x|^{-\beta}\,
    (\psi_{\infty}k_{\zeta,\sigma})(y^{-1}x)\,
    |y|^{-\alpha}\,f(y)
    \Big|^{q}\notag\\[5pt]
    &\lesssim\,
    \sum_{k\in\mathbb{N}}\,2^{-(\alpha+\beta)qk}\,
    \int_{\mathbb{X}}\diff{x}\,
    \Big|\int_{\mathbb{X}}\diff{y}\,\mathds{1}_{Z_{1}}(x,y)\,
    (\psi_{\infty}k_{\zeta,\sigma})(y^{-1}x)\,(\Psi_{k}f)(y)
    \Big|^{q}\notag\\[5pt]
    &\lesssim\,
    \sum_{k\in\mathbb{N}}\,2^{-(\alpha+\beta)qk}\,
    \|\Psi_{k}f*\psi_{1,2^{k}}k_{\zeta,\sigma}\|_{L^{q}(\mathbb{X})}^{q},
    \label{T1 estim0}
\end{align}
and similarly,
\begin{align}
    \|T_{1}^{\infty}(\zeta,\sigma,\alpha,\beta)f\|_{L^{\infty}(\mathbb{X})}\,
    \lesssim\,
    \sum_{k\in\mathbb{N}}\,2^{-(\alpha+\beta)qk}\,
    \|\Psi_{k}f*\psi_{1,2^{k}}k_{\zeta,\sigma}\|_{L^{\infty}(\mathbb{X})}.
    \label{T1 estim0sup}
\end{align}

Now, we study $\|\Psi_{k}f*\psi_{1,2^{k}}k_{\zeta,\sigma}\|_{L^{q}(\mathbb{X})}$ for all $1\le{q}\le+\infty$. In the subcase where $\frac{1}{p}$ and $\frac{1}{q}$ are way from the critical lines, in other words, for either $\zeta>|\rho|$ or $\frac{1}{p}>\frac{1}{2}-\frac{\zeta}{2|\rho|}$ and $\frac{1}{q}<\frac{1}{2}+\frac{\zeta}{2|\rho|}$, we have
\begin{align}
    \|\Psi_{k}f*\psi_{1,2^{k}}k_{\zeta,\sigma}\|_{L^{q}(\mathbb{X})}\,
    \lesssim\,
    2^{-Mk}\,\|\Psi_{k}f\|_{L^{p}(\mathbb{X})},
    \label{T1 estim1}
\end{align}
for any $k\in\mathbb{N}$ and $M>0$, according to Lemma \ref{lemma omega2} and the definitions \eqref{DefPhi} and \eqref{DefPhiD} of $\Phi$ and $\Phi^{\#}$. 

In the critical and off-diagonal case where $\frac{1}{p}=\frac{1}{2}-\frac{\zeta}{2|\rho|}>\frac{1}{q}$ or $\frac{1}{q}=\frac{1}{2}+\frac{\zeta}{2|\rho|}<\frac{1}{p}$ with $0\le\zeta\le|\rho|$, we deduce from Lemma \ref{lemma omega2}.\textit{(i)} that 
\begin{align*}
    \|\Psi_{k}f*\psi_{1,2^{k}}k_{\zeta,\sigma}\|_{L^{q}(\mathbb{X})}\,
    \lesssim\,
    \max\Big\lbrace{
        \Phi\Big(\zeta,\frac{2|\rho|}{|\rho|+\zeta},0,1,2^{k}\Big),\,\Phi(\zeta,q,0,1,2^{k})\,
    }\Big\rbrace\,
    \Big\|\Psi_{k}f\Big\|_{L^{\frac{2|\rho|}{|\rho|-\zeta}}(\mathbb{X})},
\end{align*}
provided that $\frac{1}{q}<\frac{1}{2}-\frac{\zeta}{2|\rho|}$, and 
\begin{align*}
    \Big\|\Psi_{k}f*\psi_{1,2^{k}}k_{\zeta,\sigma}\Big\|_{L^{\frac{2|\rho|}{|\rho|+\zeta}}(\mathbb{X})}\,
    \lesssim\,
    \max\Big\lbrace{
        \Phi(\zeta,p',0,1,2^{k}),\,\Phi\Big(\zeta,\frac{2|\rho|}{|\rho|+\zeta},0,1,2^{k}\Big)\,
    }\Big\rbrace\,
    \|\Psi_{k}f\|_{L^{p}(\mathbb{X})},
\end{align*}
provided that $\frac{1}{p}>\frac{1}{2}+\frac{\zeta}{2|\rho|}$, which is equivalent to $\frac{1}{p'}<\frac{1}{2}-\frac{\zeta}{2|\rho|}$. According to the definition \eqref{DefPhi} of $\Phi$, we know, in the present case, that
\begin{align*}
    \max\Big\lbrace{\,
        \Phi(\zeta,p',0,1,2^{k}),\,
        \Phi(\zeta,q,0,1,2^{k}),\,
        \Phi\Big(\zeta,\frac{2|\rho|}{|\rho|+\zeta},0,1,2^{k}\Big)
    }\Big\rbrace\,
    =\,
    \Phi\Big(\zeta,\frac{2|\rho|}{|\rho|+\zeta},0,1,2^{k}\Big),
\end{align*}
and
\begin{align*}
    2^{-(\alpha+\beta)}\,
    \Phi\Big(\zeta,\frac{2|\rho|}{|\rho|+\zeta},0,1,2^{k}\Big)\,
    =\,\mathrm{O}(1),
\end{align*}
provided that $\alpha+\beta$ satisfies \eqref{T1ResOff}. Together with estimates \eqref{T1 estim0}, \eqref{T1 estim0sup}, and \eqref{T1 estim1}, we obtain the $L^p$-$L^q$-boundedness of the operator $T_1^{\infty}(\zeta,\sigma,\alpha,\beta)$ in the off-diagonal case where $1\le{p}<{q}\le+\infty$, $\frac{1}{p}\ge\frac{1}{2}+\frac{\zeta}{2|\rho|}$, and $\frac{1}{q}\le\frac{1}{2}+\frac{\zeta}{2|\rho|}$, with the restriction \eqref{T1ResOff} on the critical lines.

We are left with the diagonal case where $\frac{1}{p}=\frac{1}{q}=\frac{1}{2}\pm\frac{\zeta}{2|\rho|}$ with $0\le\zeta\le|\rho|$. From Lemma \ref{lemma omega2}.\textit{(ii)}, we deduce similarly that
\begin{align*}
    \|\Psi_{k}f*\psi_{1,2^{k}}k_{\zeta,\sigma}\|_{L^{p}(\mathbb{X})}\,
    \lesssim\,
    \Phi^{\#}(\zeta,p,0,1,2^{k})\,
    \|\Psi_{k}f\|_{L^{p}(\mathbb{X})},
\end{align*}
with $\frac{1}{p}=\frac{1}{2}\pm\frac{\zeta}{2|\rho|}$. Due to the definition \eqref{DefPhiD} of $\Phi^{\#}$, we have
\begin{align*}
    2^{-(\alpha+\beta)}\,
    \Phi^{\#}\Big(\zeta,\frac{2|\rho|}{|\rho|\pm\zeta},0,1,2^{k}\Big)\,
    =\,\mathrm{O}(1),
\end{align*}
provided that $\alpha+\beta$ satisfies \eqref{T1ResOn}. We conclude that the operator $T_1^{\infty}(\zeta,\sigma,\alpha,\beta)$ is $L^p$-bounded for all $1\le{p}\le{\infty}$ such that $\frac{1}{2}-\frac{\zeta}{2|\rho|}\le\frac{1}{p}\le\frac{1}{2}+\frac{\zeta}{2|\rho|}$, under the restriction \eqref{T1ResOn} in the critical cases where $\frac{1}{p}=\frac{1}{2}\pm\frac{\zeta}{2|\rho|}$ with $0\le\zeta\le|\rho|$.
\end{proof}

Next, let us shift our focus to the operators $T_{j}^{\infty}$ for $4\le{j}\le7$. Building upon Lemma \ref{lemma omega1} and Lemma \ref{lemma omega2}, we will also use the local Harnack inequality for the cases $j=4$ and $j=5$ where one of the variables is closed to the origin, and  the weighted Kunze-Stein phenomenon for the cases $j=6$ and $j=7$ where both variables are away from the origin.
\begin{proposition}\label{T4T5 prop}
    Suppose that $\zeta\ge0$, $\sigma>0$ ($0<\sigma<\nu$ if $\zeta=0$), $\alpha\in\mathbb{R}$, and $\beta\in\mathbb{R}$. Then we have the following.
    \begin{enumerate}[leftmargin=*,parsep=5pt,label=(\roman*)]    
        \item For $1\le{p}\le{q}\le{\infty}$, $\frac{1}{p}\ge\frac{1}{2}-\frac{\zeta}{2|\rho|}$, and $\frac{1}{q}>\frac{\beta}{n}$ when $\beta>0$, the operator $T_{4}^{\infty}(\zeta,\sigma,\alpha,\beta)$ is bounded from $L^{p}(\mathbb{X})$ to $L^{q}(\mathbb{X})$, with the following constraints in the critical case where $\frac{1}{p}=\frac{1}{2}-\frac{\zeta}{2|\rho|}$:
        \begin{align}
            \alpha\,>\,
            \begin{cases}
                \frac{\sigma}{2}-(\frac{1}{2}-\frac{\zeta}{2|\rho|})\frac{\ell+1}{2}
                &\qquad\textnormal{if}\,\,\,0<\zeta\le|\rho|,\\[5pt]
                \sigma-\frac{\nu}{2}
                &\qquad\textnormal{if}\,\,\,\zeta=0.
            \end{cases}
            \label{T4Res}
        \end{align}

        \item For $1\le{p}\le{q}\le{\infty}$, $\frac{1}{q}\le\frac{1}{2}+\frac{\zeta}{2|\rho|}$, and $\frac{1}{p'}>\frac{\alpha}{n}$ when $\alpha>0$, the operator $T_{5}^{\infty}(\zeta,\sigma,\alpha,\beta)$ is bounded from $L^{p}(\mathbb{X})$ to $L^{q}(\mathbb{X})$, with the following constraints in the critical case where $\frac{1}{q}=\frac{1}{2}+\frac{\zeta}{2|\rho|}$:
        \begin{align}
            \beta\,>\,
            \begin{cases}
                \frac{\sigma}{2}-(\frac{1}{2}-\frac{\zeta}{2|\rho|})\frac{\ell+1}{2}
                &\qquad\textnormal{if}\,\,\,0<\zeta\le|\rho|,\\[5pt]
                \sigma-\frac{\nu}{2}
                &\qquad\textnormal{if}\,\,\,\zeta=0.
            \end{cases}
            \label{T5Res}
        \end{align}
    \end{enumerate}
\end{proposition}

\begin{proof}
This proof is based on the local Harnack inequality of the ground spherical function: if $x$ belongs to any compact subset of $G$, then we have $\varphi_{0}(y^{-1}x)\asymp\varphi_{0}(y)$ for all $y\in{G}$, see for instance \cite[Proposition 4.6.3]{GV88} or \cite[Remark 4.5]{APZ23}. 

Notice that for all $(x,y)\in{Z_{4}}$, we have $2|x|\le1\le|y|$, then $|y^{-1}x|\asymp|y|$. Therefore, for all $\zeta\ge0$ and $\sigma>0$, we have
    \begin{align*}
        k_{\zeta,\sigma}(y^{-1}x)\,
        \asymp\,
        |y^{-1}x|\,\varphi_{0}(y^{-1}x)\,e^{-\zeta|y^{-1}x|}\,
        \asymp\,
        k_{\zeta,\sigma}(y)\,
        \qquad\forall\,(x,y)\in{Z_{4}},
    \end{align*}
due to the local Harnack inequality. Hence, for $1\le{q}<+\infty$,
    \begin{align*}
         \|T_{4}^{\infty}(\zeta,\sigma,\alpha,\beta)f\|_{L^{q}(\mathbb{X})}^{q}\,
         &=\,
         \int_{\mathbb{X}}\diff{x}\,
        \Big|\int_{\mathbb{X}}\diff{y}\,
        \mathds{1}_{Z_{4}}(x,y)\,|x|^{-\beta}\,
        (\psi_{\infty}k_{\zeta,\sigma})(y^{-1}x)\,
        |y|^{-\alpha}\,f(y)
        \Big|^{q}\\[5pt]
        &\lesssim\,
        \underbrace{
        \int_{|x|\le\frac{1}{2}}\diff{x}\,|x|^{-q\beta}
        }_{=\,I_{1}}\,
        \Big|
        \underbrace{\vphantom{\int_{|x|\le\frac{1}{2}}}
        \int_{\mathbb{X}}\diff{y}\,
         (\psi_{\infty}k_{\zeta,\sigma})(y)\,|y|^{-\alpha}\,f(y)
        }_{=\,I_{2}}
        \Big|^{q},
    \end{align*}
where 
    \begin{align*}
        I_{1}\,
        \lesssim\,
        \int_{|x^{+}|\le\frac{1}{2}}\diff{x^{+}}\,
        |x^{+}|^{-q\beta+n-\ell}\,
        =\,
        \int_{0}^{\frac{1}{2}}\diff{r}\,r^{-q\beta+n-1}\,
    \end{align*}
which is finite for all $1\le{q}<\infty$ when $\beta\le0$ and for $\frac{1}{q}>\frac{\beta}{n}$ when $\beta>0$. Note that there is no positive $\beta$ that satisfies $\frac{1}{q}>\frac{\beta}{n}$ when $q=+\infty$. Then, by using the Holder inequality in $I_{2}$, we obtain, for $1\le{p}\le{q}\le+\infty$ satisfying $\frac{1}{q}>\frac{\beta}{n}$ when $\beta>0$,
\begin{align*}
    \|T_{4}^{\infty}(\zeta,\sigma,\alpha,\beta)f\|_{L^{q}(\mathbb{X})}\,
    \lesssim\,
    \|f\|_{L^{p}(\mathbb{X})}\,
    \||\cdot|^{-\alpha}\,\psi_{1,\infty}\,k_{\zeta,\sigma}
        \|_{L^{p'}(\mathbb{X})}.
\end{align*}
According to Lemma \ref{lemma omega1} and Remark \ref{Remark omega}, we know that 
\begin{align*}
    \||\cdot|^{-\alpha}\,\psi_{1,\infty}\,k_{\zeta,\sigma}
        \|_{L^{p'}(\mathbb{X})}\,
    \lesssim\,
    \Phi(\zeta,p',-\alpha,1,\infty),
\end{align*}
which is finite if one of the following conditions holds:
\begin{itemize}[parsep=5pt]
    \item $\zeta>|\rho|$;
    \item $0\le\zeta\le|\rho|$ and $\frac{1}{p'}<\frac{1}{2}+\frac{\zeta}{2|\rho|}$;
    \item $0<\zeta\le|\rho|$, $\frac{1}{p'}=\frac{1}{2}+\frac{\zeta}{2|\rho|}$, and $\alpha>\frac{\sigma}{2}-(\frac{1}{2}-\frac{\zeta}{2|\rho|})\frac{\ell+1}{2}$;
    \item $\zeta=0$, $\frac{1}{p'}=\frac{1}{2}$, and $\alpha>\sigma-\frac{\nu}{2}$.
\end{itemize}
By noting that $\frac{1}{p'}\le\frac{1}{2}+\frac{\zeta}{2|\rho|}$ is equivalent to $\frac{1}{p}\ge\frac{1}{2}-\frac{\zeta}{2|\rho|}$, we conclude that the operator $T_{4}^{\infty}(\zeta,\sigma,\alpha,\beta)$ is bounded from $L^{p}(\mathbb{X})$ to $L^{q}(\mathbb{X})$, for $1\le{p}\le{q}\le{\infty}$, $\frac{1}{p}\ge\frac{1}{2}-\frac{\zeta}{2|\rho|}$, and $\frac{1}{q}>\frac{\beta}{n}$ when $\beta>0$, with the restriction \eqref{T4Res} in the critical case where $\frac{1}{p}=\frac{1}{2}-\frac{\zeta}{2|\rho|}$. We omit the similar proof of \textit{(ii)} for the operator $T_{5}^{\infty}$.
\end{proof}

To prove the $L^{p}$-$L^q$-boundedness of the operator $T_{j}^{\infty}$ under the ideal assumptions for the remaining two cases where $j=6$ and $7$, we must optimize each parameter to its full potential. A critical step in the proof of the following proposition is to establish the sharp weighted Kunze-Stein phenomenon.
\begin{proposition}\label{T6T7 prop}
Suppose that $\zeta\ge0$, $\sigma>0$ ($0<\sigma<\nu$ if $\zeta=0$), $\alpha\in\mathbb{R}$, and $\beta\in\mathbb{R}$. Then the operator $T_{6}^{\infty}(\zeta,\sigma,\alpha,\beta)$ (resp. $T_{7}^{\infty}(\zeta,\sigma,\alpha,\beta)$) is bounded from $L^p(\mathbb{X})$ to $L^q(\mathbb{X})$ for all $1\le{p}\le{q}\le+\infty$ such that $\frac{1}{p}\ge\frac{1}{2}-\frac{\zeta}{2|\rho|}$ and $\frac{1}{q}\le\frac{1}{2}+\frac{\zeta}{2|\rho|}$, with the constraints \eqref{T4Res} (resp. \eqref{T5Res}), and either \eqref{T1ResOff} if $p\neq{q}$ or \eqref{T1ResOn} if $p=q$, in the critical cases where $\frac{1}{p}=\frac{1}{2}-\frac{\zeta}{2|\rho|}$ or $\frac{1}{q}=\frac{1}{2}+\frac{\zeta}{2|\rho|}$.
\end{proposition}

\begin{proof}
Recall that, for all $(x,y)$ in $Z_{6}$, we have $|y|\ge2|x|\ge1$ and $|y^{-1}x|\asymp|y|$, see Figure \ref{Regions}. Then
\begin{align*}
    |x|^{-\beta}|y|^{-\alpha}\,
    \lesssim\,
    |y^{-1}x|^{-\alpha-\min\lbrace{0,\beta}\rbrace},
\end{align*}
for all $\alpha,\,\beta\in\mathbb{R}$. Using the partition of unity \eqref{dyadicN} again, we obtain, for $1\le{q}<+\infty$,
\begin{align*}
    \|T_{6}^{\infty}(\zeta,\sigma,\alpha,\beta)f\|_{L^{q}(\mathbb{X})}^{q}\,
    &=\,
    \int_{\mathbb{X}}\diff{x}\,
    \Big|\int_{\mathbb{X}}\diff{y}\,
    \mathds{1}_{Z_{6}}(x,y)\,|x|^{-\beta}\,
    (\psi_{\infty}k_{\zeta,\sigma})(y^{-1}x)\,
    |y|^{-\alpha}\,f(y)
    \Big|^{q}\notag\\[5pt]
    &\lesssim\,
    \sum_{k\in\mathbb{N}}\,
    \int_{\mathbb{X}}\diff{x}\,
    \Big|\int_{\mathbb{X}}\diff{y}\,
    |y^{-1}x|^{-\alpha-\min\lbrace{0,\beta}\rbrace}\,
    (\psi_{2^{k-1},2^{k+2}}k_{\zeta,\sigma})(y^{-1}x)\,(\Psi_{k}f)(y)
    \Big|^{q}\notag\\[5pt]
    &\lesssim\,
    \sum_{k\in\mathbb{N}}\,
    \|\Psi_{k}f*(|\cdot|^{-\alpha-\min\lbrace{0,\beta}\rbrace}
    \psi_{2^{k-1},2^{k+2}}k_{\zeta,\sigma})\|_{L^{q}(\mathbb{X})}^{q}.
\end{align*}
Similar to those used in the previous proofs, we deduce from Lemma \ref{lemma omega1} and Lemma \ref{lemma omega2} that 
\begin{align}
    \|T_{6}^{\infty}(\zeta,\sigma,\alpha,\beta)f\|_{L^{q}(\mathbb{X})}^{q}\,
    \lesssim\,
    \sum_{k\in\mathbb{N}}\,
    \|\Psi_{k}f\|_{L^{p}(\mathbb{X})}^{q}\,
    \asymp\,
    \|f\|_{L^{p}(\mathbb{X})}^{q},
    \label{T6estim1}
\end{align}
in the simpler cases where either $\frac{1}{p}>\frac{1}{2}-\frac{\zeta}{2|\rho|}$ and $\frac{1}{q}<\frac{1}{2}-\frac{\zeta}{2|\rho|}$ or $\zeta>|\rho|$. The same result holds for $q=+\infty$ with similar arguments. Moreover, note that 
\begin{align*}
    \Phi\Big(\zeta, \frac{2|\rho|}{|\rho|+\zeta},-\alpha-\min\lbrace{0,\beta}\rbrace,2^{k-1},2^{k+2}\Big)\,
    =\,
    \mathrm{O}(1),
\end{align*}
provided that $\alpha$ satisfies \eqref{T4Res} and $\alpha+\beta$ satisfies \eqref{T1ResOff}. With these conditions, the estimate \eqref{T6estim1} remains valid in the critical and off-diagonal cases where $\frac{1}{p}=\frac{1}{2}-\frac{\zeta}{2|\rho|}$ or $\frac{1}{q}=\frac{1}{2}+\frac{\zeta}{2|\rho|}$ with $0\le\zeta\le|\rho|$, and $p\neq{q}$, according to Lemma \ref{lemma omega2}.\textit{(i)}. 

We are left with the critical and diagonal cases where $\frac{1}{p}=\frac{1}{q}=\frac{1}{2}\pm\frac{\zeta}{2|\rho|}$ with $0\le\zeta\le|\rho|$. We will establish the $L^{p}$-boundedness of the operator $T_{6}(\zeta,\sigma,\alpha,\beta)$ based only on conditions \eqref{T4Res} and \eqref{T1ResOn}. This follows directly from the above arguments if $\beta\le0$. However, if we extend the same arguments to the case of positive $\beta$, an additional restriction involving $\alpha$ is required. Thus, we need a more refined analysis to properly account for the contribution of $\beta$ in this case.

Assume now $\beta>0$ and $\frac{1}{p}=\frac{1}{2}\pm\frac{\zeta}{2|\rho|}$. We write
\begin{align}
    \|T_{6}^{\infty}(\zeta,\sigma,\alpha,\beta)f\|_{L^{p}(\mathbb{X})}^{p}\,
    \lesssim\,
    \sum_{k\in\mathbb{N}}\,
    \int_{\mathbb{X}}\diff{x}\,|x|^{-\beta{p}}
    \Big|
        \Psi_{k}f*(|\cdot|^{-\alpha}\,
        \psi_{2^{k-1},2^{k+2}}k_{\zeta,\sigma})(x)
    \Big|^{p}.
    \label{T6estim2}
\end{align}
For any bi-$K$-invariant function $\tilde{B}$ in $L_{\textrm{loc}}^{1}(\mathbb{X})$, it follows from the Kunze-Stein phenomenon \eqref{KSbK1} that 
\begin{align}
    \|\Psi_{k}f*\tilde{B}\|_{L^{p}(\mathbb{X})}\,
    \lesssim\,
    \|\Psi_{k}f\|_{L^{p}(\mathbb{X})}\,
    \|\tilde{B}\varphi_{i(\frac{2}{p}-1)\rho}\|_{L^{1}(\mathbb{X})},
    \label{T6estim3}
\end{align}
for every $k\in\mathbb{N}$. 

In the limiting case where $\zeta=0$ and $p=2$, using the dyadic decomposition \eqref{dyadicN}, the Hölder inequality, and the estimate \eqref{phi0} again, we have
\begin{align*}
    \|\tilde{B}\varphi_{i(\frac{2}{p}-1)\rho}\|_{L^{1}(\mathbb{X})}\,
    &=\,
    \sum_{j\in\mathbb{N}}\,\int_{\mathbb{X}}\diff{x}\,
    \Psi_{j}(x)\,\tilde{B}(x)\,\varphi_{0}(x)\\[5pt]
    &\le\,
    \sum_{j\in\mathbb{N}}\,
    \|\Psi_{j}^{\frac{1}{2}}\tilde{B}\|_{L^{2}(\mathbb{X})}\,
    \|\Psi_{j}^{\frac{1}{2}}\varphi_{0}\|_{L^{2}(\mathbb{X})}\\[5pt]
    &\lesssim\,
    \sum_{j\in\mathbb{N}}\,2^{\frac{j\nu}{2}}
    \|\Psi_{j}^{\frac{1}{2}}\tilde{B}\|_{L^{2}(\mathbb{X})}\,
    \asymp\,
    \|\langle{\,\cdot\,}\rangle^{\frac{\nu}{2}}\tilde{B}\|_{L^{2}(\mathbb{X})},
\end{align*}
where $\langle{x}\rangle=\sqrt{1+|x|^{2}}\asymp|x|$ for $|x|$ large. Together with \eqref{T6estim3} and the duality, we obtain the following weighted Kunze-Stein phenomenon:
\begin{align}
    \|\langle{\,\cdot\,}\rangle^{-\frac{\nu}{2}}(\Psi_{k}f*\tilde{B})\|_{L^{2}(\mathbb{X})}\,
    \lesssim\,
    \|\Psi_{k}f\|_{L^{2}(\mathbb{X})}\,
    \|\tilde{B}\|_{L^{2}(\mathbb{X})}.
    \label{T6estim3L2}
\end{align}
On the one hand, when $\beta\ge\frac{\nu}{2}$, we deduce from \eqref{T6estim2}, \eqref{T6estim3L2}, and Lemma \ref{lemma omega1} with $\tilde{B}(x)=|x|^{-\alpha}(\psi_{2^{k-1},2^{k+2}}k_{\zeta,\sigma})(x)$ that 
\begin{align*}
    \|T_{6}^{\infty}(0,\sigma,\alpha,\beta)f\|_{L^{2}(\mathbb{X})}^{2}\,
    &\lesssim\,
    \sum_{k\in\mathbb{N}}\,
    \|\Psi_{k}f\|_{L^{2}(\mathbb{X})}^{2}\,
    \||\cdot|^{-\alpha}\,\psi_{2^{k-1},2^{k+2}}k_{\zeta,\sigma}\|_{L^{2}(\mathbb{X})}^{2}\,
    \lesssim\,
    \|f\|_{L^{2}(\mathbb{X})}^{2},
\end{align*}
provided that  $\alpha$ satisfies \eqref{T4Res}. On the other hand, if $\beta<\frac{\nu}{2}$, we write 
\begin{align*}
    |x|^{-\beta}|y|^{-\alpha}\,
    =\,
    |x|^{-\frac{\nu}{2}}|y^{-1}x|^{\frac{\nu}{2}-\sigma}\,
    (|y^{-1}x|^{\sigma-\frac{\nu}{2}}|x|^{\frac{\nu}{2}}|x|^{-\beta}|y|^{-\alpha}).
\end{align*}
Recall that $|y|\ge2|x|\ge1$ and $|y^{-1}x|\asymp|y|$ for all $(x,y)\in{Z_{6}}$, then $|y^{-1}x|^{\sigma-\frac{\nu}{2}}|x|^{\frac{\nu}{2}}|x|^{-\beta}|y|^{-\alpha}\lesssim1$, provided that $\alpha+\beta\ge\sigma$. Therefore, the inequality \eqref{T6estim2} can read
\begin{align*}
    \|T_{6}^{\infty}(0,\sigma,\alpha,\beta)f\|_{L^{2}(\mathbb{X})}^{2}\,
    \lesssim\,
    \sum_{k\in\mathbb{N}}\,
    \int_{\mathbb{X}}\diff{x}\,|x|^{-\nu}
    \Big|
        \Psi_{k}f*(|\cdot|^{\frac{\nu}{2}-\sigma}\,
        \psi_{2^{k-1},2^{k+2}}k_{\zeta,\sigma})(x)
    \Big|^{2}.
\end{align*}
By using the weighted Kunze-Stein phenomenon \eqref{T6estim3L2} again, we conclude in the same way that the operator $T_{6}^{\infty}(0,\sigma,\alpha,\beta)$ is $L^{2}$-bounded for $\beta<\frac{\nu}{2}$ as well, provided that $\alpha+\beta\ge\sigma$, that is the condition \eqref{T1ResOn}.

In the remaining case where $\frac{1}{p}=\frac{1}{2}\pm\frac{\zeta}{2|\rho|}$ with $0<\zeta\le|\rho|$, we establish the $L^p$-boundedness of $T_{6}^{\infty}(\zeta,\sigma,\alpha,\beta)$ in the same spirit. Assume that $\frac{1}{p}=\frac{1}{2}-\frac{\zeta}{2|\rho|}$. Similar to the previous $L^2$-case, using the dyadic decomposition \eqref{dyadicN} together with the estimates \eqref{phirho} and \eqref{Lemma1theta2}, we deduce that, for any bi-$K$-invariant function $\tilde{B}\in{L_{\textrm{loc}}^{1}(\mathbb{X})}$,
\begin{align*}
    \|\tilde{B}\varphi_{i(\frac{2}{p}-1)\rho}\|_{L^{1}(\mathbb{X})}\,
    \lesssim\,
    \|\langle{\,\cdot\,}\rangle^{\frac{\ell+1}{2p}}\tilde{B}\|_{L^{p'}(\mathbb{X})},
\end{align*}
since $p>2$ as assumed, and further that
\begin{align}
    \|\langle{\,\cdot\,}\rangle^{-\frac{\ell+1}{2p}}(\Psi_{k}f*\tilde{B})\|_{L^{p}(\mathbb{X})}\,
    \lesssim\,
    \|\Psi_{k}f\|_{L^{p}(\mathbb{X})}\,
    \|\tilde{B}\|_{L^{p'}(\mathbb{X})},
    \label{T6estim3Lp}
\end{align}
by duality and \eqref{T6estim3}. Therefore, in the case where $\beta\ge\frac{\ell+1}{2p}$, we deduce from \eqref{T6estim2} and \eqref{T6estim3Lp} with $\tilde{B}(x)=|x|^{-\alpha}(\psi_{2^{k-1},2^{k+2}}k_{\zeta,\sigma})(x)$ that,
\begin{align*}
    \|T_{6}^{\infty}(\zeta,\sigma,\alpha,\beta)f\|_{L^{p}(\mathbb{X})}^{p}\,
    &\lesssim\,
    \sum_{k\in\mathbb{N}}\,
    \|\Psi_{k}f\|_{L^{p}(\mathbb{X})}^{p}\,
    \||\cdot|^{-\alpha}\,\psi_{2^{k-1},2^{k+2}}k_{\zeta,\sigma}\|_{L^{p'}(\mathbb{X})}^{p}\\[5pt]
    &\lesssim\,
    \sum_{k\in\mathbb{N}}\,
    \Phi(\zeta,p',-\alpha,2^{k-1},2^{k+2})\,
    \|\Psi_{k}f\|_{L^{2}(\mathbb{X})}^{2}\,
    \lesssim\,
    \|f\|_{L^{2}(\mathbb{X})}^{2},
\end{align*}
provided that  $\alpha$ satisfies \eqref{T4Res}, according to Lemma \ref{lemma omega1}. If $\beta<\frac{\ell+1}{2p}$, we write 
\begin{align*}
    |x|^{-\beta}|y|^{-\alpha}\,
    =\,
    |x|^{-\frac{\ell+1}{2p}}|y^{-1}x|^{\frac{\ell+1}{2p}-\frac{\sigma}{2}}\,
    (|y^{-1}x|^{\frac{\sigma}{2}-\frac{\ell+1}{2p}}
    |x|^{\frac{\ell+1}{2p}}|x|^{-\beta}|y|^{-\alpha})
\end{align*}
where $|y^{-1}x|^{\frac{\sigma}{2}-\frac{\ell+1}{2p}}|x|^{\frac{\ell+1}{2p}}|x|^{-\beta}|y|^{-\alpha}\lesssim1$ for $(x,y)\in{Z_{6}}$, $\beta<\frac{\ell+1}{2p}$ and $\alpha+\beta\ge\frac{\sigma}{2}$. Hence, we obtain, by using \eqref{T6estim3Lp} and Lemma \ref{lemma omega1} again,
\begin{align*}
    \|T_{6}^{\infty}(\zeta,\sigma,\alpha,\beta)f\|_{L^{p}(\mathbb{X})}^{p}\,
    &\lesssim\,
    \sum_{k\in\mathbb{N}}\,
    \int_{\mathbb{X}}\diff{x}\,|x|^{-\frac{\ell+1}{2}}
    \Big|
        \Psi_{k}f*(|\cdot|^{\frac{\ell+1}{2p}-\frac{\sigma}{2}}\,
        \psi_{2^{k-1},2^{k+2}}k_{\zeta,\sigma})(x)
    \Big|^{p}\\[5pt]
    &\lesssim\,
    \sum_{k\in\mathbb{N}}\,
    \|\Psi_{k}f\|_{L^{p}(\mathbb{X})}^{p}\,
    \Big\lbrace
        \int_{2^{k-1}}^{2^{k+2}}\diff{r}\,
        r^{p'(\frac{\ell+1}{2p}-\frac{\sigma}{2}+\frac{\sigma}{2}-\frac{\ell+1}{2p})-1}
    \Big\rbrace^{\frac{p}{p'}}\,
    \lesssim\,\|f\|_{L^{p}(\mathbb{X})}^{p},
\end{align*}
since $\frac{1}{p}=\frac{1}{2}-\frac{\zeta}{2|\rho|}$ by assumption. Similar results hold for $\frac{1}{p}=\frac{1}{2}+\frac{\zeta}{2|\rho|}$ by duality.

We conclude that the operator $T_{6}^{\infty}(\zeta,\sigma,\alpha,\beta)$ is bounded from $L^p(\mathbb{X})$ to $L^q(\mathbb{X})$ for all $\alpha\in\mathbb{R}$ and $\beta\in\mathbb{R}$ in the large shift case where $\zeta>|\rho|$ and the non-critical case where $\frac{1}{p}>\frac{1}{2}-\frac{\zeta}{2|\rho|}$ and $\frac{1}{q}<\frac{1}{2}+\frac{\zeta}{2|\rho|}$. In the critical cases where $\frac{1}{p}=\frac{1}{2}-\frac{\zeta}{2|\rho|}$ or $\frac{1}{q}=\frac{1}{2}+\frac{\zeta}{2|\rho|}$ with $0\le\zeta\le|\rho|$, the operator $T_{6}^{\infty}(\zeta,\sigma,\alpha,\beta)$ remains $L^p$-$L^q$-bounded, provided that $\alpha$ satisfies \eqref{T4Res}, and either $\alpha+\beta$ satisfies \eqref{T1ResOff} if $p\neq{q}$ or \eqref{T1ResOn} if $p=q$. We omit the proof for the operator $T_{7}^{\infty}(\zeta,\sigma,\alpha,\beta)$, which can be handled in the same way by symmetry.
\end{proof}

Combining Proposition \ref{T1inf prop}, Proposition \ref{T4T5 prop}, and Proposition \ref{T6T7 prop}, we obtain the following corollary for the $L^{p}$-$L^{q}$-boundedness of the operator $T^{\infty}(\zeta,\sigma,\alpha,\beta)$, which is the sufficiency part of Theorem \ref{mainthm1}.\textit{(2)}. 
\begin{corollary}
Suppose that $\zeta\ge0$, $\sigma>0$ ($0<\sigma<\nu$ if $\zeta=0$), $\alpha\in\mathbb{R}$, $\beta\in\mathbb{R}$, and $\alpha+\beta\ge0$. For all $1\le{p}\le{q}\le+\infty$ such that $\frac{1}{p'}>\frac{\alpha}{n}$ when $\alpha>0$ and $\frac{1}{q}>\frac{\beta}{n}$ when $\beta>0$, the operator $T^{\infty}(\zeta,\sigma,\alpha,\beta)$ is bounded from $L^{p}(\mathbb{X})$ to $L^{q}(\mathbb{X})$ if one of the following conditions is met:
\begin{enumerate}[leftmargin=*,parsep=5pt,label=(\roman*)]
    \item if $\zeta>|\rho|$;
        
    \item if $0\le\zeta\le|\rho|$, then $\frac{1}{p}>\frac{1}{2}-\frac{\zeta}{2|\rho|}$ and $\frac{1}{q}<\frac{1}{2}+\frac{\zeta}{2|\rho|}$;
        
    \item if $0<\zeta\le|\rho|$, $\frac{1}{q}<\frac{1}{2}-\frac{\zeta}{2|\rho|}=\frac{1}{p}$ or $\frac{1}{q}=\frac{1}{2}+\frac{\zeta}{2|\rho|}<\frac{1}{p}$, then $\min\lbrace{\alpha,\beta}\rbrace>\frac{\sigma}{2}-\frac{\ell+1}{2}(\frac{1}{2}-\frac{\zeta}{2|\rho|})$;

    \item if $0<\zeta\le|\rho|$ and $\frac{1}{p}=\frac{1}{q}=\frac{1}{2}\pm\frac{\zeta}{2|\rho|}$, then $\min\lbrace{\alpha,\beta}\rbrace>\frac{\sigma}{2}-\frac{\ell+1}{2}(\frac{1}{2}-\frac{\zeta}{2|\rho|})$ and $\alpha+\beta\ge\frac{\sigma}{2}$;
    
    \item if $\zeta=0$, $p=2<q$ or $p<2=q$, then $\min\lbrace{\alpha,\beta}\rbrace>\sigma-\frac{\nu}{2}$;
    
    \item if $\zeta=0$ and $p=q=2$, then $\min\lbrace{\alpha,\beta}\rbrace>\sigma-\frac{\nu}{2}$ and $\alpha+\beta\ge\sigma$.
\end{enumerate}
\end{corollary}

\subsection{Necessity}
To complete the proof of Theorem \ref{mainthm1}, it remains to show that all the relevant conditions are necessary. We summarize them in the following proposition.
\begin{proposition}
Suppose that $\zeta\ge0$, $\sigma>0$ ($0<\sigma<\nu$ if $\zeta=0$), $\alpha\in\mathbb{R}$, and $\beta\in\mathbb{R}$. Then we have the following.
\begin{enumerate}[leftmargin=*,parsep=5pt]
    \item The operator $T^{0}(\zeta,\sigma,\alpha,\beta)$ is not bounded from $L^{p}(\mathbb{X})$ to $L^{q}(\mathbb{X})$ for any $1\le{p}\le{q}\le+\infty$ if one of the following conditions is met:
    \begin{enumerate}[parsep=5pt,label=(\roman*)]
        \item $\alpha+\beta<0$;
        \item $\alpha+\beta>\sigma$;
        \item $\frac{1}{p}-\frac{1}{q}>\frac{\sigma-\alpha-\beta}{n}$.
    \end{enumerate}

    \item The operator $T^{\infty}(\zeta,\sigma,\alpha,\beta)$ is not bounded from $L^{p}(\mathbb{X})$ to $L^{q}(\mathbb{X})$ for any $1\le{p}\le{q}\le+\infty$ if one of the following conditions is met:
    \begin{enumerate}[parsep=5pt,label=(\roman*)]
        \item $\frac{1}{p}<\frac{1}{2}-\frac{\zeta}{2|\rho|}$;
        
        \item $\frac{1}{q}>\frac{1}{2}+\frac{\zeta}{2|\rho|}$;
        
        \item $\min\lbrace{\alpha,\beta}\rbrace\le\frac{\sigma}{2}-(\frac{1}{2}-\frac{\zeta}{2|\rho|})\frac{\ell+1}{2}$ when $\frac{1}{p}=\frac{1}{2}-\frac{\zeta}{2|\rho|}$ or $\frac{1}{q}=\frac{1}{2}+\frac{\zeta}{2|\rho|}$ with $0<\zeta\le|\rho|$;
        
        \item $\min\lbrace{\alpha,\beta}\rbrace\le\sigma-\frac{\nu}{2}$ when $\frac{1}{p}=\frac{1}{2}$ or $\frac{1}{q}=\frac{1}{2}$ in the limiting case where $\zeta=0$;
        
        \item $\alpha+\beta<\frac{\sigma}{2}$ when $0<\zeta\le|\rho|$ and $\frac{1}{p}=\frac{1}{q}=\frac{1}{2}\pm\frac{\zeta}{2|\rho|}$;

        \item $\alpha+\beta<\sigma$ when $\zeta=0$ and $\frac{1}{p}=\frac{1}{q}=\frac{1}{2}$.
    \end{enumerate}

    \item Neither operator $T^{0}(\zeta,\sigma,\alpha,\beta)$ nor $T^{\infty}(\zeta,\sigma,\alpha,\beta)$ is bounded from $L^{p}(\mathbb{X})$ to $L^{q}(\mathbb{X})$ if one of the following conditions is met:
    \begin{enumerate}[parsep=5pt,label=(\roman*)]
        \item $p>q$;
        \item $\frac{1}{p'}\le\frac{\alpha}{n}$ when $\alpha>0$;
        \item $\frac{1}{q}\le\frac{\beta}{n}$ when $\beta>0$.
    \end{enumerate}
\end{enumerate}
Then the operator $T(\zeta,\sigma,\alpha,\beta)$ is not $L^{p}$-$L^{q}$-bounded, i.e. the Stein-Weiss inequality \eqref{SW} fails, if one of the above conditions is satisfied. 
\end{proposition}

In some special cases we can relate $T(\zeta,\sigma,\alpha,\beta)$ and $T(\zeta,\sigma,0,0)=R(\zeta,\sigma)$, then eliminate some conditions using known results. For most situations, however, the sharp Stein-Weiss inequality cannot be fully described in this way. For the convenience of the reader, we will give very explicit counterexamples for each case in the following proof. The main idea is to choose functions $f\in{L^{p}(\mathbb{X})}$ and $g\in{L^{q'}(\mathbb{X})}$ such that the double integral 
\begin{align*}
    I\,
    =\,
    \int_{\mathbb{X}}\diff{x}\,
    \int_{\mathbb{X}}\diff{y}\,
    |x|^{-\beta}\,k_{\zeta,\sigma}(y^{-1}x)\,|y|^{-\alpha}\,
    f(y)\,g(x)
\end{align*}
diverges in certain regions, or to find a function $f\in{L^{p}(\mathbb{X})}$ such that $T(\zeta,\sigma,\alpha,\beta)f$ does not belong to $L^{q}(\mathbb{X})$. For this purpose, let us denote by
\begin{align} 
    \mathfrak{a}_{\rho}^{+}\,
    =\,
    \lbrace{H\in\overline{\mathfrak{a}^{+}}\,|\,
    \langle{\rho,H}\rangle\asymp|H|}\rbrace
    \label{arho}
\end{align}    
a subset of $\overline{\mathfrak{a}^{+}}$ consisting of vectors around the $\rho$-axis. It is known that $\delta(H)\asymp{e^{2\langle{\rho,H}\rangle}}$ in any proper subcone of $\overline{\mathfrak{a}^{+}}$, see, for instance, \cite[(1.2)]{Ank92}.

\begin{proof}
\textit{(1)} Throughout this proof, let $\varepsilon$ be any small positive constant. We define two positive bi-$K$-invariant functions $f_{1}$ and $g_{1}$ on $\mathbb{X}$ such that
\begin{align*}
    f_{1}(\exp{x}^{+})\,
    =\,
    g_{1}(\exp{x}^{+})\,
    =\,
    \psi_{\infty}(x^{+})\,|x^{+}|^{-\ell-\varepsilon}\,\delta(x^{+})^{-1}
\end{align*}
where $\delta(x^{+})$ is the density function defined in \eqref{density}. Note that $f_{1}\in{L}^{p}(\mathbb{X})$ and $g_{1}\in{L}^{q'}(\mathbb{X})$ for all $1\le{p,\,q}\le+\infty$. Consider the double integral $I$ over $Z_{1}$ where $|x|\asymp|y|$ and $|y^{-1}x|<1$. According to the kernel estimate \eqref{estim k0}, we know that $(\psi_{0}k_{\zeta,\sigma})(y^{-1}x)\ge1$ for all $\sigma>0$ and $(x,y)\in{Z_{1}}$. Therefore, 
\begin{align*}
    I
    &\gtrsim\,
    \int\int_{Z_{1}}\diff{x}\diff{y}\,
    |x|^{-\beta}\,|y|^{-\alpha}\,f_{1}(y)\,g_{1}(x)\\[5pt]
    &\asymp\,
    \int_{\mathfrak{a}^{+}}\diff{x^{+}}\,\psi_{\infty}(x^{+})\,
    |x^{+}|^{-\beta-\ell-\varepsilon}\,
    \int_{|y^{+}|\,\asymp\,|x^{+}|}\diff{y^{+}}\,
    |y^{+}|^{-\alpha-\ell-\varepsilon}\\[5pt]
    &\asymp\,
    \int_{1}^{+\infty}\diff{r}\,r^{-\alpha-\beta-2\varepsilon-1}.
\end{align*}
When $\alpha+\beta<0$, there exists $\varepsilon>0$ such that the last integral diverges. Hence, the operator $T(\zeta,\sigma,\alpha,\beta)$ is not bounded from $L^{p}(\mathbb{X})$ to $L^{q}(\mathbb{X})$ for any $1\le{p,\,q}\le+\infty$ if $\alpha+\beta<0$.

For cases \textit{(ii)} and \textit{(iii)}, we consider two positive bi-$K$-invariant functions $f_{2}$ and $g_{2}$ on $\mathbb{X}$ such that 
\begin{align*}
    f_{2}(x)\,=\,\psi_{0}(x)\,|x|^{-\frac{n}{p}+\varepsilon}
    \qquad\textnormal{and}\qquad
    g_{2}(x)\,=\,\psi_{0}(x)\,|x|^{-\frac{n}{q'}+\varepsilon},
\end{align*}
Then, we know that $f\in{L}^{p}(\mathbb{X})$ and $g\in{L}^{q'}(\mathbb{X})$ for all $1\le{p,\,q}\le+\infty$. We stay in $Z_{1}$ where $|x|\asymp|y|$ and $|y^{-1}x|\le1$. Consider the case where $0<\sigma<n$, we know from \eqref{estim k0} that 
\begin{align*}
    (\psi_{0}k_{\zeta,\sigma})(y^{-1}x)\,
    \asymp\,|y^{-1}x|^{\sigma-n}\,
    \gtrsim\,|x|^{\sigma-n}
    \qquad\forall\,(x,y)\in{Z_{1}}.
\end{align*}
We deduce that 
\begin{align}
    I
    &\gtrsim\,
    \int\int_{Z_{1}}\diff{x}\diff{y}\,
    |x|^{\sigma-n-\beta}\,|y|^{-\alpha}\,f_{2}(y)\,g_{2}(x)\notag\\[5pt]
    &\gtrsim\,
    \int_{\mathfrak{a}_{\rho}^{+}}\diff{x^{+}}\,
    \psi_{0}(x^{+})\,|x^{+}|^{\sigma-n-\beta-\frac{n}{q'}+\varepsilon+n-\ell}\,
    \int_{|y^{+}|\,\asymp\,|x^{+}|}\diff{y^{+}}\,|y|^{-\alpha-\frac{n}{p}+\varepsilon+n-\ell}\notag\\[5pt]
    &\asymp\,
    \int_{0}^{1}\diff{r}\,r^{\sigma-\alpha-\beta-\frac{n}{p}+\frac{n}{q}+2\varepsilon-1}.
    \label{conditionpq}
\end{align}
On the one hand, the last integral cannot be finite for all $\varepsilon>0$ if $\sigma-\alpha-\beta<n(\frac{1}{p}-\frac{1}{q})$. On the other hand, by taking $p=q$, the last integral diverges as well if $\alpha+\beta>\sigma$.

\textit{(2)} The first two cases correspond to cases where the convolution kernel does not provide enough decay to compensate for the exponential growth of the density, while the remaining four cases clarify the conditions on the weights when the exponential factors no longer play a role.

Notice that the condition $1\le{p}\le{q}\le+\infty$ covers $\frac{1}{p}\ge\frac{1}{2}-\frac{\zeta}{2|\rho|}$ and $\frac{1}{q}\le\frac{1}{2}+\frac{\zeta}{2|\rho|}$ when $\zeta>|\rho|$. Assume that $0\le\zeta\le|\rho|$. Let $f_{3}$ be a bi-$K$-invariant cut-off function on $\mathbb{X}$ such that $\supp{f_{3}}$ is contained in $\lbrace{x\in\mathbb{X}\,|\,0\le|x|\le1}\rbrace$ and $f_{3}(y)=1$ when $\frac{1}{4}\le|y|\le\frac{1}{2}$. It is clear that $f_{3}\in{L^p(\mathbb{X})}$ for all $1\le{p}\le+\infty$. According to the kernel estimate \eqref{estim kinf} and the local Harnack inequality (see the proof of Proposition \ref{T4T5 prop}), we have
\begin{align*}
    (\psi_{\infty}k_{\zeta,\sigma})(y^{-1}x)\,
    \asymp\,
    \psi_{\infty}(x)\,|x|^{m_{\zeta}}\,e^{-\zeta|x|}\,\varphi_{0}(x)\,
    \qquad\forall\,y\in\supp{f_{3}},
\end{align*}
where the exponent $m_{\zeta}$ is defined as \eqref{mzeta}. Therefore,
\begin{align}
    \|T(\zeta,\sigma,\alpha,\beta)f_{3}\|_{L^{q}(\mathbb{X})}^{q}\,
    &\gtrsim\,
    \int_{\mathbb{X}}\diff{x}\,
    \psi_{\infty}^{q}(x)\,|x|^{(-\beta+m_{\zeta})q}\,
    e^{-q\zeta|x|}\,\varphi_{0}^{q}(x)\,
    \Big|\int_{\frac{1}{4}\le|y|\le\frac{1}{2}}\diff{y}\,
    |y|^{-\alpha}\Big|^{q}\notag\\[4pt]
    &\gtrsim\,
    \int_{\mathfrak{a}_{\rho}^{+}}\diff{x^{+}}\,
    \psi_{\infty}(x)\,|x|^{(-\beta+m_{\zeta}+\frac{\nu-\ell}{2})q}\,
    e^{(2-q)\langle\rho,x^{+}\rangle}\,e^{-q\zeta|x|}.
    \label{Counterexample1}
\end{align}
Observe that \eqref{Counterexample1} diverges for $q<2$ in the limiting case where $\zeta=0$. When $0<\zeta\le|\rho|$, by writing $x^{+}\in\mathfrak{a}_{\rho}^{+}$ in polar coordinates and using \eqref{Lemma1theta2}, we deduce that, for $\ell\ge2$,
\begin{align}
    \|T(\zeta,\sigma,\alpha,\beta)f_{3}\|_{L^{q}(\mathbb{X})}^{q}\,
    &\gtrsim\,
    \int_{1}^{+\infty}\diff{r}\,r^{\ell-1-\frac{\ell-1}{2}}\,
    r^{(-\beta+m_{\zeta}+\frac{\nu-\ell}{2})q}\,
    e^{(2-q)|\rho|r-q\zeta{r}}\notag\\[5pt]
    &=\,
    \int_{1}^{+\infty}\diff{r}\,
    r^{(-\beta+\frac{\sigma-\ell-1}{2}+\frac{\ell+1}{2q})q-1}\,
    e^{2|\rho|r(1-q(\frac{1}{2}+\frac{\zeta}{2|\rho|}))},
    \label{Counterexample2}
\end{align}
which is not finite for $\frac{1}{q}>\frac{1}{2}+\frac{\zeta}{2|\rho|}$. We skip the simpler argument for $\ell=1$ and conclude that the operator $T(\zeta,\sigma,\alpha,\beta)$ cannot be $L^p$-$L^q$-bounded for $\frac{1}{q}>\frac{1}{2}+\frac{\zeta}{2|\rho|}$, and by duality, for $\frac{1}{p}<\frac{1}{2}-\frac{\zeta}{2|\rho|}$ either.

In the critical cases, we know from \eqref{Counterexample1} that
\begin{align*}
    \|T(\zeta,\sigma,\alpha,\beta)f_{3}\|_{L^{q}(\mathbb{X})}^{q}\,
    &\gtrsim\,
    \int_{\mathfrak{a}_{\rho}^{+}}\diff{x^{+}}\,
    \psi_{\infty}(x^{+})\,|x^{+}|^{2(-\beta+\sigma-\nu+\frac{\nu-\ell}{2})}\\[5pt]
    &\asymp\,
    \int_{1}^{+\infty}\diff{r}\,r^{-2\beta+2\sigma-\nu-1},
\end{align*}
when $\zeta=0$ and $q=2$. The last integral is finite if and only if $\beta>\sigma-\frac{\nu}{2}$. On the other hand, if $0<\zeta\le|\rho|$ and $\frac{1}{q}=\frac{1}{2}+\frac{\zeta}{2|\rho|}$, we deduce from \eqref{Counterexample2} that
\begin{align*}
    \|T(\zeta,\sigma,\alpha,\beta)f_{3}\|_{L^{q}(\mathbb{X})}^{q}\,
    &\gtrsim\,
    \int_{1}^{+\infty}\diff{r}\,
    r^{(-\beta+\frac{\sigma}{2}-\frac{\ell+1}{2}+\frac{\ell+1}{2q})q-1},
\end{align*}
which diverges when $\beta\le\frac{\sigma}{2}-(\frac{1}{2}-\frac{\zeta}{2|\rho|})\frac{\ell+1}{2}$. By duality, the similar conclusion holds for $\frac{1}{p}=\frac{1}{2}-\frac{\zeta}{2|\rho|}$, with constraints on $\alpha$. We have just shown that the operator $T^{\infty}(\zeta,\sigma,\alpha,\beta)$ cannot be $L^{p}$-$L^{q}$-bounded if condition \textit{(2).(iii)} or \textit{(2).(iv)} is satisfied.

For the diagonal cases \textit{(2).(v)} or \textit{(2).(vi)}, we assume for simplicity that $\beta\ge0$. Recall that $|y|\ge2|x|\ge1$ and $|y^{-1}x|\asymp|y|$ for $(x,y)$ in $Z_{6}$, see Figure \ref{Regions}. Then, under these assumptions, we have $|x|^{-\beta}|y|^{-\alpha}\gtrsim|y^{-1}x|^{-\alpha-\beta}$. As in the proof of Proposition \ref{T6T7 prop}, we obtain
\begin{align*}
    \|T_{6}^{\infty}(\zeta,\sigma,\alpha,\beta)f\|_{L^{p}(\mathbb{X})}^{p}\,
    \gtrsim\,
    \sum_{k\in\mathbb{N}}\,
    \|\Psi_{k}f*(|\cdot|^{-\alpha-\beta}\,\psi_{2^{k-1},2^{k+2}}k_{\zeta,\sigma})\|_{L^{p}(\mathbb{X})}^{p}
\end{align*}
for all $f\in{L^p(\mathbb{X})}$, where $\frac{1}{p}=\frac{1}{2}\pm\frac{\zeta}{2|\rho|}$ with $0\le\zeta\le|\rho|$ and the cut-off functions $\Psi_{k}$ are defined as \eqref{dyadicN}. According to the Kunze-Stein phenomenon \eqref{KSbK1}, the operator $T_{6}^{\infty}(\zeta,\sigma,\alpha,\beta)$ can be $L^p$-bounded only if 
\begin{align*}
    \tilde{I}(k)\,
    =\,
    \int_{\mathbb{X}}\diff{x}\,
    |x|^{-\alpha-\beta}\,\psi_{2^{k-1},2^{k+2}}(x)\,
    k_{\zeta,\sigma}(x)\,\varphi_{i(\frac{2}{p}-1)\rho}(x)
\end{align*}
is finite as $k$ tends to $+\infty$. However, according to the kernel estimate \eqref{estim kinf} and the estimate \eqref{phirho} of the spherical function, we have, for $p>2$, 
\begin{align*}
    \tilde{I}(k)\,
    &\asymp\,
    \int_{\mathfrak{a}^{+}}\diff{x^{+}}\,\delta(x^{+})\,
    \psi_{2^{k-1},2^{k+2}}(x^{+})\,|x^{+}|^{-\alpha-\beta+m_{\zeta}}\,
    e^{-\zeta|x^{+}|}\,(\varphi_{0}\varphi_{i(\frac{2}{p}-1)\rho})(\exp{x^{+}})\\[5pt]
    &\gtrsim\,
    \int_{\mathfrak{a}_{\rho}^{+}}\diff{x^{+}}\,
    \psi_{2^{k-1},2^{k+2}}(x^{+})\,|x^{+}|^{-\alpha-\beta+m_{\zeta}+\frac{\nu-\ell}{2}}\,
    e^{-\zeta|x^{+}|}\,e^{(1-\frac{2}{p})\langle{\rho,x^{+}}\rangle},
\end{align*}
which diverges as $k\rightarrow+\infty$ if the condition \textit{(2).(v)} is met. In fact, we know from \eqref{Lemma1theta2} that
\begin{align*}
    \tilde{I}(k)\,
    \gtrsim\,
    \int_{2^{k-1}}^{2^{k+2}}\diff{r}\,
        r^{-\alpha-\beta+\frac{\sigma-\nu-1}{2}+\frac{\nu-\ell}{2}-\frac{\ell-1}{2}+\ell-1}\,
    =\,
    \int_{2^{k-1}}^{2^{k+2}}\diff{r}\,r^{-\alpha-\beta+\frac{\sigma}{2}-1},
\end{align*}
for $\frac{1}{p}=\frac{1}{2}-\frac{\zeta}{2|\rho|}$ with $0<\zeta\le|\rho|$. Similarly, when $\zeta=0$ and $p=2$, by using the estimate \eqref{phi0}, we have straightforwardly that 
\begin{align*}
    \tilde{I}(k)\,
    \gtrsim\,
    \int_{2^{k-1}}^{2^{k+2}}\diff{r}\,
        r^{-\alpha-\beta+\sigma-\nu+\nu-\ell+\ell-1}\,
    =\,
    \int_{2^{k-1}}^{2^{k+2}}\diff{r}\,r^{-\alpha-\beta+\sigma-1},
\end{align*}
which diverges as $k\rightarrow+\infty$ if $\alpha+\beta<\sigma$, or in other words, the condition \textit{(2).(vi)} is met. Here, the exponent $m_{\zeta}$ is defined as \eqref{mzeta}. Similar conclusions apply to the case where $\frac{1}{p}=\frac{1}{2}+\frac{\zeta}{2|\rho|}$ by duality.

\textit{(3)} First, taking respectively $\alpha+\beta=\sigma$ in \eqref{conditionpq} and $\zeta=0$ in \textit{(2).(i)} and \textit{(2).(ii)}, we know that neither the operator $T^{0}(\zeta,\sigma,\alpha,\beta)$ nor $T^{\infty}(\zeta,\sigma,\alpha,\beta)$ can be $L^{p}$-$L^{q}$-bounded if $p>q$.

Let $f_{3}$ be the bi-$K$ invariant cut-off function defined as above. Recall that $f_{3}\in{L^{p}}(\mathbb{X})$ for any $1\le{p}\le+\infty$ and $f_{3}(y)=1$ when $\frac{1}{4}\le|y|\le\frac{1}{2}$. According to the kernel estimate \eqref{estim k0}, we know that $(\psi_{0}k_{\zeta,\sigma})(y^{-1}x)\ge1$ for all $\zeta\ge0$ and $\sigma>0$, if $|x|\le1$ and $\frac{1}{4}\le|y|\le\frac{1}{2}$, i.e., if $(x,y)\in{Z_{3}}$. Therefore, for all $1\le{q}<+\infty$,
\begin{align*}
    \|T_{3}^{0}(\zeta,\sigma,\alpha,\beta)f_{3}\|_{L^{q}(\mathbb{X})}^{q}\,
    &=\,
    \int_{\mathbb{X}}\diff{x}\,
    \Big|\int_{\mathbb{X}}\diff{y}\,
    |x|^{-\beta}\,(\psi_{0}k_{\zeta,\sigma})(y^{-1}x)\,|y|^{-\alpha}\,f_{3}(y)
    \Big|^{q}\\[5pt]
    &\ge\,
    \int_{\mathbb{X}}\diff{x}\,\psi_{0}(x)\,|x|^{-\beta{q}}
    \underbrace{
    \Big|\int_{\frac{1}{4}\,\le\,|y|\,\le\,\frac{1}{2}}\diff{y}\,|y|^{-\alpha}
    \Big|^{q}}_{=\,\const},
\end{align*}
where the inner integral is a constant for all $\alpha\in\mathbb{R}$. We deduce that 
\begin{align*}
    \|T_{3}^{0}(\zeta,\sigma,\alpha,\beta)f_{3}\|_{L^{q}(\mathbb{X})}^{q}\,
    &\gtrsim\,
    \int_{K(\exp\mathfrak{a}_{\rho}^{+})K}\diff{x}\,
    \psi_{0}(x)\,|x|^{-\beta{q}}\\[5pt]
    &\asymp\,
    \int_{\mathfrak{a}_{\rho}^{+}}\diff{x^{+}}\,
    \psi_{0}(x^{+})\,|x^{+}|^{-\beta{q}+n-\ell}\,
    \asymp\,
    \int_{0}^{1}\diff{r}\,r^{-\beta{q}+n-1},
\end{align*}
which is not finite for $\beta\ge\frac{n}{q}$. In the sup norm case, assume that $0<\sigma<n$. Note that $(\psi_{0}k_{\zeta,\sigma})(y^{-1}x)\gtrsim|x|^{-\varepsilon}$ for some $\varepsilon>0$ in the region $Z_3$. Hence, the operator $T_{3}^{0}(\zeta,\sigma,\alpha,\beta)$ is not bounded for any $\beta\ge0$. We conclude that the operator $T^{0}(\zeta,\sigma,\alpha,\beta)$ cannot be $L^{p}$-$L^{q}$-bounded for any $1\le{p}\le{q}\le+\infty$ if $\frac{1}{q}\le\frac{\beta}{n}$ when $\beta>0$. Similarly, by considering the region where $\frac{1}{4}\le|x|\le\frac{1}{2}$ and $|y|\le1$, and by using the duality, we can show that the operator $T^{0}(\zeta,\sigma,\alpha,\beta)$ is not $L^{p}$-$L^{q}$-bounded either for any $\alpha\ge\frac{n}{p'}$ when $\alpha>0$. For the operator $T^{\infty}(\zeta,\sigma,\alpha,\beta)$, we consider a cut-off function $f$ which is compactly supported in an annulus away from the origin. We conclude by applying the same argument instead in the region $Z_4$ where $|x|$ is small but $|y|$ is large.
\end{proof}

\noindent\textbf{Acknowledgement.}
The authors would like to thank Jean-Philippe Anker and Stefano Meda for their valuable discussions and detailed explanations of their classic works. The authors are supported by the FWO Odysseus 1 grant G.0H94.18N: Analysis and Partial Differential Equations and by the Methusalem programme of the Ghent University Special Research Fund (BOF) (Grant number 01M01021). The first and second authors are supported by FWO Senior Research Grant G011522N. The second author is also supported by EPSRC grants EP/R003025/2 and EP/V005529/1. The third author receives funding from the Deutsche Forschungsgemeinschaft (DFG) via SFB-TRR 358/1 2023 — 491392403 (CRC “Integral Structures in Geometry and Representation Theory”).

\printbibliography

\vspace{10pt}

\address{
    \noindent\textsc{Vishvesh Kumar:}
    \href{mailto:kumar.vishvesh@ugent.be}
    {vishvesh.kumar@ugent.be}\\
    \textsc{
    Department of Mathematics: 
    Analysis, Logic and Discrete Mathematics\\
    Ghent University, Belgium}
}\vspace{10pt}

\address{
    \noindent\textsc{Michael Ruzhansky:}
    \href{mailto:michael.ruzhansky@ugent.be}
    {michael.ruzhansky@ugent.be}\\
    \textsc{
    Department of Mathematics: 
    Analysis, Logic and Discrete Mathematics\\
    Ghent University, Belgium\\
    and\\
    School of Mathematical Sciences\\
    Queen Mary University of London, United Kingdom}
}\vspace{10pt}

\address{
    \noindent\textsc{Hong-Wei Zhang:}
    \href{mailto:zhongwei@math.uni-paderborn.de}
    {zhongwei@math.uni-paderborn.de}\\
    \textsc{
    Faculty EIM, Institute of Mathematics \\
    Paderborn University\\
    and\\
    Department of Mathematics: 
    Analysis, Logic and Discrete Mathematics\\
    Ghent University, Belgium}
}

\end{document}

%% file: HLS02.tex
\begin{tikzpicture}[line cap=round,line join=round,>=triangle 45,x=1.0cm,y=1.0cm, scale=0.33]
\clip(-3,-2) rectangle (15,11);
\draw [line width=0.5pt] (0,0)-- (0,10);
\draw [line width=0.5pt] (5.,0.)-- (10,0);
\draw [line width=0.5pt] (10,0.)-- (10,5);
\draw [line width=0.5pt] (0.,10)-- (10,10);

\draw [line width=1pt,color=cyan] (0,0)--(10,10);
\draw [line width=1pt,color=cyan] (5,0)--(10,5);
\draw [line width=1pt,color=cyan] (10,10)--(10,5);
\draw [line width=1pt,color=cyan] (5,0)--(0,0);

\fill[line width=1pt,color=cyan,fill=cyan, opacity=0.1] (0,0) -- (5,0) -- (10,5) -- (10,10) --cycle;

\draw (5,0) node[color=cyan]{$\circ$};
\draw (10,5) node[color=cyan]{$\circ$};

\draw (-0.75,-0.75) node{\small $0$};
\draw (10,-0.75) node{\small $1$};
\draw (-0.75,10) node{\small $1$};

\draw (5,-1.25) node{\small $\frac{\sigma}{n}$};
\draw (12,5) node{\small $1-\frac{\sigma}{n}$};

\draw (8,2) node[rotate=45]{\tiny $\frac{1}{q}\!\!=\!\!\frac{1}{p}\!-\!\frac{\sigma}{n}$};
\end{tikzpicture}

%% file: SW01.tex
\begin{tikzpicture}[line cap=round,line join=round,>=triangle 45,x=1.0cm,y=1.0cm, scale=0.33]
\clip(-3,-2) rectangle (15,11);
\draw [line width=0.5pt] (0,0)-- (0,10);
\draw [line width=0.5pt] (0,0.)-- (10,0);
\draw [line width=0.5pt] (10,10.)-- (10,0);
\draw [line width=0.5pt] (0.,10)-- (10,10);
\draw [line width=0.5pt] (0.,0)-- (1,1);
\draw [line width=0.5pt] (8.,8)-- (10,10);
\draw [line width=0.5pt] (5,0)-- (6,1);
\draw [line width=0.5pt] (10,5)-- (8,3);

\draw [line width=1pt,color=cyan] (1,1)--(8,8);
\draw [line width=1pt,color=cyan] (6,1)--(8,3);
\draw [line width=1pt,color=cyan,dashed] (8,8)--(8,3);
\draw [line width=1pt,color=cyan,dashed] (1,1)--(6,1);

\fill[line width=1pt,color=cyan,fill=cyan, opacity=0.1] (1,1) -- (6,1) -- (8,3) -- (8,8) --cycle;

\draw (6,1) node[color=cyan]{$\circ$};
\draw (8,3) node[color=cyan]{$\circ$};
\draw (1,1) node[color=cyan]{$\circ$};
\draw (8,8) node[color=cyan]{$\circ$};

\draw (-0.75,-0.75) node{\small $0$};
\draw (10,-0.75) node{\small $1$};
\draw (-0.75,10) node{\small $1$};

\draw (5,-1) node{\small $\frac{\sigma-\alpha-\beta}{n}$};
\draw (12.5,5) node{\small $1\!-\!\frac{\sigma-\alpha-\beta}{n}$};
\draw (10.7,1) node{\small $\frac{\beta}{n}$};
\draw (11.2,8) node{\small $1\!-\!\frac{\alpha}{n}$};

\draw (8,2) node[rotate=45]{\tiny $\frac{1}{q}\!\!=\!\!\frac{1}{p}\!-\!\frac{\sigma\!-\!\alpha\!-\!\beta}{n}$};
\end{tikzpicture}

%% file: SW03.tex
\begin{tikzpicture}[line cap=round,line join=round,>=triangle 45,x=1.0cm,y=1.0cm, scale=0.33]
\clip(-3,-2) rectangle (15,11);
\draw [line width=0.5pt] (0,0)-- (0,10);
\draw [line width=1pt,color=cyan] (10,1.)-- (10,10);
\draw [line width=0.5pt] (10,0)-- (10,1);
\draw [line width=0.5pt] (0.,10)-- (10,10);

\draw [line width=1pt,color=cyan] (1,1)--(10,10);
\draw [line width=0.5pt] (1,1)--(0,0);
\draw [line width=1pt,color=cyan,dashed] (1,1)--(10,1);
\draw [line width=0.5pt] (0,0)--(10,0);

\fill[line width=1pt,color=cyan,fill=cyan, opacity=0.1] (1,1) -- (10,1) -- (10,10) --cycle;

\draw (10,1) node[color=cyan]{$\circ$};
\draw (1,1) node[color=cyan]{$\circ$};

\draw (-0.75,-0.75) node{\small $0$};
\draw (10,-0.75) node{\small $1$};
\draw (-0.75,10) node{\small $1$};

\draw (10.7,1) node{\small $\frac{\beta}{n}$};
\draw (1,-1) node{\small $\frac{\beta}{n}$};
\end{tikzpicture}

%% file: HLSinf2.tex
\begin{tikzpicture}[line cap=round,line join=round,>=triangle 45,x=1.0cm,y=1.0cm, scale=0.33]
\clip(-3,-2) rectangle (15,11);
\draw [line width=0.5pt] (0,0)-- (0,10);
\draw [line width=0.5pt] (0.,0.)-- (2.5,0);
\draw [line width=0.5pt] (10,7.5)-- (10,10);
\draw [line width=0.5pt] (0.,10)-- (10,10);

\draw [line width=0.5pt] (0,0)--(2.5,2.5);
\draw [line width=0.5pt] (7.5,7.5)--(10,10);

\draw [line width=1pt,color=red] (2.5,2.5)--(2.5,0);
\draw [line width=1pt,color=red] (7.5,7.5)--(10,7.5);
\draw [line width=1pt,color=cyan] (2.5,0)--(10,0);
\draw [line width=1pt,color=cyan] (10,7.5)--(10,0);

\draw [line width=1pt,color=cyan] (2.5,2.5)--(7.5,7.5);

\fill[line width=0.pt,color=cyan,fill=cyan,fill opacity=0.1] (2.5,0) -- (2.5,2.5) -- (7.5,7.5) -- (10,7.5) -- (10,0) --cycle;

\draw (2.5,0) node[color=red]{$\bullet$};
\draw (10,7.5) node[color=red]{$\bullet$};
\draw (2.5,2.5) node[color=cyan]{$\circ$};
\draw (7.5,7.5) node[color=cyan]{$\circ$};
\draw (2.5,0) node[color=cyan]{$\circ$};
\draw (10,7.5) node[color=cyan]{$\circ$};

\draw (-0.75,-0.75) node{\small $0$};
\draw (10,-0.75) node{\small $1$};
\draw (-0.75,10) node{\small $1$};

\draw (2.5,-1.25) node{\small $\frac12-\frac{\zeta}{2|\rho|}$};
\draw (12,7.5) node{\small $\frac12\!+\!\frac{\zeta}{2|\rho|}$};

\end{tikzpicture}

%% file: SWinf2.tex
\begin{tikzpicture}[line cap=round,line join=round,>=triangle 45,x=1.0cm,y=1.0cm, scale=0.33]
    \clip(-3,-2) rectangle (15,11);
    \draw [line width=0.5pt] (0,0)-- (0,10);
    \draw [line width=0.5pt] (0,0.)-- (10,0);
    \draw [line width=0.5pt] (10,10.)-- (10,0);
    \draw [line width=0.5pt] (0.,10)-- (10,10);
    \draw [line width=0.5pt] (0.,0)-- (3,3);
    \draw [line width=0.5pt] (7,7)-- (10,10);

    \draw [line width=1pt,color=cyan] (3,3)--(7,7);
    
    \draw [line width=1pt,color=cyan,dashed] (3,1)--(8,1);
    \draw [line width=1pt,color=cyan,dashed] (8,1)--(8,7);

    \draw [line width=1pt,color=red] (3,3)--(3,1);
    \draw [line width=1pt,color=red] (7,7)--(8,7);
    
    \fill[line width=1pt,color=cyan,fill=cyan,fill opacity=0.1] (3,1) -- (3,3) -- (7,7) -- (8,7) -- (8,1) -- (3,1) --cycle;
    
    \draw (3,3) node[color=red]{$\bullet$};
    \draw (7,7) node[color=red]{$\bullet$};

    \draw (3,3) node[color=cyan]{$\circ$};
    \draw (3,1) node[color=cyan]{$\circ$};
    \draw (7,7) node[color=cyan]{$\circ$};
    \draw (8,7) node[color=cyan]{$\circ$};
    
    \draw (-0.75,-0.75) node{\small $0$};
    \draw (10,-0.75) node{\small $1$};
    \draw (-0.75,10) node{\small $1$};

    \draw (10.7,1) node{\small $\frac{\beta}{n}$};
    \draw (8,-1) node{\small $1-\frac{\alpha}{n}$};
    \draw (12,7) node{\small $\frac{1}{2}\!+\!\frac{\zeta}{2|\rho|}$};
    \draw (3,-1) node{\small $\frac{1}{2}\!-\!\frac{\zeta}{2|\rho|}$};
    
    \end{tikzpicture}

%% file: SWinf1.tex
\begin{tikzpicture}[line cap=round,line join=round,>=triangle 45,x=1.0cm,y=1.0cm, scale=0.33]
\clip(-3,-2) rectangle (15,11);
\draw [line width=0.5pt] (0,0)-- (0,10);
\draw [line width=0.5pt] (0,0.)-- (10,0);
\draw [line width=0.5pt] (10,0.)-- (10,1);
\draw [line width=0.5pt] (0.,10)-- (10,10);
\draw [line width=0.5pt] (0.,0)-- (1,1);

\draw [line width=1pt,color=cyan] (10,10.)-- (10,1);
\draw [line width=1pt,color=cyan] (1,1)--(10,10);
\draw [line width=1pt,color=cyan,dashed] (1,1)--(10,1);

\fill[line width=1pt,color=cyan,fill=cyan,fill opacity=0.1] (1,1) -- (10,1) -- (10,10) --cycle;

\draw (1,1) node[color=cyan]{$\circ$};
\draw (10,1) node[color=cyan]{$\circ$};

\draw (-0.75,-0.75) node{\small $0$};
\draw (10,-0.75) node{\small $1$};
\draw (-0.75,10) node{\small $1$};

\draw (10.7,1) node{\small $\frac{\beta}{n}$};
\draw (1,-1) node{\small $\frac{\beta}{n}$};
\end{tikzpicture}

%% file: HLSinf3.tex
\begin{tikzpicture}[line cap=round,line join=round,>=triangle 45,x=1.0cm,y=1.0cm, scale=0.33]
\clip(-2,-2) rectangle (11,11);
\draw [line width=0.5pt] (0,0)-- (0,10);
\draw [line width=0.5pt] (0.,0.)-- (5,0);
\draw [line width=0.5pt] (10,5)-- (10,10);
\draw [line width=0.5pt] (0.,10)-- (10,10);

\draw [line width=0.5pt] (0,0)--(10,10);

\draw [line width=1pt,color=cyan] (5.,0.)-- (10,0);
\draw [line width=1pt,color=cyan] (10,0.)-- (10,5);

\fill[line width=1pt,color=cyan,fill=cyan,fill opacity=0.1] (5,0) -- (10,0) -- (10,5) -- (5,5)--cycle;

\draw [line width=1pt,color=cyan] (5,0)--(5,5);
\draw [line width=1pt,color=cyan] (10,5)--(5,5);

\draw (-0.75,-0.75) node{\small $0$};
\draw (10,-0.75) node{\small $1$};
\draw (-0.75,10) node{\small $1$};

\draw (5,5) node[color=cyan]{$\circ$};

\draw (5,-1.25) node{\small $\frac{1}{2}$};
\draw (10.5,5) node{\small $\frac{1}{2}$};

\end{tikzpicture}

%% file: SWinf3.tex
\begin{tikzpicture}[line cap=round,line join=round,>=triangle 45,x=1.0cm,y=1.0cm, scale=0.33]
    \clip(-3,-2) rectangle (15,11);
    \draw [line width=0.5pt] (0,0)-- (0,10);
    \draw [line width=0.5pt] (0,0.)-- (10,0);
    \draw [line width=0.5pt] (10,10.)-- (10,0);
    \draw [line width=0.5pt] (0.,10)-- (10,10);
    \draw [line width=0.5pt] (0.,0)-- (10,10);
    
    \draw [line width=1pt,color=cyan,dashed] (5,2)--(9,2);
    \draw [line width=1pt,color=cyan,dashed] (9,2)--(9,5);

    \draw [line width=1pt,color=cyan] (5,5)--(5,2);
    \draw [line width=1pt,color=cyan] (5,5)--(9,5);
    
    \fill[line width=1pt,color=cyan,fill=cyan,fill opacity=0.1] (5,5) -- (5,2) -- (9,2) -- (9,5) --cycle;
    
    \draw (5,5) node[color=red]{$\bullet$};

    \draw (5,2) node[color=cyan]{$\circ$};
    \draw (9,5) node[color=cyan]{$\circ$};
    \draw (5,5) node[color=cyan]{$\circ$};
    
    \draw (-0.75,-0.75) node{\small $0$};
    \draw (10,-0.75) node{\small $1$};
    \draw (-0.75,10) node{\small $1$};

    \draw (10.7,2) node{\small $\frac{\beta}{n}$};
    \draw (8,-1) node{\small $1\!-\!\frac{\alpha}{n}$};
    \draw (10.7,5) node{\small $\frac{1}{2}$};
    \draw (5,-1) node{\small $\frac{1}{2}$};
    
    \end{tikzpicture}

%% file: MainTHM.tex
\begin{tikzpicture}[line cap=round,line join=round,>=triangle 45,x=1.0cm,y=1.0cm, scale=0.75]
\clip(-1,-1) rectangle (12,12);
\draw [line width=0.5pt] (0,0)-- (0,10);
\draw [line width=0.5pt] (0.,0.)-- (10,0);
\draw [line width=0.5pt] (10,0)-- (10,10);
\draw [line width=0.5pt] (0.,10)-- (10,10);
\draw [line width=0.5pt] (0,0)--(2.5,2.5);
\draw [line width=0.5pt] (7.5,7.5)--(10,10);

\draw [line width=1pt,color=cyan] (2.5,2.5)--(7.5,7.5);
\draw [line width=1.5pt,color=green] (6,1)--(9,4);

\draw [line width=1pt,color=cyan,dashed] (2.5,1)-- (6,1);
\draw [line width=1pt,color=cyan,dashed] (9,4)-- (9,7.5);

\draw [line width=1pt,color=red] (2.5,2.5)--(2.5,1);
\draw [line width=1pt,color=red] (7.5,7.5)--(9,7.5);

\fill[line width=0.pt,color=cyan,fill=cyan,fill opacity=0.1] (2.5,2.5) -- (2.5,1) -- (6,1) -- (9,4) -- (9,7.5) -- (7.5,7.5) --cycle;

\draw (2.5,2.5) node[color=red]{\large $\bullet$};
\draw (7.5,7.5) node[color=red]{\large $\bullet$};
\draw (2.5,2.5) node[color=cyan]{\Large $\circ$};
\draw (7.5,7.5) node[color=cyan]{\Large $\circ$};
\draw (6,1) node[color=cyan]{\Large $\circ$};
\draw (9,4) node[color=cyan]{\Large $\circ$};
\draw (2.5,1) node[color=cyan]{\Large $\circ$};
\draw (9,7.5) node[color=cyan]{\Large $\circ$};

\draw (-0.25,-0.25) node{$O$};
\draw (10,-0.5) node{\large $1$};
\draw (-0.25,10) node{\large $1$};
\draw (11, 7.5) node{\large $\frac{1}{2}+\frac{\zeta}{2|\rho|}$};
\draw (10.25, 4) node{\large $\frac{\beta}{n}$};
\draw (2.5, -0.5) node{\large $\frac{1}{2}-\frac{\zeta}{2|\rho|}$};
\draw (6, -0.5) node{\large $1-\frac{\alpha}{n}$};

\draw (8,2.25) node[rotate=45]{\small $\frac{1}{q}-\frac{1}{p}=\frac{\alpha+\beta-\sigma}{n}$};
\end{tikzpicture}

%% file: Regions.tex
\begin{tikzpicture}[line cap=round,line join=round,>=triangle 45,x=1.0cm,y=1.0cm, scale=0.75]
\clip(-1,-1) rectangle (10,10);
\draw [line width=0.5pt] (0,0)-- (0,10);
\draw [line width=0.5pt] (0.,0.)-- (10,0);

\draw [line width=1.pt,dashed] (0,0)-- (5,10);
\draw [line width=1.pt,dashed] (0,0)-- (10,5);
\draw [line width=1.pt,dashed] (0,3)-- (1.5,3);
\draw [line width=1.pt,dashed] (1.5,3)-- (1.5,10);
\draw [line width=1.pt,dashed] (3,0)-- (3,1.5);
\draw [line width=1.pt,dashed] (3,1.5)-- (10,1.5);

\draw (-.5,-.5) node{$O$};
\draw (3,-.5) node{$1$};
\draw (9.5,-.5) node{$|x|$};
\draw (-.5,3) node{$1$};
\draw (-.5,9.5) node{$|y|$};

\draw (5,5) node{$Z_{1}$};
\draw (0.5,2.5) node{$Z_{2}$};
\draw (2.5,0.5) node{$Z_{3}$};
\draw (0.7,7) node{$Z_{4}$};
\draw (7,0.7) node{$Z_{5}$};
\draw (2.5,7) node{$Z_{6}$};
\draw (7,2.5) node{$Z_{7}$};

\draw (8.5,4.8) node[rotate=30]{$|y|=\frac{|x|}{2}$};
\draw (4.8,8.5) node[rotate=60]{$|y|=2|x|$};

\end{tikzpicture}

%% file: Lemma1.tex
\begin{tikzpicture}[line cap=round,line join=round,>=triangle 45,x=1.0cm,y=1.0cm, scale=0.75]
\clip(-1,-1) rectangle (12,12);
\draw [line width=0.5pt] (0,0)-- (0,10);
\draw [line width=0.5pt] (0.,0.)-- (2.5,0);
\draw [line width=0.5pt] (10,7.5)-- (10,10);
\draw [line width=0.5pt] (0.,10)-- (10,10);
\draw [line width=1pt,color=cyan] (2.5,0.)-- (10,0);
\draw [line width=1pt,color=cyan] (10,0.)-- (10,7.5);

\draw [line width=1pt,dashed] (0,0)--(2.5,2.5);
\draw [line width=1pt,dashed] (7.5,7.5)--(10,10);
\draw [line width=0.5pt,color=cyan,dashed] (5,5)--(10,0);
\draw [line width=1pt,color=cyan] (2.5,2.5)--(2.5,0);
\draw [line width=1pt,color=cyan] (7.5,7.5)--(10,7.5);
\draw [line width=0.5pt,color=cyan,dashed] (5,5)--(5,0);

\draw [line width=1pt,color=red] (2.5,2.5)--(7.5,7.5);

\fill[line width=0.pt,color=cyan,fill=cyan,fill opacity=0.1] (2.5,0) -- (2.5,2.5) -- (5,5) -- (10,0) --cycle;

\draw (2.5,2.5) node[color=red]{\large $\bullet$};
\draw (5,5) node[color=red]{\large $\bullet$};
\draw (7.5,7.5) node[color=red]{\large $\bullet$};
\draw (2.5,2.5) node[color=cyan]{\Large $\circ$};
\draw (5,5) node[color=cyan]{\Large $\circ$};
\draw (7.5,7.5) node[color=cyan]{\Large $\circ$};

\draw (-0.25,-0.25) node{$O$};
\draw (9.75,-0.5) node{\large $\frac{1}{p}$};
\draw (-0.25,10) node{\large $\frac{1}{q}$};
\draw (11, 7.5) node{\large $\frac{1}{2}+\frac{\zeta}{2|\rho|}$};
\draw (10.25, 5) node{\large $\frac{1}{2}$};
\draw (2.5, -0.5) node{\large $\frac{1}{2}-\frac{\zeta}{2|\rho|}$};
\draw (5, -0.5) node{\large $\frac{1}{2}$};

\draw (8.2,2.5) node[rotate=-45]{$\frac{1}{q}=\frac{1}{p'}$};
\draw (4.5,5.5) node[rotate=45]{$\frac{1}{q}=\frac{1}{p}$};
\end{tikzpicture}